\documentclass{article}
\usepackage{hyperref}
\usepackage[autostyle=true, english=american]{csquotes}
\usepackage{amsmath}
\usepackage{amsthm}
\usepackage{amssymb}
\usepackage{mathtools}
\usepackage[edge length=0.7cm, root radius=0.08cm]{dynkin-diagrams}
\usepackage{setspace}
\usepackage{tikz,tikz-cd}
\usepackage{multicol}
\usepackage{makecell}
\usepackage{booktabs}
\usepackage{subcaption}
\usepackage{fancyhdr}
\usepackage[shortlabels]{enumitem}

\usetikzlibrary{positioning}
\usetikzlibrary{calc}

\newcommand*{\map}[5]{\ifblank{#1}{}{#1\colon} \ifblank{#2}{}{#2 \rightarrow #3} \ifblank{#4}{}{\ifblank{#2}{}{,} #4 \mapsto #5}}
\newcommand{\mapdot}{\mathord{\:\cdot\:}}
\newcommand{\Set}[1]{\{#1\}}
\newcommand{\defl}{\coloneqq}
\newcommand{\union}{\cup}
\newcommand{\braidword}[1]{P_{#1}}
\newcommand{\invosym}{*}
\newcommand{\invosymroot}[1]{\invosym_{#1}}
\newcommand{\invo}[1]{#1^{\invosym}}
\newcommand{\invoroot}[2]{#1^{\invosymroot{#2}}}
\DeclarePairedDelimiterX{\brackets}[1]{(}{)}{#1}
\DeclarePairedDelimiterX{\abs}[1]{\lvert}{\rvert}{\ifblank{#1}{\:\cdot\:}{#1}}
\DeclarePairedDelimiterX{\gen}[1]{\langle}{\rangle}{#1}
\DeclarePairedDelimiterX{\commutator}[2]{[}{]}{\ifblank{#1}{\:\cdot\:}{#1}, \ifblank{#2}{\:\cdot\:}{#2}}
\newcommand{\commpart}[4][]{\commutator[#1]{#2}{#3}_{#4}}
\newcommand{\assoc}[3]{[#1,#2,#3]}
\DeclareMathOperator{\ad}{ad}
\newlength{\fixmidfigure}
\newcommand{\premidfigure}{\setlength{\fixmidfigure}{\lastskip}\addvspace{-\lastskip}}
\newcommand{\postmidfigure}{\addvspace{\fixmidfigure}}
\newcommand{\roots}{\Phi}
\newcommand{\rootbase}{\Delta}

\newcommand{\rootbaseE}{\rootbase_E}
\newcommand{\possys}{\Pi}
\newcommand{\refl}[1]{\sigma_{#1}}
\newcommand{\reflbr}[1]{\sigma(#1)}
\newcommand{\oprootint}[1]{\mathord{]#1[}}
\newcommand{\oprootintcry}[1]{\mathord{]#1[}_{\text{cry}}}
\newcommand{\subsys}[1]{\Psi_{#1}}
\newcommand{\htworootsym}{\chi}
\newcommand{\htworoot}[1]{\htworootsym_{#1}}

\DeclareMathOperator{\Weyl}{Weyl}
\newcommand{\rootbaseH}{\rootbase_H}
\newcommand{\rootbaseX}{\rootbase_X}

\newcommand{\possysH}{\possys_H}
\newcommand{\possysX}{\possys_X}
\newcommand{\paratwopos}{E_2^+}
\newcommand{\gold}{\tau}
\newcommand{\hcoord}[3]{[#1,#2,#3]}
\newcommand{\rootcoord}[1]{[#1]}

\newcommand{\goldfoldproj}{\tilde{\pi}}
\newcommand{\foldproj}{\pi}
\newcommand{\weylemb}{u}
\newcommand{\rootgr}[1]{U_{#1}}
\newcommand{\rootgrH}[1]{U_{#1}^H}
\newcommand{\rootgrGH}[1]{U_{#1}^{GH}}
\newcommand{\rootgrD}[1]{U_{#1}^D}

\newcommand{\rootgrX}[1]{U_{#1}^X}
\newcommand{\rootgrst}[1]{U_{#1}^{\text{St}}}
\newcommand{\rootgrHex}[1]{\tilde{U}_{#1}^H}
\newcommand{\rootgrDex}[1]{\tilde{U}_{#1}^D}
\newcommand{\rootgrAex}[1]{\tilde{U}_{#1}^A}
\newcommand{\rootgrEex}[1]{\tilde{U}_{#1}^E}

\newcommand{\rootgrunfold}[1]{U_{#1}^X}
\newcommand{\invset}[1]{U_{#1}^\sharp}
\newcommand{\risom}[1]{\theta_{#1}}
\newcommand{\risomH}[1]{\theta_{#1}^{H}}

\newcommand{\risomX}[1]{\theta_{#1}^X}
\newcommand{\risomXnaive}[1]{\theta_{#1}^{X, \text{n}}}
\newcommand{\risomHex}[1]{\tilde{\theta}_{#1}^{H}}
\newcommand{\risomDex}[1]{\tilde{\theta}_{#1}^D}
\newcommand{\risomDexnaive}[1]{\tilde{\theta}_{#1}^{D, \text{n}}}
\newcommand{\risomEex}[1]{\tilde{\theta}_{#1}^E}
\newcommand{\risomEexnaive}[1]{\tilde{\theta}_{#1}^{E, \text{n}}}

\newcommand{\risomXex}[1]{\tilde{\theta}_{#1}^X}
\newcommand{\risomst}[1]{\theta_{#1}^{\text{St}}}
\newcommand{\weylset}[1]{M_{#1}}
\newcommand{\commmap}[2]{\ifblank{#2}{\psi_{#1}}{\psi_{#1}^{#2}}}
\newcommand{\commmapex}[2]{\ifblank{#2}{\psi_{#1}}{\psi_{#1}^{#2}}}
\newcommand{\comone}{f}
\newcommand{\comtwo}{g}
\newcommand{\comthrone}{h_1}
\newcommand{\comthrtwo}{h_2}
\newcommand{\switch}{s}
\DeclareMathOperator{\El}{E}

\newcommand{\inverparsym}{\eta}
\newcommand{\inverpar}[2]{\inverparsym_{#1,#2}}
\newcommand{\inverparbr}[2]{\inverparsym(#1, #2)}
\newcommand{\twistgroup}{A}
\newcommand{\word}[1]{\bar{#1}}

\newcommand{\rwrule}[1]{\ifblank{#1}{\varphi}{\varphi_{#1}}}
\newcommand{\blumap}[1]{\ifblank{#1}{\gamma}{\gamma_{#1}}}
\newcommand{\longel}{w_0}
\newcommand{\ring}{\mathcal{R}}
\newcommand{\prodring}{\mathcal{S}}
\DeclareMathOperator{\ringproj}{pr}
\newcommand{\subisom}{\varphi}
\newcommand{\ringcoord}[2]{[#1,#2]}
\newcommand{\invring}[1]{#1^\times}
\newcommand{\scalmap}{\phi}
\newcommand{\Z}{\mathbb{Z}}
\newcommand{\C}{\mathbb{C}}

\newcommand{\R}{\mathbb{R}}
\newcommand{\N}{\mathbb{N}}
\newcommand{\inv}{^{-1}}

\newcommand{\ep}{\varepsilon}
\newcommand{\la}{\langle}
\newcommand{\ra}{\rangle}
\DeclareMathOperator{\id}{id}
\DeclareMathOperator{\Chev}{Ch}
\DeclareMathOperator{\Stein}{St}
\DeclareMathOperator{\GL}{GL}
\DeclareMathOperator{\SL}{SL}

\DeclareMathOperator{\rank}{rank}

\numberwithin{equation}{section}

\theoremstyle{plain}
\newtheorem{theorem}[equation]{Theorem}
\newtheorem{prop}[equation]{Proposition}
\newtheorem{lemma}[equation]{Lemma}

\theoremstyle{definition}
\newtheorem{definition}[equation]{Definition}
\newtheorem{convention}[equation]{Convention}
\newtheorem{notation}[equation]{Notation}
\newtheorem{rem}[equation]{Remark}
\newtheorem{note}[equation]{Note}

\newtheorem{example}[equation]{Example}

\newtheorem{algorithm}[equation]{Algorithm}
\newtheorem{conclusion}[equation]{Conclusion}
\newtheorem{goal}[equation]{Goal}

\allowdisplaybreaks

\makeatletter
\g@addto@macro\bfseries{\boldmath}
\makeatother

\title{Root graded groups of type \texorpdfstring{$H_3$}{H3} and \texorpdfstring{$H_4$}{H4}}

\author{Lennart Berg and Torben Wiedemann}

\fancypagestyle{copy}{\fancyhf{}\fancyfoot[C]{© 2025. This manuscript version is made available under the CC-BY-NC-ND 4.0 license \url{https://creativecommons.org/licenses/by-nc-nd/4.0/}.\\
Published in the Journal of Algebra, Volume 682: \url{https://doi.org/10.1016/j.jalgebra.2025.06.007}}}

\begin{document}
	\setlength{\abovedisplayskip}{0.5\abovedisplayskip}
	\setlength{\belowdisplayskip}{0.5\belowdisplayskip}
	\setlist{itemsep=0pt, topsep=2pt}

	\maketitle
	\thispagestyle{copy}

\begin{abstract}
	Using the well-known realisation of the root system $ H_4 $ as a folding of $ E_8 $, one can construct examples of $ H_4 $-graded groups from Chevalley groups of type $ E_8 $. Such Chevalley groups are defined over a commutative ring $ \ring $, and the root groups of the resulting $ H_4$-grading are coordinatised by $ \ring \times \ring $. We show that every $ H_4 $-graded group arises as the folding of an $ E_8 $-graded group, or in other words, that it is coordinatised by $ \ring \times \ring $ for some commutative ring $ \ring $. We also prove similar assertions for $ (D_6, H_3) $ in place of $ (E_8, H_4) $.
\end{abstract}

\section{Introduction}
\label{sec:intro}

In his seminal paper \cite{Chev-Tohoku}, Chevalley introduced the groups which are today known as \emph{Chevalley groups}. By construction, any such group is equipped with a family of subgroups $ (\rootgr{\alpha})_{\alpha \in \roots} $ indexed by a finite crystallographic reduced root system $ \roots $. These root groups are coordinatised by a commutative ring and the commutator of two root groups adheres to the celebrated \emph{Chevalley commutator formula}.

Motivated by these observations, we say that a \emph{$ \roots $-grading} of a group $ G $ is a family of non-trivial subgroups $ (U_\alpha)_{\alpha \in \roots} $ satisfying some (purely combinatorial) properties of the root groups in Chevalley groups. Here $ \roots $ is any finite root system which is not assumed to be crystallographic or reduced. Our definition of root graded groups was introduced in \cite{torben} and encompasses several existing concepts in the literature as special cases. In particular, it generalises Tits' notion of RGD-systems (see \cite[p.~258]{Tits-TwinBuildingsKacMoody}) and Shi's notion of groups graded by a root system (see \cite{Shi1993}). An overview of the historic background of root graded groups is given in the preface of \cite{torben}.

The main result of \cite{torben} says that if $ \roots $ is irreducible, crystallographic and of rank at least $ 3 $, then every $ \roots $-graded group is coordinatised by an algebraic structure (whose nature depends on $ \roots $). This algebraic structure is a division structure if and only if the corresponding $ \roots $-graded group is an RGD-system.

In this paper, we prove a similar result for the non-crystallo\-graphic root systems $ H_3 $ and $ H_4 $, thereby extending the coordinatisation results of \cite{torben} to the class of all irreducible finite root systems of rank at least~3. Specifically, we show in~\ref{thm:prod-param} that any root graded group $ G $ of type $ H_3 $ or $ H_4 $ is coordinatised by $ \ring \times \ring $ where $ \ring $ is a commutative ring (which depends on $ G $). Since $\ring \times \ring$ cannot be a division ring, it follows that there do not exist RGD-systems of types $H_3$ or $H_4$. This is a classical result of Tits \cite[Hauptsatz]{Tits-EndSpiegWeyl}, who proved it by showing that there are no RGD-systems of type $H_2$ \cite[Theorème~1]{Tits-NonEx}.

In the following, let $\ell \in \Set{3,4}$. Since $ H_\ell $ is not crystallographic, there are no Chevalley groups of type $ H_\ell $, so it is not evident a priori that there actually exist $ H_\ell $-graded groups. However, using the well-known realisation of the root system $ H_4 $ as a folding of $ E_8 $, examples of $ H_4 $-graded groups can be constructed from $E_8$-graded groups (such as, for example, Chevalley groups of type $ E_8 $). We will give a detailed description of this construction, which generalises the construction of the twisted Chevalley groups (see, for example, \cite[Chapter~13]{Carter-Chev}).
In a similar way, one obtains examples of $ H_3 $-gradings as foldings of $ D_6 $-gradings. It is a consequence of our main result that, in fact, all root gradings of type $ H_3 $ or $ H_4 $ arise as foldings of $D_6$- or $E_8$-gradings (see~\ref{thm:from-folding}).

The proof of our main result relies on the two main tools from \cite{torben}: The first one is the \emph{parametrisation theorem}, which yields canonical isomorphisms between the root groups of an arbitrary $ H_\ell $-graded group $ G $ and hence provides us with a parametrisation of the root groups by a single abelian group $ (\prodring, +) $. Our second main tool is the \emph{blueprint technique}, which is a computational method that allows us to equip $ (\prodring, +) $ with the structure of a commutative ring decomposing as $ \prodring = \prodring_1 \times \prodring_2 $ for isomorphic ideals $ \prodring_1 $ and $ \prodring_2 $. The blueprint technique follows the ideas of Ronan-Tits on the construction of buildings in \cite{titsronan}, but in a \enquote{reversed} setting.

In \cite[9.16]{mwpent}, it is shown that an analogue of our coordinatisation result for $ H_\ell $-graded groups holds for $ H_2 $-graded groups which arise from a Tits pentagon that satisfies a certain thickness condition. By \cite[Section~5]{mw19r2}, the $ H_2 $-graded groups arising from a Tits pentagon are precisely those which satisfy a certain stability condition. Our main result is true regardless of any stability assumptions, but in return, it uses in a crucial way that the root systems $ H_3 $ and $ H_4 $ are of rank higher than~2.

This paper is organised as follows. In Sections~\ref{sec:prelim} and \ref{sec:rgg}, we assemble some basic facts about root systems and introduce the notion of root graded groups. In Section~\ref{sec:fold}, we construct examples of $ H_\ell $-graded groups via foldings and explicitly compute the commutator formulas in these groups.
The groundwork for the coordinatisation of $ H_\ell $-graded groups is laid in Section~\ref{sec:h3}, where we investigate the commutator relations and Weyl elements in such groups. The main result of this section is~\ref{prop:h3-weyl-summary}, which provides formulas for the action of squares of Weyl elements on all root groups. These formulas allow us to apply the parametrisation theorem, which is proven in \cite{torben} for arbitrary root systems. In Section~\ref{sec:gpt}, we formulate the parametrisation theorem for the root systems $ H_3 $ and $H_4$ and, for the reader's convenience, sketch a simplified version of its proof in this special case.
We also define the notion of coordinatisations of $ H_\ell $-graded groups, which is studied in more detail in Section~\ref{sec:stein}.
In Section~\ref{sec:bp}, we describe the general idea of the blueprint technique and then, in Section~\ref{sec:bph3}, apply it to the special case of $ H_\ell $-graded groups. This leads to our main result, the coordinatisation theorem (Theorem~\ref{thm:prod-param}). As a corollary, we obtain in Section~\ref{sec:from-fold} that every root grading of type $H_3$ or $H_4$ is the folding of a grading of type $D_6$ or $E_8$, respectively (Theorem~\ref{thm:from-folding}).

\subsection*{Acknowledgements}

We are grateful to Bernhard Mühlherr for numerous insightful discussions and for his encouragement. We also want to thank the referee for carefully reading our paper and for making many valuable suggestions. The work of the first author was supported by the Studienstiftung des deutschen Volkes and the work of the second author was supported by the DFG Grant MU1281/7-1. 

\section{Preliminaries}
\label{sec:prelim}

\begin{convention}
\label{conv:comm}
Let $G$ be a group. For any subset $U $ of $ G$ we denote by $\la U\ra$ the subgroup of $G$ generated by $U$. Furthermore, for all $g,h\in G$, we put $g^h\coloneqq h\inv gh$ and $[g,h]\coloneqq g\inv h\inv gh$. 
\end{convention}

\begin{rem}
\label{rem:commrel}
With the conventions in \ref{conv:comm}, we have the following relations for all $g,g_1,g_2,h,h_1,h_2$ in any group $ G $:
	\begin{enumerate}[label=(\roman*)]
		\item \label{rem:commrel:conj}$g^h=[h,g\inv]g=g[g,h]$.
		\item $gh=hg^h$.
		\item \label{rem:commrel:comm}$gh=hg[g,h]=[h\inv,g\inv]hg$.
		\item $[g,h]\inv=[h,g]$.
		\item \label{rem:commrel:add-left}$[g_1g_2,h]=[g_1,h]^{g_2}[g_2,h]$.
		\item \label{rem:commrel:add-right}$[g,h_1h_2]=[g,h_2][g,h_1]^{h_2}$.
	\end{enumerate}
\end{rem}

\begin{convention}
	We only consider rings with identity in this paper, and homomorphisms of rings preserve the identity. Further, all rings are assumed to be associative except when we specifically refer to \enquote{nonassociative rings}.
\end{convention}

Throughout this paper, we will mostly use the standard terminology for finite root systems that are not necessarily crystallographic. However, though $ H_3 $ and $H_4$ are reduced root systems, it will be practical to consider the non-reduced root systems $ GH_3 $ and $GH_4$ (which we will introduce in~\ref{def:GH}) at some points. Hence we use the following definition of root systems.

\begin{definition}\label{def:rootsys}
	A \emph{root system} is a finite subset $ \roots $ of a Euclidean space $ (V, \cdot) $ such that $ \roots^{\reflbr{\alpha}} = \roots $ for all $ \alpha \in \roots $ where
	\[ \map{\reflbr{\alpha} \defl \refl{\alpha}}{V}{V}{v}{v^{\reflbr{\alpha}} \defl v - 2 \frac{\alpha \cdot v}{\alpha \cdot \alpha} \alpha} \]
	denotes the reflection along $ \alpha^\perp $. It is called \emph{reduced} if $ \roots \cap \R \alpha = \{ \alpha, -\alpha\} $ for all $ \alpha \in \roots $ and its Weyl group is $ \Weyl(\roots) \defl \gen{\refl{\alpha} \mid \alpha \in \roots} $. It is called \emph{crystallographic} if for all $\alpha, \beta \in \roots$, the number $2\frac{\alpha \cdot \beta}{\beta \cdot \beta}$ is an integer.
\end{definition}

\begin{rem}\label{rem:refl-right}
	Note our convention in~\ref{def:rootsys} that reflections act on $V$ from the right. Consequently, we have $\refl{v^w} = \refl{v}^w$ for all $v \in V \setminus \Set{0}$ and all orthogonal $w \in \GL(V)$, and not the formula $\refl{w(v)} = w \circ \refl{v} \circ w^{-1}$ that holds for left actions.
\end{rem}

\begin{notation}
\label{nota:subsystempsi}
	Let $ \roots $ be a (fixed) root system and let $ A \subseteq \roots $. Then $ \subsys{A} \defl \gen{A}_\R \cap \roots $ is called the \emph{root subsystem (of $ \roots $) spanned by $ A $}. In general, a \emph{root subsystem of $ \roots $} is the root subsystem spanned by some subset of $ \roots $. If $ A = \Set{\alpha_1, \ldots, \alpha_n} $, we will write $ \subsys{\alpha_1, \ldots, \alpha_n} $ for $ \subsys{A} $. If $ X_n $ is the type of a root system, we will often write \enquote{$ X_n $-subsystem} in place of \enquote{root subsystem of type~$ X_n $}.
\end{notation}

\begin{definition}\label{def:pair}
	Let $ \roots $, $ \roots' $ be root systems such that $ \roots' $ is of rank 2. A pair $ (\alpha, \beta) \in \roots^2 $ is called a \emph{$ \roots' $-pair (in $ \roots $)} if the root subsystem $ \subsys{\alpha,\beta} $ of $ \roots $ is of type $ \roots' $ and $ (\alpha, \beta) $ is a root base of this subsystem.
\end{definition}

\begin{notation}\label{def:paratwopos}
	Let $\roots$ be a root system, let $\rootbase$ be a root base of $\roots$ and denote by $\possys$ the corresponding positive subsystem of $\roots$. We put
	\[ \paratwopos(\rootbase) \defl \paratwopos(\rootbase, \roots) \defl \roots \cap \Set{a\delta + b\delta' \mid a,b \in \R_{\ge 0}, \delta, \delta' \in \rootbase} = \possys \cap \bigcup_{\delta, \delta' \in \rootbase} \subsys{\delta, \delta'}. \]
\end{notation}

\begin{rem}\label{rem:paratwopos}
	Let $\roots$ be a root system, let $\rootbase$ be a root base of $\roots$ and let $\alpha, \beta \in \roots$ be non-proportional. Then there exists $w \in \Weyl(\roots)$ such that $\alpha^w$ and $\beta^w$ both lie in $\paratwopos(\rootbase, \roots)$. For this reason, we will only describe the commutator relations in $H_\ell$-graded groups (for $\ell \in \Set{3,4}$) for pairs of roots in $\paratwopos(\rootbase,H_\ell)$ where $ \rootbase $ is some fixed root base of $ H_\ell $.
\end{rem}

\begin{notation}\label{def:rootcoord}
	Let $\roots$ be a root system of rank $n$ and let $\rootbase = (\delta_1, \ldots, \delta_n)$ be an ordered root base of $\roots$. For all $a \in \R^n$, we put $\rootcoord{a_1, \ldots, a_n}_\rootbase \defl \sum_{i=1}^n a_i \delta_i$. When $\rootbase$ is clear from the context, we will simply write $\rootcoord{a_1, \ldots, a_n}$ for $\rootcoord{a_1, \ldots, a_n}_\rootbase$.
\end{notation}

We briefly recall the notion of root intervals, which plays an important role in the theory of root gradings.

\begin{definition}\label{def:ri}
	Let $\roots$ be a root system and let $ \alpha, \beta \in \roots $ be non-proportional. The \emph{root interval between $\alpha$ and $\beta$} and the \emph{crystallographic root interval between $\alpha$ and $\beta$} are
	\[ \oprootint{\alpha, \beta} \defl \roots \cap \Set{a \alpha + b \beta \mid a,b \in \R_{>0}}, \quad \oprootintcry{\alpha, \beta} \defl \roots \cap \Set{a \alpha + b \beta \mid a,b \in \Z_{>0}}, \]
	respectively. Further, $ \alpha, \beta $ are called \emph{(crystallographically) adjacent} if $ \oprootint{\alpha, \beta} $ (respectively, $ \oprootintcry{\alpha, \beta} $) is empty.
\end{definition}

\begin{lemma}\label{def:interval-ord}
	Let $ \roots $ be a reduced root system, let $ \beta, \gamma \in \roots $ be non-proportional and denote by $ k $ the cardinality of $ \oprootint{\beta, \gamma} $. Then there exists a unique way to label the roots in $ \oprootint{\beta, \gamma} $ by $ \alpha_1, \ldots, \alpha_k $ such that $ \oprootint{\alpha_i, \alpha_j} = \Set{\alpha_p \mid i < p < j} $ for all $ i<j \in \Set{0, \ldots, k+1} $ where $ \alpha_0 \defl \beta $ and $ \alpha_{k+1} \defl \gamma $. The tuple $ (\alpha_1, \ldots, \alpha_k) $ is called the \emph{interval ordering of $ \oprootint{\beta, \gamma} $ starting from $ \beta $}.
\end{lemma}

We end this section with the definition of $ H_2 $, $ H_3 $ and $ H_4 $ and some of their basic properties.

\begin{notation}
\label{nota:h2}
We denote by $H_2$, $ H_3 $ and $ H_4 $ the root systems which correspond to the Coxeter diagrams \dynkin[Coxeter, backwards]{H}{2}, \dynkin[Coxeter, backwards]{H}{3} and \dynkin[Coxeter, backwards]{H}{4}, respectively. They are uniquely determined up to rescaling all roots by a common scalar, and we use the convention that all roots in $ H_2 $, $ H_3 $ and $ H_4 $ have length~$ 1 $. For an explicit construction of $H_3$ and $H_4$, see \cite[2.13]{HumphreysCox}. An ordered root base $ (\rho_1, \rho_2, \rho_3) $ of $ H_3 $ (respectively, $(\rho_0, \rho_1, \rho_2, \rho_3)$ of $H_4$) is said to be \emph{in standard order} if $(\rho_2, \rho_3)$ is an $H_2$-pair and every other pair of the form $(\rho_i, \rho_{i+1})$ is an $A_2$-pair. When a root base $(\rho_0, \rho_1, \rho_2, \rho_3)$ of $ H_4 $ in standard order is fixed, we call $ \subsys{\rho_1, \rho_2, \rho_3} $ the \emph{canonical $ H_3 $-subsystem of $ H_4 $}.
\end{notation}

\begin{rem}\label{lem:WH3transitive}
	The cardinalities of the root systems $ H_2 $, $ H_3 $ and $H_4$ are $ 10 $, $ 30 $ and $120$, respectively. Their Weyl groups have respective orders $10$, $120$ and $120 \cdot 120$. Further, a brief computation shows that $\Weyl(H_\ell)$ acts transitively on $H_\ell$ for $\ell=2$, which implies that the same holds for $\ell \in \Set{3,4}$ as well.
\end{rem}

\begin{rem}
\label{rem:noncryst}
In the literature, the root system $ H_2 $ is sometimes referred to as $I_2(5)$. Note that it is not crystallographic: If $ (\rho_2, \rho_3) $ denotes a root base of $ H_2 $, then $2\frac{\rho_2\cdot\rho_3}{\rho_3\cdot\rho_3}=2\cos(\frac{\pi}{5})=\frac{1+\sqrt{5}}{2}$ is not an integer. This number is known as the \emph{golden ratio}, and we will denote it by $\tau$ throughout the paper. Similarly, $ H_3 $ and $ H_4 $ are not crystallographic because they contain $H_2$.
\end{rem}

\begin{notation}\label{def:GH}
	Let $\ell \in \Set{2,3,4}$. The \emph{root system $GH_\ell$} is $GH_\ell \defl H_\ell \cup \gold H_\ell$ where $\gold \defl \frac{1+\sqrt{5}}{2} $ denotes the golden ratio. We also call it the \emph{golden version of $H_\ell$}. It is a non-reduced root system with the same Weyl group as $H_\ell$ and with two roots in each ray, whose lengths differ by a factor of $\gold$. When we speak of $H_\ell$ as a subset of $GH_\ell$, we will always mean the set of indivisible roots in $GH_\ell$ (which is not a subsystem in the sense of~\ref{nota:subsystempsi}). A \emph{root base of $GH_\ell$ (in standard order)} is a root base of the subset $H_\ell$ (in standard order). The map $\map{\scalmap}{GH_\ell}{H_\ell}{\alpha}{\frac{1}{\lVert \alpha \rVert} \alpha}$ which sends each root in $GH_\ell$ to the indivisible root in its ray is called the \emph{scaling map}.
\end{notation}

\begin{rem}
	See Figure~\ref{fig:gh2} for a graphic depiction of $ GH_2 $ and \cite{mwpent} for more details on $ GH_2 $, which inspired our definition of $ GH_3 $ and $ GH_4 $.
\end{rem}

\premidfigure
\begin{figure}[htb]
\centering
	\begin{tikzpicture}[scale=1.25]
		\foreach \angle [count=\ani] in {0, 36, 72, 108, 144, 180, 216, 252, 288, 324}
		{
			\begin{scope}[rotate=\angle]
				\draw [arrows = {-Stealth[scale=0.75]}] (0,0) -- (1,0) node(short-\ani){};
				\draw [arrows = {-Stealth[scale=0.75]}] (0,0) -- (1.618033988749895,0) node (\ani) {};
			\end{scope}
		}
		
		\foreach \root/\angle [count=\ani] in {\alpha/0, \beta/36, \gamma/72, \delta/108, \ep/144}{
			\begin{scope}[font=\footnotesize, rotate=\angle]
				\node at (1.775, 0){$ \tau \root $};
				\node at (0.8, 0.15){$ \root $};
			\end{scope}
		}
		
		\begin{scope}[densely dotted, thin]
			\draw (1) -- (4);
			\draw (2) -- (5);
			\draw (3) -- (6);
			\draw (4) -- (7);
			\draw (5) -- (8);
			\draw (6) -- (9);
			\draw (7) -- (10);
			\draw (8) -- (1);
			\draw (9) -- (2);
			\draw (10) -- (3);
		\end{scope}
		
		\begin{scope}[densely dashed, thin]
			\draw (1) -- (2);
			\draw (2) -- (3);
			\draw (3) -- (4);
			\draw (4) -- (5);
			\draw (5) -- (6);
			\draw (6) -- (7);
			\draw (7) -- (8);
			\draw (8) -- (9);
			\draw (9) -- (10);
			\draw (10) -- (1);
		\end{scope}
	\end{tikzpicture}
\caption{The root system $GH_2$ with an $ H_2 $-quintuple $ (\alpha, \beta, \gamma, \delta, \ep) $.}
\label{fig:gh2}
\end{figure}
\postmidfigure

\begin{rem}\label{rem:gap}
	Throughout this paper, we will make several assertions which follow from a straightforward but lengthy computation. The interested reader may find GAP~\cite{GAP} code which verifies these assertions in \cite{repo}.
\end{rem}

\begin{prop}
\label{prop:222}
	Let $\alpha$ be any root in $ H_3 $. Then the following hold:
		\begin{enumerate}[(i)]
			\item $\alpha$ is contained in exactly $2$ subsystems of type $H_2$.
			\item $\alpha$ is contained in exactly $2$ subsystems of type $A_2$.
			\item \label{prop:222:A1xA1}$\alpha$ is contained in exactly $2$ subsystems of type $A_1\times A_1$.
		\end{enumerate}
\end{prop}
\begin{proof}
	By~\ref{lem:WH3transitive}, we can choose $\rho_2, \rho_3 \in H_3$ such that $(\alpha, \rho_2, \rho_3)$ is a root base of $H_3$ in standard order. Now the assertion follows from a straightforward inspection of $ H_3 $ (see~\ref{rem:gap} for details).
\end{proof}

\begin{definition}
\label{def:i25quintuple}
 If $ (\alpha, \ep) $ is an $ H_2 $-pair in $ \roots $, then the \emph{$ H_2 $-quintuple associated to $ (\alpha, \ep) $} is the tuple $(\alpha,\beta,\gamma,\delta,\ep)$ where $\beta\coloneqq\tau\alpha+\ep$, $\gamma\coloneqq\tau\alpha+\tau\ep$ and $\delta\coloneqq\alpha+\tau\ep$. Note that $ (\beta, \gamma, \delta) $ is the interval ordering of $ \oprootint{\alpha, \ep} $ starting from $ \alpha $. Similarly, the \emph{$ A_2 $-triple associated to an $ A_2 $-pair $ (\alpha,\gamma) $} is $ (\alpha, \alpha+\gamma, \gamma) $.
\end{definition}

\section{Root graded groups}
\label{sec:rgg}

\begin{notation}
	Throughout this section, we denote by $ \roots $ an arbitrary reduced root system with root base $ \rootbase $. From after~\ref{def:risom} on, we denote by $ (G, (\rootgr{\alpha})_{\alpha \in \roots}) $ a $ \roots $-graded group.
\end{notation}

In this section, we introduce $\roots$-graded groups in the sense of \cite[2.5.2]{torben}. They generalise Tits' notion of RGD-systems (see~\ref{rem:rgd} and~\cite[p.~258]{Tits-TwinBuildingsKacMoody}).

\begin{notation}
	If $ G $ is a group, $(U_\alpha)_{\alpha\in\Phi}$ is a family of subgroups and $ A $ is a subset of $ \roots $, we put $U_A\coloneqq\la U_\alpha\mid\alpha\in A\ra$.
\end{notation}

\begin{definition}
\label{def:comm-rel}
Let $G$ be a group and let $(U_\alpha)_{\alpha\in\Phi}$ be a family of non-trivial subgroups. We say that $ G $ has \emph{$ \roots $-commutator relations with root groups $ (\rootgr{\alpha})_{\alpha \in \roots} $} if $[U_\alpha,U_\beta]\subseteq U_{\oprootint{\alpha, \beta}}$ for all non-proportional $ \alpha, \beta \in \roots $.
If, moreover, we have $[U_\alpha,U_\beta]\subseteq U_{\oprootintcry{\alpha, \beta}}$ for all non-proportional $ \alpha, \beta \in \roots $, we say that $ G $ has \emph{crystallographic $ \roots $-commutator relations with root groups $ (\rootgr{\alpha})_{\alpha \in \roots} $}.
\end{definition}

\begin{definition}\label{def:rgg}
	Let $G$ be a group. A \emph{(crystallographic) $ \roots $-grading} is a family $(U_\alpha)_{\alpha\in\Phi}$ of non-trivial subgroups satisfying the following conditions.
	\begin{enumerate}[label=(\roman*)]
		\item The group $G$ is generated by $(U_\alpha)_{\alpha\in\Phi}$.
		\item The group $ G $ has (crystallographic) $ \roots $-commutator relations with root groups $(U_\alpha)_{\alpha\in\Phi}$.
		\item For all $\alpha\in\Phi$, the set
		\[M_\alpha\coloneqq\{w_\alpha\in U_{-\alpha}U_\alpha U_{-\alpha}\mid U_\beta^{w_\alpha}=U_{\beta^{\reflbr{\alpha}}}\text{ for all }\beta\in\Phi\}\]
		of so-called \emph{$ \alpha $-Weyl elements} is non-empty.
		\item \label{def:rgg:pos}$U_\Pi\cap U_{-\alpha}=\{1_G\}$ for all positive systems of roots $\Pi$ and all $ \alpha \in \Pi $.
	\end{enumerate}
We will also refer to $ (G, (\rootgr{\alpha})_{\alpha \in \roots}) $ or $G$ as a \emph{(crystallographic) $\Phi$-graded group} if these conditions are satisfied. Further, we define
	\[U_\alpha^\sharp\coloneqq\{u_\alpha\in U_\alpha\mid\text{There exist }u_{-\alpha},u_{-\alpha}'\in U_{-\alpha}\text{ such that }u_{-\alpha}u_\alpha u_{-\alpha}'\in M_\alpha\}\]
	for all $ \alpha \in \roots $.
\end{definition}

\begin{definition}\label{def:risom}
	Let $(G, (\rootgr{\alpha})_{\alpha \in \roots})$ a $\roots$-graded group and let $\alpha \in \roots$. A \emph{root isomorphism (for $\alpha$)} is an isomorphism $\map{\risom{\alpha}}{M}{\rootgr{\alpha}}{}{}$ for a group $M$. Given such an isomorphism, we also say that \emph{$\rootgr{\alpha}$ is parametrised by $M$}.
\end{definition}

For the rest of this section, we denote by $ (G, (\rootgr{\alpha})_{\alpha \in \roots}) $ an arbitrary $ \roots $-graded group.

\begin{rem}\label{rem:rootgr-abel}
	Note that the root groups in a $\roots$-graded group are not assumed to be abelian. However, we will see in~\ref{prop:a2gradprop}~\ref{prop:a2gradprop:abel} that $\rootgr{\alpha}$ is necessarily abelian if $\alpha$ is contained in an $A_2$-subsystem. In particular, root gradings of types $H_3$ and $H_4$ have abelian root groups.
\end{rem}

\begin{rem}\label{rem:pos-trans}
	Since the Weyl group acts transitively on the set of positive systems in $ \roots $, the existence of Weyl elements implies that it suffices to require that Axiom~\ref{def:rgg}~\ref{def:rgg:pos} holds for some positive system $ \Pi $.
\end{rem}

We now collect some basic properties of root graded groups.

\begin{prop}[{\cite[2.5.11]{torben}}]
\label{prop:parabolicsubgrading}
Let $\Psi$ be a root subsystem of $ \roots $ and put $G_\Psi\coloneqq\la U_\gamma\mid\gamma\in \Psi\ra$. Then $(U_\gamma)_{\gamma\in\Psi}$ is a $\Psi$-grading of $G_\Psi$.
\end{prop}

\begin{lemma}[{\cite[2.1.13]{torben}}]
\label{lem:oneadditivity}
Let $\alpha,\gamma\in\Phi$ such that the root interval $\oprootint{\alpha, \gamma}$ contains exactly one element and let $x_\alpha,x_\alpha'\in U_\alpha$, $y_\gamma,y_\gamma'\in U_\gamma$. Then
	\[[x_\alpha x_\alpha',y_\gamma]=[x_\alpha,y_\gamma][x_\alpha',y_\gamma]\hspace*{0.4cm}\text{and}\hspace*{0.4cm}[x_\alpha,y_\gamma y_\gamma']=[x_\alpha,y_\gamma][x_\alpha,y_\gamma'].\]
In particular, $[x_\alpha\inv,y_\gamma]=[x_\alpha,y_\gamma]\inv=[x_\alpha,y_\gamma\inv]$.
\end{lemma}
\begin{proof}
	Denote by $ \beta $ the unique root in $ \oprootint{\alpha,\gamma} $. By~\ref{rem:commrel},
	\[ [x_\alpha x_\alpha',y_\gamma]=[x_\alpha,y_\gamma]^{x_\alpha'}[x_\alpha',y_\gamma] \quad \text{and} \quad [x_\alpha,y_\gamma y_\gamma']=[x_\alpha,y_\gamma][x_\alpha,y_\gamma']^{y_\gamma'}. \]
	Since $U_\beta$ commutes with $U_\alpha$ and $U_\gamma$, the assertion follows.
\end{proof}

\begin{prop}[{\cite[2.2.6]{torben}}]
\label{prop:weylelements}
Let $\alpha\in\Phi$. Then the following hold:
\begin{enumerate}[label=(\roman*)]
	\item $M_\alpha=M_\alpha\inv\coloneqq\{w_\alpha\inv\mid w_\alpha\in M_\alpha\}$ and $U_\alpha^\sharp=(U_\alpha^\sharp)\inv$.
	\item \label{prop:weylelements:conj}$M_\alpha^{w_\beta}=M_{\beta^{\reflbr{\alpha}}}$ and $(U_\alpha^\sharp)^{w_\beta}=U_{\beta^{\reflbr{\alpha}}}^\sharp$ for all $\beta\in\Phi$ and $w_\beta\in M_\beta$.
	\item \label{prop:weylelements:minus}Let $w_\alpha$ be an $\alpha$-Weyl element with factorization $w_\alpha=a_{-\alpha} b_\alpha c_{-\alpha}$ for some $b_\alpha\in U_\alpha$, $a_{-\alpha},c_{-\alpha}\in U_{-\alpha}$. Then $w_\alpha=c_{-\alpha}^{w_\alpha\inv}a_{-\alpha}b_\alpha$ and $w_\alpha=b_\alpha c_{-\alpha}a_{-\alpha}^{w_\alpha}$. In particular, $w_\alpha$ is also a ($-\alpha$)-Weyl element and $ M_\alpha = M_{-\alpha} $.
	\item Let $b_\alpha\in U_\alpha^\sharp$ and choose $a_{-\alpha},c_{-\alpha}\in U_{-\alpha}$ such that $a_{-\alpha}b_\alpha c_{-\alpha}$ is an $\alpha$-Weyl element. Then $a_{-\alpha},c_{-\alpha}\in U_{-\alpha}^\sharp$.
\end{enumerate}
\end{prop}

The following two assertions are consequences of Axiom~\ref{def:rgg}~\ref{def:rgg:pos}. The second, more general assertion will only be needed in~\ref{lem:stein-epi} and~\ref{prop:stein-is-rgg}.

\begin{prop}[{\cite[2.5.9, 2.1.19]{torben}}]
\label{prop:prodmapbij}
	Let $ \beta, \gamma \in \roots $ be non-proportional and denote by $ (\alpha_1, \ldots, \alpha_k) $ the interval ordering of $ \oprootint{\beta, \gamma} $ starting from $ \beta $ (as in \ref{def:interval-ord}). Then the product map
	\[ \map{}{\rootgr{\alpha_1} \times \cdots \times \rootgr{\alpha_k}}{\rootgr{\oprootint{\beta, \gamma}}}{(x_1, \ldots, x_k)}{x_1 \cdots x_k} \]
	is bijective.
\end{prop}

\begin{prop}\label{prop:prodmapbij:possys}
	Let $ \possys $ be a positive system in $ \roots $ and put $ k \defl \abs{\possys} $. Then there exists an ordering $ \alpha_1,\ldots, \alpha_k $ of the roots in $ \possys $ such that for every group $ G' $ which has $ \roots $-commutator relations with root groups $ (\rootgr{\alpha}')_{\alpha \in \roots} $, the following properties hold:
	\begin{enumerate}[(a)]
		\item \label{prop:prodmapbij:possys:surj}The product map
		\[ \map{}{\rootgr{\alpha_1}' \times \cdots \times \rootgr{\alpha_k}'}{\rootgr{\possys}'}{(x_1, \ldots, x_k)}{x_1 \cdots x_k} \]
		is surjective.
		
		\item \label{prop:prodmapbij:possys:bij}If $ (\rootgr{\alpha}')_{\alpha \in \roots} $ satisfies Axiom~\ref{def:rgg}~\ref{def:rgg:pos}, then the product map above is bijective.
	\end{enumerate}
\end{prop}
\begin{proof}
	By \cite[2.3.23]{torben}, there exists an ordering $ \alpha_1,\ldots, \alpha_k $ which is extremal in the sense of \cite[2.3.16]{torben}. For any such ordering, the product map is surjective by \cite[2.3.21,~2.4.14]{torben}. By~\cite[2.4.17]{torben}, the product map is even bijective for extremal orderings if Axiom~\ref{def:rgg}~\ref{def:rgg:pos} is satisfied.
\end{proof}

Using~\ref{prop:prodmapbij}, we can introduce the following notation, which will be used frequently.

\begin{notation}\label{nota:commpart}
	Let $ \beta, \gamma \in \roots $ be non-proportional and denote by $ (\alpha_1, \ldots, \alpha_k) $ the interval ordering of $ \oprootint{\beta, \gamma} $ starting from $ \beta $ (as in \ref{def:interval-ord}). Let $ x_\beta \in \rootgr{\beta} $ and $ x_\gamma \in \rootgr{\gamma} $. By~\ref{prop:prodmapbij}, there exist unique elements $ x_i \in \rootgr{\alpha_i} $ for all $ i \in \Set{1, \ldots, k} $ such that $ [x_\beta, x_\gamma] = x_1 \cdots x_k $. We put $ \commpart{x_\beta}{x_\gamma}{\alpha_i} \defl x_i $ for all $ i \in \Set{1, \ldots, n} $. When we write $ \commpart{x_\beta}{x_\gamma}{\alpha_i}^g $ for some $ g \in G $ or $ \commpart{x_\beta}{x_\gamma}{\alpha_i}^{-1} $, we always mean $ (\commpart{x_\beta}{x_\gamma}{\alpha_i})^g $ or $ (\commpart{x_\beta}{x_\gamma}{\alpha_i})^{-1} $, respectively.
\end{notation}

\begin{rem}\label{rem:commpart-notation-rem}
	In general, $ \commpart{x_\beta}{x_\gamma}{\alpha_i} $ depends on the (ordered) tuple $ (x_\beta, x_\gamma) $ and not only on the element $ \commutator{x_\beta}{x_\gamma} $ because the interval ordering $ (\alpha_1, \ldots, \alpha_k) $ depends on the starting root $ \beta $.
\end{rem}

\begin{lemma}[{\cite[2.1.18, 2.2.9]{torben}}]\label{lem:commpart}
	Let $ \alpha, \gamma \in \roots $ be non-proportional, let $ \beta \in \oprootint{\alpha, \gamma} $ and let $ x_\alpha \in \rootgr{\alpha} $, $ x_\gamma \in \rootgr{\gamma} $. Then for any root $ \delta $ and any $ \delta $-Weyl element $ w_\delta $, we have
	\[ \commpart{x_\alpha}{x_\gamma}{\beta}^{w_\delta} = \commpart{x_\alpha^{w_\delta}}{x_\gamma^{w_\delta}}{\beta^{\reflbr{\delta}}} \quad \text{and} \quad \commpart{x_\alpha}{x_\gamma}{\beta}^{w_\delta^2} = \commpart{x_\alpha^{w_\delta^2}}{x_\gamma^{w_\delta^2}}{\beta}. \]
	Further, $ \commpart{x_\alpha}{x_\gamma}{\beta}^{-1} = \commpart{x_\gamma}{x_\alpha}{\beta} $.
\end{lemma}

\begin{note}
	The assertion that $ \commpart{x_\alpha}{x_\gamma}{\beta}^{-1} = \commpart{x_\gamma}{x_\alpha}{\beta} $ is stated in \cite[2.1.18]{torben} only for the case that $ \oprootint{\alpha, \gamma} $ has at most two elements, but the same proof remains valid in general. This fact is related to~\ref{rem:commpart-notation-rem}.
\end{note}

The following notation is a variation of~\ref{nota:commpart}.

\begin{definition}\label{def:commmap}
	Assume that $ (G, (\rootgr{\alpha})_{\alpha \in \roots}) $ is equipped with a family of root isomorphisms $ (\map{\risomX{\alpha}}{M}{\rootgrX{\alpha}}{}{})_{\alpha \in X} $ for a group $ M $.
	For all non-proportional roots $ \alpha, \gamma \in \roots $ and all $ \beta \in \oprootint{\alpha, \gamma} $, we define a map $ \map{\commmap{\alpha, \gamma}{\beta}}{M \times M}{M}{}{} $ by the property that
	\[ \commpart{\risom{\alpha}(x)}{\risom{\gamma}(y)}{\beta} = \risom{\beta}\brackets[\big]{\commmap{\alpha, \gamma}{\beta}\brackets{x,y}} \]
	for all $ x,y \in M $. It is called the \emph{commutation map of $ G $ at position $ (\alpha, \gamma, \beta) $ (with respect to $ (\risom{\delta})_{\delta \in \roots} $)}. If $ \oprootint{\alpha, \gamma} $ has only one element $ \beta $, we will write $ \commmapex{\alpha, \gamma}{} $ in place of $ \commmapex{\alpha, \gamma}{\beta} $.
\end{definition}

The well-known \emph{braid relations} in Coxeter groups also hold for Weyl elements in root graded groups. This result goes back to \cite[(6.9)]{MoufangPolygons}, where it is stated and proven in a less general context. However, essentially the same proof remains valid in our setting, as is verified in \cite[2.2.34]{torben}.

\begin{prop}[{\cite[(6.9)]{MoufangPolygons}}]\label{braidrel}
	Let $ \alpha, \beta \in \rootbase $ and denote by $ m $ the order of $ \refl{\alpha} \refl{\beta} $ in the Weyl group. Then for all Weyl elements $ w_\alpha \in \weylset{\alpha} $, $ w_\beta \in \weylset{\beta} $, we have
	\[ w_\alpha w_\beta w_\alpha \cdots = w_\beta w_\alpha w_\beta \cdots \]
	where the products on both sides have length $ m $. In other words, $ (w_\alpha w_\beta)^{m/2} = (w_\beta w_\alpha)^{m/2} $ if $ m $ is even and $ (w_\alpha w_\beta)^{(m-1)/2} w_\alpha = (w_\beta w_\alpha)^{(m-1)/2} w_\beta $ if $ m $ is odd.
\end{prop}

\section{Foldings}
\label{sec:fold}

\begin{notation}\label{not:root-base-order}
	In this section, we consider the root systems $ A_4 $, $ D_6 $ and $ E_8 $ with the labelling of simple roots given by Figure~\ref{fig:fold-diags}. Note that our ordering of simple roots in $ E_8 $ is non-standard: It is chosen so that we can identify the simple roots $ \delta_1, \ldots, \delta_6 $ in $ D_6 $ with the simple roots $ \delta_1, \ldots, \delta_6 $ in $ E_8 $. We will also use the standard representations
	\begin{align*}
		A_{n-1} &= \Set{e_i - e_j \mid i \ne j \in \Set{1, \ldots, n}}, \quad D_n = \Set{\pm e_i \pm e_j \mid i < j \in \Set{1, \ldots, n}},
	\end{align*}
	where $n \ge 4$ and $ (e_1, \ldots, e_n) $ denotes the standard basis of $ \R^n $ (with the standard inner product). Further, we fix a root base $\rootbaseH = (\rho_0, \rho_1, \rho_2,  \rho_3)$ of $H_4$ in standard order and we make the identifications $H_2 = \subsys{\rho_2, \rho_3}$ and $H_3 = \subsys{\rho_1, \rho_2, \rho_3}$. We will also consider $H_2$, $H_3$ and $H_4$ as subsets of $GH_2$, $GH_3$ and $GH_4$.
	
	Certain arguments and remarks will apply to all three root systems $A_4$, $D_6$ and $E_8$ simultaneously. In such situations, we simply write $X$ to denote an arbitrary root system $X \in \Set{A_4, D_6, E_8}$. We also put $\ell \defl \rank(X)/2$ (because we will construct a folding $\map{}{X}{H_\ell}{}{}$) and denote by $\rootbaseX$ the ordered root base of $X$ given by Figure~\ref{fig:fold-diags}.
\end{notation}

\premidfigure
\begin{figure}[htb]
	\centering\begin{subfigure}{0.32\linewidth}
		\centering\begin{tikzpicture}[scale=0.3]
			\filldraw [thick] (-1.5,1) circle [radius=3pt];
			\filldraw [thick] (1.5,1) circle [radius=3pt];
			\filldraw [thick] (-1.5,-1) circle [radius=3pt];
			\filldraw [thick] (1.5,-1) circle [radius=3pt];
			\draw (-1.5,1)--(1.5,1);
			\draw (-1.5,-1)--(1.5,1);
			\draw (-1.5,-1)--(1.5,-1);
			\draw (-1.5,0) ellipse [x radius=15pt, y radius=44pt];
			\draw (1.5,0) ellipse [x radius=15pt, y radius=44pt];
			\begin{scope}[font=\footnotesize]
				\node at (-1.5, 2.2){$ \delta_2 $};
				\node at (1.5, 2.2){$ \delta_3 $};
				\node at (-1.5, -2.2){$ \delta_4 $};
				\node at (1.5, -2.2){$ \delta_5 $};
			\end{scope}
		\end{tikzpicture}
		\caption{The folding $ A_4 \rightarrow H_2 $.}
		\label{fig:fold-diags:H2}
	\end{subfigure}%
	\begin{subfigure}{0.32\linewidth}
		\centering\begin{tikzpicture}[scale=0.3]
			\filldraw [thick] (0,1) circle [radius=3pt];
			\filldraw [thick] (0,-1) circle [radius=3pt];
			\filldraw [thick] (3,1) circle [radius=3pt];
			\filldraw [thick] (3,-1) circle [radius=3pt];
			\filldraw [thick] (-3,1) circle [radius=3pt];
			\filldraw [thick] (-3,-1) circle [radius=3pt];
			\draw (0,1)--(3,1);
			\draw (0,-1)--(3,1);
			\draw (0,-1)--(3,-1);
			\draw (-3,1)--(0,1);
			\draw (-3,-1)--(0,-1);
			\draw (0,0) ellipse [x radius=15pt, y radius=44pt];
			\draw (3,0) ellipse [x radius=15pt, y radius=44pt];
			\draw (-3,0) ellipse [x radius=15pt, y radius=44pt];
			\begin{scope}[font=\footnotesize]
				\node at (-3, 2.2){$ \delta_1 $};
				\node at (0, 2.2){$ \delta_2 $};
				\node at (3, 2.2){$ \delta_3 $};
				\node at (-3, -2.2){$ \delta_6 $};
				\node at (0, -2.2){$ \delta_4 $};
				\node at (3, -2.2){$ \delta_5 $};
			\end{scope}
		\end{tikzpicture}
		\caption{The folding $ D_6 \rightarrow H_3 $.}
		\label{fig:fold-diags:H3}
	\end{subfigure}%
	\begin{subfigure}{0.32\linewidth}
		\centering\begin{tikzpicture}[scale=0.3]
			\filldraw [thick] (0,1) circle [radius=3pt];
			\filldraw [thick] (0,-1) circle [radius=3pt];
			\filldraw [thick] (3,1) circle [radius=3pt];
			\filldraw [thick] (3,-1) circle [radius=3pt];
			\filldraw [thick] (-3,1) circle [radius=3pt];
			\filldraw [thick] (-3,-1) circle [radius=3pt];
			\filldraw [thick] (-6,1) circle [radius=3pt];
			\filldraw [thick] (-6,-1) circle [radius=3pt];
			\draw (-6,1) -- (-3,1) -- (0,1) -- (3,1) -- (0,-1) -- (-3,-1) -- (-6,-1);
			\draw (3,-1) -- (0,-1);
			\draw (0,0) ellipse [x radius=15pt, y radius=44pt];
			\draw (3,0) ellipse [x radius=15pt, y radius=44pt];
			\draw (-3,0) ellipse [x radius=15pt, y radius=44pt];
			\draw (-6,0) ellipse [x radius=15pt, y radius=44pt];
			\begin{scope}[font=\footnotesize]
				\node at (-6, 2.2){$ \delta_7 $};
				\node at (-3, 2.2){$ \delta_1 $};
				\node at (0, 2.2){$ \delta_2 $};
				\node at (3, 2.2){$ \delta_3 $};
				\node at (0, -2.2){$ \delta_4 $};
				\node at (-3, -2.2){$ \delta_6 $};
				\node at (-6, -2.2){$ \delta_8 $};
				\node at (3, -2.2){$ \delta_5 $};
			\end{scope}
		\end{tikzpicture}
		\caption{The folding $ E_8 \rightarrow H_4 $.}
		\label{fig:fold-diags:H4}
	\end{subfigure}
	\caption{Three foldings.}
	\label{fig:fold-diags}
\end{figure}
\postmidfigure

Given a root system $ \roots $ and an automorphism $ \rho $ of its Dynkin diagram, one can construct a new root system $ \roots' $ by projecting $ \roots $ orthogonally to the fixed space of an isometry $ \tau $ which is related to $ \rho $ (with $ \tau = \rho $ if $ \roots $ is simply-laced, which is the only case we will be interested in). In particular, we have a map $ \map{}{\roots}{\roots'}{}{} $ which we call a \emph{folding}, and we will also refer to $ \roots' $ as a \emph{folding of $ \roots $}. Further, any Chevalley group of type $ \roots $ gives rise to a so-called \emph{twisted Chevalley group of type $ \roots' $}. This procedure is described in \cite[Chapter~13]{Carter-Chev}. In a similar way, every $ \roots $-grading $ (\rootgr{\alpha})_{\alpha \in \roots} $ of a group $ G $ gives rise to a $ \roots' $-grading $ (\rootgr{\alpha'}')_{\alpha' \in \roots'} $ of (the same group) $ G $.

A more general class of foldings, which do not arise from diagram automorphisms, is described in \cite{coxcox}. In this section, we illustrate how, in the notation of~\ref{not:root-base-order}, the root system $H_\ell$ can be realised in this way as a folding of $X \in \Set{A_4, D_6, E_8}$ (see~\ref{pro:fold-abstract-description}).
We will also describe a procedure to construct $H_\ell$-graded groups from $X$-graded groups in a way which generalises the aforementioned construction of twisted Chevalley groups (see~\ref{def:rootgr-fold} and~\ref{def:param-folding}).
This general procedure can be applied to the example of Chevalley groups of type $X$, and we will compute the commutator relations in the foldings of these specific examples (see~\ref{lem:excommrel}) as well as some more structure constants (see~\ref{rem:fold-ex-parmap}).

The main result of this section is~\ref{folding-is-rgg}, which shows that the folding procedure actually produces an $H_\ell$-grading in the sense of~\ref{def:rgg}.
Some technicalities of the proof are only dealt with in the appendix, which is independent of the following sections of this paper.
However, even without the appendix, the arguments (and computations) in this section are sufficient to prove that the foldings of the aforementioned Chevalley groups are $H_\ell$-gradings.

Since the root groups of a Chevalley group (of type $X$) are parametrised by a commutative ring $\ring$ (in the sense of~\ref{def:risom}), the root groups of the corresponding $H_\ell$-grading are parametrised by $\ring \times \ring$ in a natural way. Hence we say that the folding of a Chevalley group of type $X$ is \emph{coordinatised by $\ring \times \ring$}, a notion that we will make precise in~\ref{def:stand-coord}. The main result of this paper is that if $\ell \ge 3$, then every $ H_\ell $-graded group has a coordinatisation by $ \ring \times \ring $ for a commutative ring $ \ring $ (Theorem~\ref{thm:prod-param}). This means, essentially, that the commutator relations in an arbitrary $H_\ell$ grading have the same form as in~\ref{lem:excommrel}. It is a consequence of this result that for $\ell \ge 3$, every $ H_\ell $-graded group arises from an $ X $-graded group (but not necessarily from a Chevalley group, see~\ref{rem:not-chev}) by the folding procedure described in this section (Theorem~\ref{thm:from-folding}). Further, we point out that by our construction in this section, there exists an $ H_\ell $-graded group coordinatised by $\ring \times \ring$ for every commutative ring $\ring$ (and even for every ring $\ring$ if $\ell = 2$).

We will obtain the root system $GH_\ell$ as a natural intermediate step in the folding $\map{}{X}{H_\ell}{}{}$. Essentially, each root in $H_\ell$ is associated to two roots in $GH_\ell$, with short roots in $GH_\ell$ corresponding to the left factor of $\ring \times \ring$ in our coordinatisation and long roots corresponding to the right factor. Hence the root system $GH_\ell$ allows us to see the structure of an $H_\ell$-graded group $G$ more clearly if $G$ is obtained via folding (which, by the previous paragraph, is always the case for $\ell \ge 3$). See~\ref{note:preim-order} and~\ref{rem:exgh3commrel} for more details.

We point out that our focus in this section lies on the cases $ X=D_6 $ and $ X=E_8 $. The case $ X=A_4 $ is included mainly because this higher generality does not introduce any additional difficulties.

\begin{note}
	While Chevalley groups over fields are well-known objects (see, for instance, \cite{Carter-Chev,Steinberg-ChevGroups}), the notion of Chevalley groups over commutative rings is much less standard. For an introduction, see \cite{Vavilov-Plotkin}.
	
	To be more precise, the objects that we call \enquote{Chevalley groups} are called \enquote{elementary Chevalley groups} in \cite{Vavilov-Plotkin}. Chevalley groups of type $\roots$ have a natural $\roots$-grading $(\rootgr{\alpha})_{\alpha \in \roots}$. Specifically, (elementary) Chevalley groups are generated by their root groups by definition, they have $\roots$-commutator relations by the Chevalley commutator formula and Weyl elements exist by \cite[(13.1)]{Vavilov-Plotkin}. While we are not aware of a reference which explicitly states Axiom~\ref{def:rgg}~\ref{def:rgg:pos} for Chevalley groups over rings, the proof in \cite[Corollary~3 of Lemma~17]{Steinberg-ChevGroups} for Chevalley groups over fields remains valid in this generality.
	
	We also point out that we provide explicit matrix representation of (certain) Chevalley groups of types $A_4$ and $D_6$ in~\ref{ex:a4h2folding} and~\ref{ex:chev-D6}. Hence our treatment of these cases is self-contained.
\end{note}

As an introductory example, we begin with the construction of an $H_2$-graded group from a Chevalley group of type $A_4$.

\begin{example}
\label{ex:a4h2folding}
Let $ \ring $ be a commutative ring. For all $1\leq i\neq j\leq 5$, we denote by $e_{ij}$ the $ (5 \times 5) $-matrix with $1_\ring$ at position $(i,j)$ and $0_\ring$ at every other position, and we define a root group
	\[ \rootgrAex{\alpha} \coloneqq \rootgrAex{ij}\coloneqq \{\mathord{\id}+ae_{ij}\mid a\in \ring\} \le \SL_5(\ring).\]
We denote by $\Chev(A_4, \ring)$ the group generated by these root groups, which is a (simply connected) Chevalley group of type $ A_4 $. Then $(\rootgrAex{\alpha})_{\alpha \in A_4}$ is an $ A_4 $-grading of $\Chev(A_4, \ring)$.

Denote the ten roots in $ H_2 $ by $ \htworoot{1}, \ldots, \htworoot{10} $ where $ (\htworoot{1}, \htworoot{2}, \htworoot{3}, \htworoot{4}, \htworoot{5}) $ is an (arbitrarily chosen) $ H_2 $-quintuple and $ \htworoot{i+5} \defl -\htworoot{i} $ for all $ i \in \Set{1, \ldots, 5} $. We define root groups $ \rootgrHex{\htworoot{i}} \defl \rootgrHex{i} $ for all $ i \in \Set{1, \ldots, 10} $ by
	\begin{align*}
		&\rootgrHex{1} \coloneqq \rootgrAex{12}\rootgrAex{34}, \: \rootgrHex{2} \coloneqq \rootgrAex{35}\rootgrAex{14},\: \rootgrHex{3} \coloneqq \rootgrAex{24}\rootgrAex{15}, \: \rootgrHex{4} \coloneqq \rootgrAex{13}\rootgrAex{25},\: \rootgrHex{5} \coloneqq \rootgrAex{45}\rootgrAex{23}, \\
		&\rootgrHex{6} \coloneqq \rootgrAex{21}\rootgrAex{43},\: \rootgrHex{7} \coloneqq \rootgrAex{53}\rootgrAex{41},\:	\rootgrHex{8} \coloneqq \rootgrAex{42}\rootgrAex{51}, \: \rootgrHex{9} \coloneqq \rootgrAex{31}\rootgrAex{52},\: \rootgrHex{10} \coloneqq \rootgrAex{54}\rootgrAex{32}.
	\end{align*}
By a straightforward computation, which we leave out because this example serves only motivational purposes, $ (\rootgrHex{\htworootsym})_{\htworootsym \in H_2} $ is an $ H_2 $-grading of $ \El_5(\ring) $. We call $ (\rootgrHex{\htworootsym})_{\htworootsym \in H_2} $ a \emph{folding of $(\rootgrAex{\alpha})_{\alpha \in A_4}$}. Note that the surrounding group $ \El_5(\ring) $ is the same for both gradings: only the family of root groups changes. Note further that for each $ i \in \{1, \ldots, 10\} $, the two root groups in the definition of $ \rootgrHex{i} $ commute, so that $ \rootgrHex{i} $ is isomorphic to $ \ring \times \ring $

To pave the way for a generalisation, we rephrase our construction in a more abstract way. Denote by
\[ \delta_2 \defl e_1 - e_2, \quad \delta_3 \defl e_2 - e_3, \quad \delta_4 - e_3-e_4, \quad \delta_5 \defl e_4-e_5 \]
a root base of $A_4$ and
define a linear map
\[ \map{\bar{\foldproj}}{\gen{A_4}_\R}{\gen{H_2}_\R}{}{}, \delta_2, \delta_4 \mapsto \htworoot{1}, \delta_3, \delta_6 \mapsto \htworoot{5}. \]
Then a straightforward computation (which, again, we leave out) shows that $\bar{\foldproj}$ induces a surjection $\map{\foldproj}{A_4}{H_2}{}{}$ with fibres of cardinality~$2$, and we have
\[ \rootgrHex{\beta} = \gen{\rootgrAex{\alpha} \mid \alpha \in \foldproj^{-1}(\beta)} \]
for all $\alpha \in H_2 $. We encode this folding by the partition $ \Set{\Set{\delta_2, \delta_4}, \Set{\delta_3, \delta_5}} $ of the root base $\Set{\delta_2, \delta_3, \delta_4, \delta_5}$. We also say that $ (\rootgrHex{\htworootsym})_{\htworootsym \in H_2} $ is \emph{the} folding of $(\rootgrAex{\alpha})_{\alpha \in A_4}$ described by the diagram in Figure~\ref{fig:fold-diags:H2}.
\end{example}

\begin{notation}
	The tilde on the root groups in~\ref{ex:a4h2folding} represents the fact that $(\rootgrAex{\alpha})_{\alpha \in A_4}$ is a \emph{specific example} of an $A_4$-grading and not an \emph{arbitrary} $A_4$-grading. We will use a similar notation in the later Examples~\ref{ex:chev-D6} and~\ref{ex:chev-E8}.
\end{notation}

\begin{rem}
	Everything in~\ref{ex:a4h2folding} remains correct if the ring $ \ring $ is not assumed to be commutative, except that the resulting group is then no longer called a Chevalley group, and that there is no longer a standard notion of $ \SL_5(\ring) $.
\end{rem}

\begin{rem}\label{rem:fold-A4-cox}
	At the level of Coxeter groups, $ (\gen{s_1, s_2}, \Set{s_1, s_2}) $ is a Coxeter system of type $ H_2 $ where $ s_1 \defl \refl{e_1 - e_2} \refl{e_3 - e_4} $ and $ s_2 \defl \refl{e_2 - e_3} \refl{e_4 - e_5} $. This realisation of $ \Weyl(H_2) $ as a subgroup of $ \Weyl(A_4) $ is a special case of the construction in \cite{coxcox}: Given a Coxeter system $ (W,S) $ and a partition $ I = \Set{I_1, \ldots, I_k} $ of $ S $ which has the so-called property of being \emph{admissible}, the pair $ (\gen{R}, R) $ is a Coxeter system where $ R = \Set{r_1, \ldots, r_k} $ and $ r_j $ is the longest element in the Coxeter system $ (\gen{I_j}, I_j) $ for all $ j \in \Set{1, \ldots, k} $. In all cases we consider, each set $ I_j $ has cardinality $ 2 $ and its elements $ s_j $, $ t_j $ commute (that is, $ (\gen{I_j}, I_j) $ is of type $ A_1 \times A_1 $), so that $ r_j = s_j t_j $.
\end{rem}

\begin{rem}\label{rem:fold-H-weyl}
	The same procedure as in~\ref{rem:fold-A4-cox} allows us to realise $\Weyl(H_3)$ and $\Weyl(H_4)$ as subgroups of $\Weyl(D_6)$ and $\Weyl(E_8)$, respectively. Namely, let $\rootbaseE = (\delta_1, \ldots, \delta_8)$ be an ordered root base of $E_8$ as in \ref{not:root-base-order}. Then by \cite[1.1, 5.4]{coxcox}, $ (\gen{s_0, s_1, s_2, s_3}, R) $ is a Coxeter system of type $ H_4 $ where $ s_0 \defl \refl{\delta_7} \refl{\delta_8} $, $ s_1 \defl \refl{\delta_1} \refl{\delta_6} $, $ s_2 \defl \refl{\delta_2} \refl{\delta_4} $, $ s_3 \defl \refl{\delta_3} \refl{\delta_5} $ and $ R \defl \Set{s_0, s_1, s_2, s_3} $. Similarly, $ (\gen{s_1, s_2, s_3}, \Set{s_1, s_2, s_3}) $ is a Coxeter system of type $ H_3 $. These foldings are encoded by the diagrams in Figures~\ref{fig:fold-diags:H3} and~\ref{fig:fold-diags:H4}.
\end{rem}

Remarks~\ref{rem:fold-A4-cox} and~\ref{rem:fold-H-weyl} describe the folding $\map{}{X}{H_\ell}{}{}$ at the level of Coxeter groups -- that is, as an embedding $\map{}{\Weyl(H_\ell)}{\Weyl(X)}{}{}$. To obtain a map $\map{\foldproj}{X}{H_\ell}{}{}$, we will first project onto a suitable subspace and then scale to the unit circle. The first step yields a map $ \map{\goldfoldproj}{X}{GH_\ell}{}{} $, which is described in the following proposition.

\begin{prop}\label{pro:fold-abstract-description}
	The following hold:
	\begin{enumerate}[(a)]
		\item \label{pro:fold-abstract-description:char}There exists a unique bijection $\map{\goldfoldproj}{E_8}{GH_4}{}{}$ with
		\begin{align*}
			\goldfoldproj(\delta_7) &= \rho_0, & \goldfoldproj(\delta_1) &= \rho_1, & \goldfoldproj(\delta_2) &= \rho_2, & \goldfoldproj(\delta_3) &= \gold \rho_3, \\*
			\goldfoldproj(\delta_8) &= \gold \rho_0, & \goldfoldproj(\delta_6) &= \gold \rho_1, & \goldfoldproj(\delta_4) &= \gold \rho_2, & \goldfoldproj(\delta_5) &= \rho_3
		\end{align*}
		and such that $\goldfoldproj$ is induced by a linear map $\map{}{\gen{E_8}_\R}{\gen{GH_4}_\R}{}{}$ on the surrounding vector spaces.
		
		\item \label{pro:fold-abstract-description:H3}The map $\goldfoldproj$ satisfies $\goldfoldproj(D_6) = GH_3$ and $\goldfoldproj(A_4) = GH_2$. Here $D_6$, $A_4$ and $GH_3$, $GH_2$ are embedded into $E_8$, $GH_4$, respectively, as in~\ref{not:root-base-order}.
	
		\item \label{pro:fold-abstract-description:refl}We have
		\[ \goldfoldproj(\alpha^{\reflbr{\delta_i} \reflbr{\delta_j}}) = \goldfoldproj(\alpha)^{\reflbr{\goldfoldproj(\delta_i)}} = \goldfoldproj(\alpha)^{\reflbr{\goldfoldproj(\delta_j)}} \]
		for all $\alpha \in E_8$ and all $(i,j) \in \Set{(7,8), (1,6), (2,4), (5,3)}$.
		
		\item \label{pro:fold-abstract-description:weyl-hom}For each $w \in \Weyl(H_4)$, there exists a unique $\weylemb(w) \in \Weyl(E_8)$ with $\goldfoldproj(\alpha^{\weylemb(w)}) = \goldfoldproj(\alpha)^w$ for all $\alpha \in E_8$. This defines an injective homomorphism $\map{}{\Weyl(H_4)}{\Weyl(E_8)}{w}{\weylemb(w)}$ which maps $\Weyl(H_3)$ into $\Weyl(D_6)$.
		
		\item \label{pro:fold-abstract-description:foldproj}The map $\map{\foldproj \defl \scalmap \circ \goldfoldproj}{E_8}{H_4}{}{}$ (where $\map{\scalmap}{GH_4}{H_4}{}{}$ denotes the scaling map from~\ref{def:GH}) is a surjection with $\abs{\foldproj^{-1}(\beta)} = 2$ for all $\beta \in H_4$ and $\goldfoldproj(\alpha^{\weylemb(w)}) = \goldfoldproj(\alpha)^w$ for all $\alpha \in E_8$ and $w \in \Weyl(H_4)$.
	\end{enumerate}

\end{prop}
\begin{proof}
	Put $W\coloneqq\langle\delta_1,\ldots,\delta_8\rangle_\R$ and define subspaces
		\begin{align*}
			W_1 &\coloneqq \langle \delta_8 - \gold\delta_7,\delta_6-\tau\delta_1,\delta_4-\tau\delta_2, \delta_3-\tau\delta_5 \rangle_\R, \\
			W_2 &\coloneqq \langle \delta_7 + \gold\delta_8,\delta_1+\tau\delta_6,\delta_2+\tau\delta_4, \delta_5+\tau\delta_3 \rangle_\R
		\end{align*}
	where $\tau\coloneqq\frac{1+\sqrt{5}}{2}$ denotes the golden ratio. Denote by $\goldfoldproj:W\rightarrow W_2$ the orthogonal projection onto $W_2$ and define $ R $ as in~\ref{rem:fold-H-weyl}. A straightforward computation (see~\ref{rem:gap}) shows that $W=W_1\oplus W_2$ is an orthogonal, $R$-invariant decomposition of $ W $, that $ \goldfoldproj $ is injective on $ E_8 $ (so that $ \abs{\goldfoldproj(E_8)} = 240 $) and that $ \goldfoldproj(E_8) $ is a non-reduced root system in $ W_2 $. More precisely, we have $ \goldfoldproj(E_8) = Y \cup \tau Y $ where $ Y $ denotes the set of indivisible roots in $ \goldfoldproj(E_8) $. We conclude that $ Y $ is of type $ H_4 $ because this is the only reduced root system of rank at most $ 4 $ with $ 120 $ roots. Hence $ \goldfoldproj(E_8) $ is (up to scaling) a root system of type $ GH_4 $. Another computation shows that $(\goldfoldproj(\delta_7), \goldfoldproj(\delta_1), \goldfoldproj(\delta_2), \goldfoldproj(\delta_5))$ is a root base of $ \goldfoldproj(E_8) $ (in the sense of~\ref{def:GH}) and that $\goldfoldproj(\delta_j) = \gold \goldfoldproj(\delta_i)$ for all $(i,j) \in \Set{(7,8), (1,6), (2,4), (5,3)}$. By identifying $ \goldfoldproj(E_8) $ with $GH_4$ via the root bases $(\goldfoldproj(\delta_7), \goldfoldproj(\delta_1), \goldfoldproj(\delta_2), \goldfoldproj(\delta_5))$ and $(\rho_0, \rho_1, \rho_2, \rho_3)$, we obtain a map $\map{}{E_8}{GH_4}{}{}$ which is induced by a linear map and which acts on $(\delta_1, \ldots, \delta_8)$ in the desired way. Since a linear map is determined by its action on a basis, $\goldfoldproj$ is necessarily unique. This proves~\ref{pro:fold-abstract-description:char}.
	
	Since $\goldfoldproj(\Set{\delta_1, \ldots, \delta_6}) \subseteq GH_3$ and $\goldfoldproj$ is (induced by a) linear map, we have $\goldfoldproj(D_6) \subseteq GH_3$. Now let $\beta \in GH_3$. Since $\goldfoldproj$ is bijective, there exists $\alpha \in E_8$ with $\goldfoldproj(\alpha) = \beta$. Since $(\delta_1, \ldots, \delta_8)$ is a root base and $\goldfoldproj$ is linear, we may without loss of generality assume that $\alpha$ is positive and write $\alpha = \sum_{i=1}^8 a_i \delta_i$ with $a_1, \ldots, a_8 \in \Z_{>0}$. Then
	\begin{align*}
		\beta = \goldfoldproj(\alpha) &= (a_7 + \gold a_8) \rho_0 + (a_1 + \gold a_6) \rho_1 + (a_2 + \gold a_4) \rho_2 + (a_5 + \gold a_3) \rho_3.
	\end{align*}
	Since $\beta \in GH_3$, we conclude that $a_7 = -\gold a_8$. This implies that $a_7 = a_8 = 0$. Thus $\alpha \in D_6$, so $\goldfoldproj(D_6) = GH_3$. In a similar way, one can show $\goldfoldproj(A_4) = GH_2$. This proves~\ref{pro:fold-abstract-description:H3}.
	
	For~\ref{pro:fold-abstract-description:refl}, let $(i,j) \in \Set{(7,8), (1,6), (2,4), (5,3)}$, let $\alpha \in E_8$ and put $r \defl \reflbr{\delta_i} \reflbr{\delta_j}$. Note that the map $\goldfoldproj$ in our explicit construction is an orthogonal projection, so that $\goldfoldproj(\alpha^r) = \goldfoldproj(\alpha)^r$. It remains to check that $r$, $\reflbr{\goldfoldproj(\delta_i)}$ and $\reflbr{\goldfoldproj(\delta_j)}$ act identically on any root in $GH_4$, which is a straightforward computation (see, again,~\ref{rem:gap}).
	
	For~\ref{pro:fold-abstract-description:weyl-hom}, let $w \in H_4$. The existence of $\weylemb(w)$ follows from~\ref{pro:fold-abstract-description:refl} because $\Weyl(H_4)$ is generated by $\Set{\reflbr{\rho_i} \mid i \in \Set{0,1,2,3}} = \Set{\reflbr{\goldfoldproj(\delta_j)} \mid j \in \Set{1,\ldots, 8}}$. Since $\goldfoldproj$ is bijective, $\weylemb(w)$ is uniquely determined by the desired property and $\map{}{}{}{w}{\weylemb(w)}$ is an injective homomorphism. Further, our construction shows that $\weylemb(w) \in \Weyl(D_6)$ if $w \in \Weyl(H_3)$.
	
	Finally, assertion~\ref{pro:fold-abstract-description:foldproj} follows from the fact that $\scalmap$ is surjective with fibres of cardinality~2 and that $\scalmap(\alpha^w) = \scalmap(\alpha)^w$ for all $\alpha \in GH_4$ and $w \in \Weyl(H_4)$.
\end{proof}

\begin{note}
	The ordering of $\delta_1, \ldots, \delta_8$ in~\ref{pro:fold-abstract-description}~\ref{pro:fold-abstract-description:char} corresponds to the way the diagram in Figure~\ref{fig:fold-diags:H4} is drawn. Further, the monomorphism $\map{}{}{}{w}{\weylemb(w)}$ in~\ref{pro:fold-abstract-description}~\ref{pro:fold-abstract-description:weyl-hom} is precisely the embedding of Weyl groups in~\ref{rem:fold-H-weyl}. Note also that the choice of $\goldfoldproj$ depends on the fixed choices of the root bases $(\rho_0, \rho_1, \rho_2, \rho_3)$ and $(\delta_1, \ldots, \delta_8)$.
\end{note}

\begin{note}
	The generators of $W_1$ in the proof of~\ref{pro:fold-abstract-description} encode precisely the relations $\delta_8 \equiv \gold \delta_7$, $\delta_6 \equiv \gold \delta_1$, $\delta_4 \equiv \gold \delta_2$ and $\delta_3 \equiv \gold \delta_5$ that we want to realise modulo $\goldfoldproj$. Further, note that the orbit of $ (\delta_1, \delta_6) $ under $ \weylemb(\Weyl(H_4)) $ contains $ (\delta_2, \delta_4) $, $ (\delta_5, \delta_3) $ and $ (\delta_7, \delta_8) $, as can be seen from the following diagram:
	\[ \begin{tikzcd}[column sep=large]
		{(\delta_7, \delta_8)} \arrow[r, mapsto, "r_1 r_0"] & {(\delta_1, \delta_6)} \arrow[r, mapsto, "r_2 r_1"] & {(\delta_2, \delta_4)} \arrow[r, mapsto, "r_2 r_3 r_2 r_3"] & {(\delta_5, \delta_3)}.
	\end{tikzcd} \]
	Together with the observation from~\ref{pro:fold-abstract-description} that $ W_1 $ is $ R $-invariant, it follows that $ W_1 $ is actually the smallest $ R $-invariant subspace of $ W $ containing any of its four generators. In particular, it is not possible to choose generators of $ W_1 $ so that $ \goldfoldproj(\delta_2) $ and $ \goldfoldproj(\delta_3) $ have the same length in $ GH_4 $.
\end{note}

\begin{rem}\label{rem:goldfoldproj-restrict}
	We can also obtain bijections $\map{}{D_6}{GH_3}{}{}$ and $\map{}{A_4}{GH_2}{}{}$ as well as surjections $\map{}{D_6}{H_3}{}{}$ and $\map{}{A_4}{H_2}{}{}$ with similar properties by restricting $\map{\goldfoldproj}{E_8}{GH_4}{}{}$ and $\map{\foldproj}{E_8}{H_4}{}{}$ to the subsystems $\subsys{\delta_1, \ldots, \delta_6}$ and $\subsys{\delta_2, \ldots, \delta_5}$ of $E_8$.
\end{rem}

\begin{rem}\label{rem:H3-preimages}
	Figures~\ref{fig:foldrootse8-1} and~\ref{fig:foldrootse8-2} in the appendix depict all tuples $ (\beta, \alpha_1, \alpha_2) \in H_4 \times E_8 \times E_8 $ with $\beta$ positive, $\scalmap\goldfoldproj(\alpha_1) = \beta = \scalmap\goldfoldproj(\alpha_2)$ and $\goldfoldproj(\alpha_2) = \gold \goldfoldproj(\alpha_1)$. (Note that then $(-\beta, -\alpha_1, -\alpha_2)$ is a tuple with the same properties except that $\beta$ is negative.) Here we use Notation~\ref{def:rootcoord} to depict roots via their coordinates with respect to the root bases $ (\rho_0, \rho_1, \rho_2, \rho_3) $ and $ (\delta_1, \ldots, \delta_8) $. For the reader's convenience, the subtable consisting of all tuples $(\beta, \alpha_1, \alpha_2) \in H_3 \times D_6 \times D_6$ is given in Figure~\ref{fig:foldrootsd6}.
\end{rem}

\premidfigure
\begin{figure}[htb]
	\centering\begin{tabular}[t]{ccc}
		\toprule
		$ \beta $ & $ \alpha_1 $ & $ \alpha_2 $ \\
		\midrule
		$ \hcoord{1}{0}{0} $ & $ e_1 - e_2 $ & $ e_5 + e_6 $ \\
		$ \hcoord{0}{1}{0} $ & $ e_2 - e_3 $ & $ e_4 - e_5 $ \\
		$ \hcoord{0}{0}{1} $ & $ e_5 - e_6 $ & $ e_3 - e_4 $ \\
		$ \hcoord{1}{1}{0} $ & $ e_1 - e_3 $ & $ e_4 + e_6 $ \\
		$ \hcoord{1}{1}{\gold} $ & $ e_1 - e_4 $ & $ e_3 + e_5 $ \\
		$ \hcoord{1}{\gold^2}{\gold} $ & $ e_1 - e_5 $ & $ e_2 + e_4 $ \\
		$ \hcoord{1}{\gold^2}{\gold^2} $ & $ e_1 - e_6 $ & $ e_2 + e_3 $ \\
		$ \hcoord{0}{1}{\gold} $ & $ e_2 - e_4 $ & $ e_3 - e_6 $ \\
		\bottomrule
	\end{tabular}\quad\begin{tabular}[t]{ccc}
		\toprule
		$ \beta $ & $ \alpha_1 $ & $ \alpha_2 $ \\
		\midrule
		$ \hcoord{\gold}{\gold^2}{\gold} $ & $ e_2 + e_6 $ & $ e_1 + e_4 $ \\
		$ \hcoord{\gold}{2\gold}{\gold^2} $ & $ e_3 + e_4 $ & $ e_1 + e_2 $ \\
		$ \hcoord{0}{\gold}{\gold} $ & $ e_3 - e_5 $ & $ e_2 - e_6 $ \\
		$ \hcoord{\gold}{\gold}{\gold} $ & $ e_3 + e_6 $ & $ e_1 + e_5 $ \\
		$ \hcoord{\gold}{\gold}{1} $ & $ e_4 + e_5 $ & $ e_1 + e_6 $ \\
		$ \hcoord{0}{\gold}{1} $ & $ e_4 - e_6 $ & $ e_2 - e_5 $ \\
		$ \hcoord{\gold}{\gold^2}{\gold^2} $ & $ e_2 + e_5 $ & $ e_1 + e_3 $ \\
		\bottomrule
	\end{tabular}
	\caption{The preimages of the positive roots in $ H_3 $ under $ \foldproj $. See~\ref{rem:H3-preimages}.}
	\label{fig:foldrootsd6}
\end{figure}
\postmidfigure

\begin{notation}\label{not:fold-goldproj}
	For the rest of this section and in addition to~\ref{not:root-base-order}, we fix the bijection $\map{\goldfoldproj}{E_8}{GH_4}{}{}$ and the surjection $\map{\foldproj}{E_8}{H_4}{}{}$ from~\ref{pro:fold-abstract-description}. If $X \in \Set{A_4, D_6}$, we will always consider the restrictions of $\foldproj$ and $\goldfoldproj$ to $ X $ in place of $\foldproj$ and $\goldfoldproj$, as described in~\ref{rem:goldfoldproj-restrict}.
\end{notation}

Before we can use~\ref{pro:fold-abstract-description} to construct $H_\ell$-graded groups, we collect some more useful properties of $\foldproj$.

\begin{lemma}\label{lem:preim-trans}
	The Weyl group of $ H_\ell $, embedded into $\Weyl(X)$ via~\ref{pro:fold-abstract-description}~\ref{pro:fold-abstract-description:weyl-hom}, acts transitively on the set $\Set{\foldproj^{-1}(\beta) \mid \beta \in H_\ell}$. More precisely, $\foldproj^{-1}(\beta^w) = \foldproj^{-1}(\beta)^{\weylemb(w)}$ for all $w \in \Weyl(H_\ell)$.
\end{lemma}
\begin{proof}
	Let $\beta \in  H_\ell$, $w \in \Weyl(H_\ell)$ be arbitrary. Then $\foldproj^{-1}(\beta)^{\weylemb(w)} \subseteq \foldproj^{-1}(\beta^w)$ by~\ref{pro:fold-abstract-description}~\ref{pro:fold-abstract-description:foldproj}. As both sets have cardinality $2$, we conclude that $\foldproj^{-1}(\beta)^{\weylemb(w)} = \foldproj^{-1}(\beta^w)$. It follows that $\Weyl(H_\ell)$ acts on $\Set{\foldproj^{-1}(\beta) \mid \beta \in H_\ell}$, and this action is transitive by~\ref{lem:WH3transitive}.
\end{proof}

\begin{lemma}\label{lem:preim-descr}
	For all $\beta \in H_\ell$, the preimage $\foldproj^{-1}(\beta)$ is of the form $\Set{\alpha_1, \alpha_2}$ for unique orthogonal roots $\alpha_1, \alpha_2 \in X$ with $\goldfoldproj(\alpha_2) = \gold \goldfoldproj(\alpha_1)$.
\end{lemma}
\begin{proof}
	Let $\beta \in H_\ell$. Since $\goldfoldproj$ is bijective and $\foldproj = \scalmap \circ \goldfoldproj$, we have $\foldproj^{-1}(\beta)=\Set{\alpha_1, \alpha_2}$ for unique distinct roots $\alpha_1, \alpha_2 \in X$ with $\goldfoldproj(\alpha_2) = \gold \goldfoldproj(\alpha_1)$. If $\beta \in \Set{\rho_1, \rho_2, \rho_3}$ (or $\beta = \rho_0$ if $\ell = 4$), then $\alpha_1$, $\alpha_2$ are given by ~\ref{pro:fold-abstract-description}~\ref{pro:fold-abstract-description:char}, and we can observe that they are orthogonal. Now it follows from~\ref{lem:preim-trans} that $\alpha_1, \alpha_2$ are orthogonal for arbitrary $\beta \in H_\ell$.
\end{proof}

\begin{rem}\label{note:preim-order}
	In~\ref{lem:preim-descr}, we implicitly introduce an order on $ \foldproj^{-1}(\beta) $ by declaring that $ \alpha_1 $ is the root in $ \foldproj^{-1}(\beta) $ whose image in $ GH_\ell $ is short. It is compatible with the action of $ \weylemb(\Weyl(H_4)) $: If $ \foldproj(\alpha) $ is short in $ GH_3 $, then so is $ \foldproj(\alpha^{\weylemb(w)}) = \foldproj(\alpha)^w $ (for $ \alpha \in E_8 $ and $ w \in \Weyl(H_4) $).
	
	We will use this order in~\ref{def:param-folding} to decide which factor of a root group will be the \enquote{left one} and which one the \enquote{right one}. For the computations in~\ref{rem:fold-ex-parmap} (specifically, in~\eqref{eq:fold-ex-parmap}), it will be crucial that this choice is performed in a way that is compatible with the action of the Weyl elements.
\end{rem}

\begin{lemma}\label{lem:weylemb-on-refl}
	Let $\beta \in H_4$ and let $\alpha_1, \alpha_2$ denote the two roots in $\foldproj^{-1}(\beta)$. Then $\weylemb(\refl{\beta}) = \refl{\alpha_1} \refl{\alpha_2} = \refl{\alpha_2} \refl{\alpha_1}$.
\end{lemma}
\begin{proof}
	If $\beta \in \Set{\rho_0, \ldots, \rho_3}$, the assertion holds by~\ref{pro:fold-abstract-description}~\ref{pro:fold-abstract-description:refl}. For arbitrary $\beta$, there exists $w \in \Weyl(H_4)$ with $\beta = \rho_1^w$ by~\ref{lem:WH3transitive}. Hence, by~\ref{rem:refl-right}, $\refl{\beta} = \refl{\rho_1}^w$ and
	\[ \weylemb(\refl{\beta}) = \weylemb(\reflbr{\rho_1})^{\weylemb(w)} = (\refl{\delta_1} \refl{\delta_6})^{\weylemb(w)} = \refl{\delta_1^{\weylemb(w)}} \refl{\delta_6^{\weylemb(w)}}. \]
	Since
	\[ \Set{\alpha_1, \alpha_2} = \foldproj^{-1}(\beta) = \foldproj^{-1}(\rho_1^w) = \foldproj^{-1}(\rho_1)^{\weylemb(w)} = \Set{\delta_1, \delta_6}^{\weylemb(w)} \]
	by~\ref{lem:preim-trans}, we infer that $\weylemb(\refl{\beta}) = \refl{\alpha_1} \refl{\alpha_2}$.
\end{proof}

We can now turn to the folding procedure for root graded groups.

\begin{notation}
	For the rest of this section, we denote by $G$ an arbitrary group and by $ (\rootgrX{\alpha})_{\alpha \in X} $ an $ X $-grading of $G$.
\end{notation}

\begin{definition}\label{def:rootgr-fold}
	The \emph{folding of $ (\rootgrX{\alpha})_{\alpha \in X} $ (with respect to the partition given by Figure~\ref{fig:fold-diags})} is the family $ (\rootgrH{\beta})_{\beta \in H_\ell} $ defined by
	\[ \rootgrH{\beta} \defl \rootgrX{\foldproj^{-1}(\Set{\beta})} = \gen{\rootgrX{\alpha} \mid \alpha \in X, \foldproj(\alpha) = \beta} \]
	for all $ \beta\in H_\ell $.
\end{definition}

In~\ref{folding-is-rgg}, we will show that the folding $ (\rootgrH{\beta})_{\beta \in H_\ell} $ is an $H_\ell$-grading of $G$.

\begin{rem}\label{rem:fold-rootgr-commute}
	For any root $ \beta \in H_\ell $, the two roots $ \alpha_1, \alpha_2 $ in $ \foldproj^{-1}(\beta) $ are orthogonal by~\ref{lem:preim-descr} and hence adjacent. Hence $ \rootgrH{\beta} = \rootgrX{\alpha_1} \rootgrX{\alpha_2} = \rootgrX{\alpha_2} \rootgrX{\alpha_1} $, and this group is abelian. In particular, the construction in~\ref{ex:a4h2folding} is a special case of~\ref{def:rootgr-fold}.
\end{rem}

In the examples of Chevalley groups, each root group is parametrised by the additive group of a commutative ring $\ring$. This leads to a natural parametrisation of the root groups in the folding by $\ring \times \ring$.

\begin{definition}\label{def:param-folding}
	Let $ (\ring, +) $ be a group and let $ (\map{\risomX{\alpha}}{\ring}{\rootgrX{\alpha}}{}{})_{\alpha \in X} $ be a family of root isomorphisms (in the sense of~\ref{def:risom}). The \emph{folding of $ (\map{\risomX{\alpha}}{\ring}{\rootgrX{\alpha}}{}{})_{\alpha \in X} $ (with respect to the partition given by Figure~\ref{fig:fold-diags})} is the family $ (\map{\risomH{\beta}}{\ring \times \ring}{\rootgrH{\beta}}{}{})_{\beta \in H_\ell} $ defined by
	\[ \risomH{\beta}(a,b) \defl \risomX{\alpha_1}(a) \risomX{\alpha_2}(b) \]
	for all $ \beta \in H_\ell $ where $ (\alpha_1, \alpha_2) $ is the unique pair of roots in $ X $ with $ \foldproj(\alpha_1) = \beta = \foldproj(\alpha_2) $ and $ \goldfoldproj(\alpha_2) = \gold \goldfoldproj(\alpha_1) $.
\end{definition}

Note that the group $(\ring, +)$ in~\ref{def:param-folding} is necessarily abelian by~\ref{rem:rootgr-abel}.

In the following, we give two more examples of $X$-graded groups, which give rise to foldings in the sense of~\ref{def:rootgr-fold} and~\ref{def:param-folding}. They are equipped with families $(\risomXnaive{\alpha})_{\alpha \in X}$ of root isomorphisms where the letter \enquote{n} stands for \enquote{naive}. In~\ref{rem:ex-twist-signs}, we will define \emph{twists} of these families which ensure that certain structure constants (or rather, signs) behave in a specific way.

\begin{example}\label{ex:chev-D6}
	Let $\ring$ be a commutative ring and denote by $e_{ij}\in \ring^{6\times 6}$ the matrix with $1_\ring$ at position $(i,j)$ and $0_\ring$ at every other position. For each root $ \alpha \in D_6 $, we define a root homomorphism $ \map{\risomDexnaive{\alpha}}{(\ring,+)}{\GL_{12}(\ring)}{}{} $ by the formulas
	\[ \risomDexnaive{e_i - e_j}(\lambda) \defl \id + \lambda \begin{pmatrix}
		e_{ij} & 0 \\
		0 & -e_{ji}
	\end{pmatrix} \]
	for all $ i \ne j \in \Set{1, \ldots, 6} $ and
	\[ \risomDexnaive{e_i + e_j}(\lambda) \defl \id + \lambda \begin{pmatrix}
		0 & e_{ij}-e_{ji} \\
		0 & 0
	\end{pmatrix}, \quad \risomDexnaive{-e_i - e_j}(\lambda) \defl \id + \lambda \begin{pmatrix}
		0 & 0 \\
		e_{ji}-e_{ij} & 0
	\end{pmatrix} \]
	for all $ i<j \in \Set{1, \ldots, 6} $, where $ \lambda \in \ring $. For each $ \beta \in D_6 $, we define a root group $ \rootgrDex{\beta} \defl \risomDexnaive{\beta}(\ring) $ and we denote by $ \Chev(D_6, \ring) $ the group generated by all root groups. This is a (simply connected) Chevalley group of type $ D_6 $, and hence a $ D_6 $-graded group. It is equipped with a family $ (\risomDexnaive{\alpha})_{\alpha \in D_6} $ of root isomorphisms.
\end{example}

\begin{example}\label{ex:chev-E8}
	Let $\ring$ be a commutative ring, let $ L $ be a complex Lie algebra of type $ E_8 $ and let $ (x_\alpha)_{\alpha \in E_8} \cup (h_i)_{1 \le i \le 8} $ be a Chevalley basis of $ L $ in the sense of \cite[25.2]{hum}. Then by~\cite[25.5]{hum}, $ L(\Z) \defl \gen{x_\alpha, h_i \mid \alpha \in E_8, 1 \le i \le 8}_\Z $ is invariant under $ \ad_{x_\alpha}^n / n! $ for all $ \alpha \in E_8 $ and $ n \in \N_0 $. Further, $ \ad_{x_\alpha} $ is nilpotent for all $ \alpha \in E_8 $. Hence for all $ \alpha \in E_8 $, we may define a homomorphism
	\[ \map{\risomEexnaive{\alpha}}{(\ring, +)}{\GL_\ring\brackets[\big]{L(\Z) \otimes \ring}}{\lambda}{\sum_{n=0}^\infty \frac{\ad_{x_\alpha}^n}{n!} \otimes \lambda^n \id_\ring} \]
	and a corresponding root group $ \rootgrEex{\alpha} \defl \risomEexnaive{\alpha}(\ring) $, which acts on $L(\Z) \otimes \ring$ on the left. (For details, see~\cite[Section~3]{Vavilov-Plotkin}.) The group generated by $ (\rootgrEex{\alpha})_{\alpha \in E_8} $ is an (adjoint) Chevalley group of type $ E_8 $, which we denote by $ \Chev(E_8, \ring) $. It is equipped with an $ E_8 $-grading $ (\rootgrEex{\alpha})_{\alpha \in E_8} $ and a family $ (\risomEexnaive{\alpha})_{\alpha \in E_8} $ of root isomorphisms. Note that it is a matrix group of degree $\dim_\C L = \abs{E_8} + \rank(E_8) = 248$.
\end{example}

We now show that $(\rootgrH{\beta})_{\beta \in H_\ell}$ is an $H_\ell$-grading of $G$. At first, we prove that there exist Weyl elements with respect to $(\rootgrH{\beta})_{\beta \in H_\ell}$.

\begin{lemma}\label{lem:fold-weyl}
	Let $ \beta \in H_\ell $ and write $\foldproj^{-1}(\beta) = \Set{\alpha_1, \alpha_2}$ where $ \alpha_1, \alpha_2 \in X $. For each $ i \in \Set{1,2} $, let $ w_{\alpha_i} $ be an $ \alpha_i $-Weyl element in $ (G, (\rootgrX{\gamma})_{\gamma \in X}) $. Then $ w_\beta \defl w_{\alpha_1} w_{\alpha_2} = w_{\alpha_2} w_{\alpha_1} $ is a $ \beta $-Weyl element in $ (G, (\rootgrH{\delta})_{\delta \in H_\ell}) $.
\end{lemma}
\begin{proof}
	By (the arguments in)~\ref{rem:fold-rootgr-commute}, $\gen{\rootgrX{\alpha_1}, \rootgrX{-\alpha_1}}$ commutes with $\gen{\rootgrX{\alpha_2}, \rootgrX{-\alpha_2}}$. Hence
	\[ w_\beta \in \rootgrX{-\alpha_1} \rootgrX{-\alpha_2} \rootgrX{\alpha_1} \rootgrX{\alpha_2} \rootgrX{-\alpha_1} \rootgrX{-\alpha_2} = \rootgrH{-\beta} \rootgrH{\beta} \rootgrH{-\beta}. \]
	Now let $\delta \in H_\ell$ be arbitrary and write $\foldproj^{-1}(\delta) = \Set{\gamma_1, \gamma_2}$ for $\gamma_1, \gamma_2 \in X$. Then
	\begin{align*}
		(\rootgrH{\delta})^{w_\beta} &= (\rootgrX{\gamma_1} \rootgrX{\gamma_2})^{w_{\alpha_1} w_{\alpha_2}} = \rootgrX{\gamma_1^{\reflbr{\alpha_1} \reflbr{\alpha_2}}} \rootgrX{\gamma_2^{\reflbr{\alpha_1} \reflbr{\alpha_2}}} = \rootgrX{\gamma_1^{\weylemb(\refl{\beta})}} \rootgrX{\gamma_2^{\weylemb(\refl{\beta})}}
	\end{align*}
	by \ref{lem:weylemb-on-refl}. Since $\Set{\gamma_1^{\weylemb(\refl{\beta})}, \gamma_2^{\weylemb(\refl{\beta})}} = \foldproj^{-1}(\delta)^{\weylemb(\refl{\beta})} = \foldproj^{-1}(\delta^{\refl{\beta}})$ by~\ref{lem:WH3transitive}, we infer that $(\rootgrH{\delta})^{w_\beta} = \rootgrH{\delta^{\reflbr{\beta}}}$. Hence $w_\beta$ is indeed a $\beta$-Weyl element.
\end{proof}

While~\ref{lem:fold-weyl} assures the existence of Weyl elements, we want to go further and compute explicit formulas for the actions of certain Weyl elements in the Chevalley groups considered previously. This boils down to computing certain structure constants, or rather, signs. It is part of our main result~\ref{thm:prod-param} that every $H_\ell$-graded group (for $\ell \ge 3$) can be coordinatised in such a way that the same structural signs as for foldings of Chevalley groups are satisfied. However, the structural signs that we would obtain from the families $(\risomDexnaive{\alpha})_{\alpha \in D_6}$ and $(\risomEexnaive{\alpha})_{\alpha \in E_8}$ in~\ref{ex:chev-D6} and~\ref{ex:chev-E8} do not behave as we would like (see~\ref{rem:ex-twist-signs} for details). Hence we will twist these families in the way described as follows, and then compute the desired signs in~\ref{rem:fold-ex-parmap}.

\begin{definition}\label{def:param-twist}
	Let $ T $ be a subset of a positive system of $ X $, let $ \ring $ be an abelian group and let $ (\map{\risomXnaive{\alpha}}{\ring}{\rootgrX{\alpha}}{}{})_{\alpha \in X} $ be a family of root isomorphisms. The \emph{$ T $-twist of $ (\risomXnaive{\alpha})_{\alpha \in X} $} is the family $ (\risomX{\alpha})_{\alpha \in X} $ defined by
	\[ \map{\risomX{\alpha}}{\ring}{\rootgrX{\alpha}}{\lambda}{\begin{cases}
		\risomXnaive{\alpha}(-\lambda) & \text{if } \alpha \in T \cup (-T), \\
		\risomXnaive{\alpha}(\lambda) & \text{if } \alpha \notin T \cup (-T).
	\end{cases}} \]
\end{definition}

\begin{rem}\label{rem:ex-twist-signs}
	For the group $ \Chev(D_6, \ring) $ in~\ref{ex:chev-D6}, we will choose
	\[ T(D_6) \defl \Set{\delta_6 = e_5+e_6}. \]
	This choice will ensure that the formula for the commutator relation of the pair $(\rootgrHex{\rho_1},\rootgrHex{\rho_2})$ has a certain form. See~\ref{rem:fold-comm-sub} for details.
	
	For the group $ \Chev(E_8, \ring) $ in~\ref{ex:chev-E8}, the situation is more complicated because, even without twisting, different choices of the Chevalley basis result in different structure constants of the Chevalley group. For our computation of explicit signs, we will use the Chevalley basis that is provided by the GAP \cite{GAP} command \texttt{SimpleLieAlgebra("{}E"{}, 8, Rationals)} and the $ T(E_8) $-twist defined by
	\[ T(E_8) \defl \left\lbrace \begin{gathered}
		\rootcoord{0, 0, 0, 1, 1, 1, 0, 0}, \rootcoord{0, 1, 1, 1, 1, 1, 0, 0}, \rootcoord{0, 0, 1, 2, 1, 1, 0, 0}, \\
		\rootcoord{0, 0, 0, 0, 0, 1, 0, 0}, \rootcoord{1, 1, 1, 1, 1, 1, 0, 0}, \rootcoord{0, 0, 1, 1, 1, 1, 0, 0}, \\
		\rootcoord{1, 1, 1, 2, 1, 1, 0, 0}, \rootcoord{0, 1, 1, 2, 1, 1, 0, 0}, \rootcoord{0, 1, 2, 2, 1, 1, 0, 0}, \\
		\rootcoord{1, 1, 2, 2, 1, 1, 0, 0}, \rootcoord{1, 2, 2, 2, 1, 1, 0, 0}
	\end{gathered} \right\rbrace. \]
	Here we use Notation~\ref{def:rootcoord}. This choice ensures that the canonical $D_6$-graded subgroup of $ \Chev(E_8, \ring) $ has the same \enquote{structure constants} as the $D_6$-graded group in~\ref{ex:chev-D6}. For details, see~\ref{rem:fold-ex-parmap} and~\ref{rem:fold-comm-sub}.
	
	In the following, we will write $ (\risomDex{\alpha})_{\alpha \in D_6} $ and $ (\risomEex{\alpha})_{\alpha \in E_8} $ for the twists of the families $ (\risomDexnaive{\alpha})_{\alpha \in D_6} $ and $ (\risomEexnaive{\alpha})_{\alpha \in E_8} $ (from~\ref{ex:chev-D6} and~\ref{ex:chev-E8}) by the sets specified above. Further, we will write $ (\risomHex{\beta})_{\beta \in H_\ell} $ for the corresponding folding in the sense of~\ref{def:param-folding}.
\end{rem}

\begin{rem}\label{rem:fold-ex-parmap}
	Let $ \ring $ be a commutative ring and let $X = E_8$. We consider the group $\Chev(E_8, \ring)$ with the $E_8$-grading $ (\rootgrEex{\alpha})_{\alpha \in E_8} $ from~\ref{ex:chev-E8} and the families $(\risomEex{\alpha})_{\alpha \in E_8}$ and $ (\risomHex{\beta})_{\beta \in H_4} $ of root isomorphisms from~\ref{rem:ex-twist-signs}. It is a well-known fact about Chevalley groups (see \cite[Lemma~20~(b)]{Steinberg-ChevGroups} and \cite[3.3.6]{torben} for details) that for all invertible $r \in \ring$ and all $\alpha \in E_8$, the element
	\[ w_\alpha(r) \defl \risomEex{-\alpha}(-r^{-1}) \risomEex{\alpha}(r) \risomEex{-\alpha}(-r^{-1}) \]
	is an $\alpha$-Weyl element in $ (G, (\rootgrEex{\alpha})_{\alpha \in E_8}) $. Hence it follows from~\ref{lem:fold-weyl} that for all invertible $r,s \in \ring$ and all $\alpha_1, \alpha_2 \in E_8$ with $\goldfoldproj(\alpha_2) = \gold \goldfoldproj(\alpha_1)$ and $\beta \defl \foldproj(\alpha_1) = \foldproj(\alpha_2)$, the element
	\begin{align*}
		w_{\beta}(r,s) \defl& \mathrel{} w_{\alpha_1}(r) w_{\alpha_2}(s) = w_{\alpha_2}(s) w_{\alpha_1}(r) \\*
		=& \mathrel{} \risomHex{-\beta}(-r^{-1}, -s^{-1}) \risomHex{\beta}(r,s) \risomHex{-\beta}(-r^{-1}, -s^{-1})
	\end{align*}
	is a $\beta$-Weyl element in $(G, (\rootgrHex{\gamma})_{\gamma \in H_\ell})$.
	Put $w_\rho \defl w_\rho(1_\ring, 1_\ring)$ for all $\rho \in \Set{\rho_0, \rho_1, \rho_2, \rho_3}$. A straightforward but long matrix computation (see~\ref{rem:gap}) shows that for all $\beta \in H_4$ and all $\rho \in \Set{\rho_0, \rho_1, \rho_2, \rho_3}$, there exist signs $\epsilon_{\beta, \delta}, \bar{\epsilon}_{\beta, \delta} \in \Set{\pm 1}$ such that
	\begin{equation}\label{eq:fold-ex-parmap}
		\risomHex{\beta}(r,s)^{w_\rho} = \risomHex{\beta^{\reflbr{\rho}}}(\epsilon_{\beta, \delta} r, \bar{\epsilon}_{\beta, \delta} s)
	\end{equation}
	for all $r,s \in \ring$. If $2_\ring \ne 0_\ring$, then the signs $\epsilon_{\beta, \delta}, \bar{\epsilon}_{\beta, \delta}$ are uniquely determined by~\eqref{eq:fold-ex-parmap} and they do not depend on $\ring$. If, on the other hand, $2_\ring = 0_\ring$, then $\epsilon_{\beta, \delta}, \bar{\epsilon}_{\beta, \delta}$ may be chosen arbitrarily. We define
	\[ \map{\inverparsym}{H_4 \times \Set{\rho_0, \rho_1, \rho_2, \rho_3}}{\Set{\pm 1} \times \Set{\pm 1}}{(\beta, \delta)}{\inverpar{\beta}{\delta} = \inverparbr{\beta}{\delta} \defl (\epsilon_{\beta, \delta}, \bar{\epsilon}_{\beta, \delta})} \]
	where $\epsilon_{\beta, \delta}, \bar{\epsilon}_{\beta, \delta}$ denote the unique choice of signs that satisfies~\eqref{eq:fold-ex-parmap} for all commutative rings $\ring$. The map $\inverparsym$ is called	the \emph{standard parity map for $H_4$}. It will be studied in more detail in Section~\ref{sec:gpt}. For all $\beta \in H_4$ and $\delta \in \{\rho_0, \rho_1, \rho_2, \rho_3\}$, the value of $\inverparbr{\beta}{\delta}$ can be computed in a straightforward way, and this yields that $\inverparbr{-\beta}{\delta} = \inverparbr{\beta}{\delta}$ (see~\ref{rem:gap}).
	All values of $\inverparbr{\beta}{\delta}$ for positive $\beta$ are given in Figures~\ref{fig:parmapex-H4-1} and~\ref{fig:parmapex-H4-2} in the appendix. Note that, while $\inverparsym$ describes important structure constants of $G$, it can be defined without reference to $G$ by simply listing all values of $\inverparsym$, as in the aforementioned tables.
	
	We can apply similar computations to the Chevalley group of type $D_6$ from~\ref{ex:chev-D6} to obtain another parity map $\map{\inverparsym'}{H_3 \times \Set{\rho_1, \rho_2, \rho_3}}{\Set{\pm 1} \times \Set{\pm 1}}{}{}$. The set $T(E_8)$ in~\ref{rem:ex-twist-signs} has been chosen precisely in a way to ensure that $\inverparbr{\beta}{\delta} = \inverparsym'(\beta, \delta)$ for all $\delta \in \Set{\rho_1, \rho_2, \rho_3}$ and all $\beta \in H_4 \cap \gen{\rho_1, \rho_2, \rho_3}_\R$. In other words, $\inverparsym'$ is the restriction of $\inverparsym$ to $H_3$. See~\ref{rem:gap} for a verification of this fact. For the reader's convenience, the values of $ \inverparsym' $ are given in Figure~\ref{fig:parmapex-H3}. In the following, we will simply write $\inverparsym$ for~$\inverparsym'$.
\end{rem}

\premidfigure
\begin{figure}[htb]
	\centering
	\begin{tabular}{cccc}
		\toprule
		$ \alpha $ & $\inverparbr{\alpha}{\rho_1}$ & $\inverparbr{\alpha}{\rho_2}$ & $\inverparbr{\alpha}{\rho_3}$ \\
		\midrule
		$ \hcoord{0}{0}{1} $ & $(1,1)$ & $(-1,-1)$ & $(-1,-1)$ \\
		$ \hcoord{0}{\gold}{1} $ & $(1,1)$ & $(1,-1)$ & $(1,1)$ \\
		$ \hcoord{0}{\gold}{\gold} $ & $(-1,-1)$ & $(1,1)$ & $(1,-1)$ \\
		$ \hcoord{\gold}{\gold}{1} $ & $(-1,1)$ & $(1,1)$ & $(-1,-1)$ \\
		$ \hcoord{0}{1}{\gold} $ & $(-1,1)$ & $(1,-1)$ & $(-1,-1)$ \\
		$ \hcoord{\gold}{\gold}{\gold} $ & $(1,-1)$ & $(-1,-1)$ & $(-1,1)$ \\
		$ \hcoord{0}{1}{0} $ & $(-1,-1)$ & $(-1,-1)$ & $(1,-1)$ \\
		$ \hcoord{1}{1}{\gold} $ & $(1,-1)$ & $(1,1)$ & $(-1,1)$ \\
		$ \hcoord{\gold}{\gold^2}{\gold} $ & $(-1,1)$ & $(1,1)$ & $(-1,-1)$ \\
		$ \hcoord{1}{1}{0} $ & $(1,1)$ & $(-1,-1)$ & $(1,1)$ \\
		$ \hcoord{1}{\gold^2}{\gold} $ & $(-1,-1)$ & $(-1,1)$ & $(1,-1)$ \\
		$ \hcoord{\gold}{\gold^2}{\gold^2} $ & $(1,1)$ & $(-1,-1)$ & $(1,1)$ \\
		$ \hcoord{1}{0}{0} $ & $(-1,-1)$ & $(1,1)$ & $(1,1)$ \\
		$ \hcoord{1}{\gold^2}{\gold^2} $ & $(1,-1)$ & $(1,1)$ & $(-1,1)$ \\
		$ \hcoord{\gold}{2\gold}{\gold^2} $ & $(1,1)$ & $(-1,1)$ & $(1,1)$ \\
		\bottomrule
	\end{tabular}
	\caption{Values of the parity map $ \map{\inverparsym'}{H_3 \times \Set{\rho_1, \rho_2, \rho_3}}{\Set{\pm 1} \times \Set{\pm 1}}{}{} $ in~\ref{rem:fold-ex-parmap}. The values $ \inverpar{\alpha}{\delta} $ for negative roots $ \alpha $ are given by the formula $ \inverpar{-\alpha}{\delta} = \inverpar{\alpha}{\delta} $.}
	\label{fig:parmapex-H3}
\end{figure}
\postmidfigure

Recall the definition of commutation maps in~\ref{def:commmap}. In the situation of the following~\ref{lem:excommrel}, we have $ M = \ring \times \ring $.

\begin{lemma}
	\label{lem:excommrel}
	Let $ \ring $ be a commutative ring, let $ X \in \Set{D_6, E_8} $ and consider the group $ \Chev(X, \ring) $ from~\ref{ex:chev-D6} or~\ref{ex:chev-E8}. Further, let $ (\risomHex{\beta})_{\beta \in H_\ell} $ be the family of root isomorphisms from~\ref{rem:ex-twist-signs}. Put $\alpha\coloneqq\rho_2$, $\ep\coloneqq\rho_3 $ and denote the corresponding $H_2$-quintuple by $(\alpha,\beta,\gamma,\delta,\ep)$.
	Then for all non-proportional non-adjacent roots $\zeta, \xi \in \paratwopos(\rootbaseH)$ (as defined in~\ref{def:paratwopos}) and all $\rho \in \oprootint{\zeta, \xi}$, the commutation map $ \commmapex{\zeta, \xi}{\rho} $ is given by Figure~\ref{fig:excommrel}.
\end{lemma}
\begin{proof}
	This is a straightforward (matrix) computation. See~\ref{rem:gap} for details.
\end{proof}

\premidfigure
\begin{figure}[htb]
	\centering$ \begin{aligned}
		&\commmapex{\rho_0, \rho_1}{}\brackets[\big]{(a,b),(c,d)} = (ac, bd), & & \commmapex{\rho_1, \rho_0}{}\brackets[\big]{(a,b), (c,d)} = (-ac, -bd), \\[2mm]
		&\commmapex{\rho_1, \rho_2}{}\brackets[\big]{(a,b),(c,d)} = (ac, bd), & & \commmapex{\rho_2, \rho_1}{}\brackets[\big]{(a,b), (c,d)} = (-ac, -bd), \\[2mm]
		&\commmapex{\alpha,\gamma}{}\brackets[\big]{(a,b),(c,d)}=(0,ac), & &\commmapex{\gamma,\alpha}{}\brackets[\big]{(a,b),(c,d)}=(0,-ac),\\
		&\commmapex{\beta,\delta}{}\brackets[\big]{(a,b),(c,d)}=(0,-ac), & &\commmapex{\delta,\beta}{}\brackets[\big]{(a,b),(c,d)}=(0,ac),\\
		&\commmapex{\gamma,\ep}{}\brackets[\big]{(a,b),(c,d)}=(0,ac), & &\commmapex{\ep,\gamma}{}\brackets[\big]{(a,b),(c,d)}=(0,-ac),\\[2mm]
		&\commmapex{\alpha,\delta}{\beta}\brackets[\big]{(a,b),(c,d)}=(0,-bc), & &\commmapex{\delta,\alpha}{\beta}\brackets[\big]{(a,b),(c,d)}=(0,ad),\\
		&\commmapex{\alpha,\delta}{\gamma}\brackets[\big]{(a,b),(c,d)}=(0,ad), & &\commmapex{\delta,\alpha}{\gamma}\brackets[\big]{(a,b),(c,d)}=(0,-bc),\\
		&\commmapex{\beta,\ep}{\gamma}\brackets[\big]{(a,b),(c,d)}=(0,bc), & &\commmapex{\ep,\beta}{\gamma}\brackets[\big]{(a,b),(c,d)}=(0,-ad),\\
		&\commmapex{\beta,\ep}{\delta}\brackets[\big]{(a,b),(c,d)}=(0,-ad), & &\commmapex{\ep,\beta}{\delta}\brackets[\big]{(a,b),(c,d)}=(0,bc),\\[2mm]
		&\commmapex{\alpha,\ep}{\beta}\brackets[\big]{(a,b),(c,d)}=(bc,abd), & &\commmapex{\ep,\alpha}{\beta}\brackets[\big]{(a,b),(c,d)}=(-ad,-bcd),\\
		&\commmapex{\alpha,\ep}{\gamma}\brackets[\big]{(a,b),(c,d)}=(-bd,abcd), & &\commmapex{\ep,\alpha}{\gamma}\brackets[\big]{(a,b),(c,d)}=(bd,-abcd),\\
		&\commmapex{\alpha,\ep}{\delta}\brackets[\big]{(a,b),(c,d)}=(ad,-bcd), & &\commmapex{\ep,\alpha}{\delta}\brackets[\big]{(a,b),(c,d)}=(-bc,abd).
	\end{aligned} $
\caption{The commutator relations in the folding of a Chevalley group. See~\ref{lem:excommrel} for details.}
\label{fig:excommrel}
\end{figure}
\postmidfigure

\begin{rem}\label{rem:fold-comm-sub}
	Similar to the behaviour of the parity map in~\ref{rem:fold-ex-parmap}, we observe that the commutator relations in Figure~\ref{fig:excommrel} on the canonical $H_3$-subsystem of $H_4$ hold for the group from~\ref{ex:chev-D6} as well as for the group from~\ref{ex:chev-E8}. Again, this fact is due to our choice of $T(E_8)$ in~\ref{rem:ex-twist-signs}. Further, the inclusion of $ \delta_1 $ in $ T(D_6) $ and of $ \delta_6 $ in $ T(E_8) $ ensures that the formula for $ \commmapex{\rho_1,\rho_2}{} $ does not involve signs.
\end{rem}

\begin{rem}\label{rem:chev-comm-sign}
	The commutator relations of the $X$-grading of $\Chev(X,\ring)$, which we will need in~\ref{lem:D6-commrel-parabolic}, are easier to state: For all roots $\xi, \zeta \in X$ for which $\xi+\zeta \in X$, there exists a sign $c_{\xi, \zeta} \in \Set{\pm 1}$ such that
	\[ \commutator{\risomXex{\xi}(r)}{\risomXex{\zeta}(s)} = \risomXex{\xi +\zeta}(c_{\xi, \zeta} rs) \]
	for all $r,s \in \ring$. See~\cite[(9.4)]{Vavilov-Plotkin} for details, but note that \cite{Vavilov-Plotkin} uses the non-twisted family $(\risomXnaive{\xi})_{\xi \in X}$ of root isomorphisms.
\end{rem}

\begin{prop}\label{folding-is-rgg}
	Let $(G, (\rootgrX{\alpha})_{\alpha \in X})$ be any $X$-graded group. Assume that for all positive systems $\possys$ in $X$ and all $\alpha_1, \alpha_2 \in \possys$ with $\foldproj(-\alpha_1) = \foldproj(-\alpha_2)$, we have $\rootgrX{\possys} \cap \rootgrX{\alpha_1} \rootgrX{\alpha_2} = \Set{1_G}$. Then the folding $(\rootgrH{\beta})_{\beta \in H_\ell}$ defined in~\ref{def:rootgr-fold} is an $H_\ell$-grading of $G$.
\end{prop}
\begin{proof}
	Since every root group in $(\rootgrX{\alpha})_{\alpha \in X}$ is contained in a root group in $(\rootgrH{\beta})_{\beta \in H_\ell}$, the groups $(\rootgrH{\beta})_{\beta \in H_\ell}$ generate $G$. Further, $(G, (\rootgrH{\beta})_{\beta \in H_\ell})$ has Weyl elements by~\ref{lem:fold-weyl}.
	
	Denote by $\possysX$ and $\possysH$ the positive systems in $X$ and $H_\ell$ corresponding to the root bases $\rootbaseX$ and $\rootbaseH$, respectively. Since $\foldproj(\rootbaseX) = \rootbaseH$ by~\ref{pro:fold-abstract-description}~\ref{pro:fold-abstract-description:char}, we have $ \foldproj(\possysX) = \possysH $ and thus $\rootgrX{\possysX} = \rootgrH{\possysH}$. Hence our additional assumption on $\rootgrX{\possys} \cap \rootgrX{-\alpha_1} \rootgrX{-\alpha_2}$ says precisely that $(\rootgrH{\beta})_{\beta \in H_\ell}$ satisfies Axiom~\ref{def:rgg}~\ref{def:rgg:pos}.
	
	It remains to show that $G$ has $H_\ell$-commutator relations with root groups $(\rootgrH{\beta})_{\beta \in H_\ell}$. For the special case of the Chevalley groups in~\ref{ex:chev-D6} and~\ref{ex:chev-E8}, this is known by the matrix computations in~\ref{lem:excommrel}. For arbitrary $X$-graded groups, this is more technical and will be proven in~\ref{pro:fold-comm}~\ref{pro:fold-comm:H} in the appendix.
\end{proof}

\begin{rem}
	For any root system $\roots$, every Chevalley group $(G, (\rootgr{\alpha})_{\alpha \in \roots})$ satisfies that $\rootgr{\possys} \cap \rootgr{-\possys} = \Set{1_G}$. Further, every $X$-graded group with this property satisfies the additional assumption on $\rootgrX{\possys} \cap \rootgrX{-\alpha_1} \rootgrX{-\alpha_2}$ in~\ref{folding-is-rgg}. In fact, we are not aware of any examples of root graded groups which do not satisfy $\rootgr{\possys} \cap \rootgr{-\possys} = \Set{1_G}$.
\end{rem}

\section{The action of Weyl elements in \texorpdfstring{$H_3$}{H3}-graded groups}
\label{sec:h3}

\begin{notation}\label{conv:braidweyl}
	In this section, $ (G, (\rootgr{\alpha})_{\alpha \in H_3}) $ denotes an $ H_3 $-graded group, $ \Psi $ is a fixed $ H_2 $-subsystem of $ H_3 $ and $(\alpha,\beta,\gamma,\delta,\ep)$ is a fixed $ H_2 $-quintuple in $ \Psi $.
\end{notation}

The goal of this section is~\ref{prop:h3-weyl-summary}, which provides formulas for the actions of squares of Weyl elements on all root groups in $ G $. Namely, we will see that these actions \enquote{factor through} an action of the group $ \twistgroup \defl \{\pm 1\}^2 $, the so-called \emph{twisting group of $ G $} (see~\ref{def:twistgrp}). The existence of a twisting group is the main prerequisite that is necessary to apply the parametrisation theorem in Section~\ref{sec:gpt}.

\begin{rem}
	Throughout this section, we will often implicitly use~\ref{prop:parabolicsubgrading} to reduce statements about $ H_3 $-graded groups to statements about its rank-2 subgroups. Similarly, we will see in~\ref{rem:weyl-h4} that all results in this section transfer to the case of $ H_4 $-graded groups without difficulty. Note that all rank-2 subsystems of $ H_3 $ are of type $ H_2 $, $ A_2 $ or $ A_1 \times A_1 $. In particular, pairs of orthogonal roots are necessarily adjacent. Further, we will frequently apply many statements in this section to the $ H_2 $-quintuple $ (\epsilon, \delta, \gamma, \beta, \alpha) $ in place of $ (\alpha, \beta, \gamma, \delta, \epsilon) $, which is legitimate because $ (\alpha, \beta, \gamma, \delta, \epsilon) $ is chosen arbitrarily.
\end{rem}

The following result summarises what is known from the case of $ A_2 $-gradings.

%

\begin{prop}\label{prop:a2gradprop}
	Let $\zeta\in H_3$ be any root. Then the following hold: 
		\begin{enumerate}[label=(\roman*)]
			\item \label{prop:a2gradprop:abel}The root group $U_\zeta$ is abelian.
			\item \label{prop:a2gradprop:weyl-left-right}For each $b_\zeta\in U_\zeta^\sharp$, there exists a unique $a_{-\zeta}\in U_{-\zeta}$ such that $a_{-\zeta}b_\zeta a_{-\zeta}$ is an $\zeta$-Weyl element.
			\item \label{prop:a2gradprop:factor-unique}If $w_\zeta$ is an $\zeta$-Weyl element, then there exist unique $a_{-\zeta}\in U_{-\zeta}$ and $b_\zeta\in U_\zeta$ such that $w_\zeta=a_{-\zeta}b_\zeta a_{-\zeta}$.
			\item \label{prop:a2gradprop:conj-formula}Let $w_\zeta=a_{-\zeta}b_\zeta a_{-\zeta}$ be a $\zeta$-Weyl element with $ a_{-\zeta} \in \rootgr{-\zeta} $ and $ b_\zeta \in \rootgr{\zeta} $. Let $\xi\in H_3$ such that $(\zeta, \xi)$ is an $ A_2 $-pair and let $x_\xi\in U_\xi$. Then $x_\xi^{w_\zeta}=[x_\xi,b_\zeta]$ and $x_\xi = \commutator{\commutator{x_\xi}{b_\zeta}}{a_{-\zeta}}^{-1}$. In particular, $\rootgr{\xi+\zeta}=\rootgr{\xi}^{w_\zeta} = \commutator{\rootgr{\xi}}{\rootgr{\zeta}}$.
			\item \label{prop:a2gradprop:balanced}Let $w_\zeta=a_{-\zeta}b_\zeta a_{-\zeta}$ be a $\zeta$-Weyl element with $ a_{-\zeta} \in \rootgr{-\zeta} $ and $ b_\zeta \in \rootgr{\zeta} $. Then $ (a_{-\zeta}, b_{\zeta}, a_{-\zeta}) $ is balanced in the sense of \cite[2.2.14]{torben}. That is, we have $ a_{-\zeta}^{w_\zeta} = b_\zeta $, $ b_\zeta^{w_\zeta} = a_{-\zeta} $ and $ w_\zeta = b_\zeta a_{-\zeta} b_\zeta $.
		\end{enumerate}		 
\end{prop}
\begin{proof}
	Since $ \alpha $ lies in an $ A_2 $-subsystem by~\ref{prop:222}, this follows from the corresponding assertions about $ A_2 $-graded groups in \cite[5.4.3,~5.4.9,~5.4.10]{torben}.
\end{proof}

\begin{prop}
	\label{thm:weyla2}
	Let $\zeta,\xi\in H_3$, let $w_\zeta$ be an $\zeta$-Weyl element and put $ w \coloneqq w_\zeta^2 $. Then the following hold:
	\begin{enumerate}[(i)]
		\item If $ \subsys{\zeta, \xi} $ is of type $ A_2 $, then $w$ acts on $U_\xi$ by inversion. That is, $x_\xi^{w}=x_\xi\inv$ for all $x_\xi\in U_\xi$.
		\item If $ \subsys{\zeta, \xi} $ is of type $ A_1 $ or $ A_1 \times A_1 $, then $w$ acts trivially on $U_\xi$. That is, $x_\xi^{w}=x_\xi$ for all $x_\xi\in U_\xi$.
	\end{enumerate}
\end{prop}
\begin{proof}
	If $ \subsys{\zeta, \xi} $ is of type $ A_1 \times A_1 $, then $ \rootgr{\xi} $ is centralised by $ \rootgr{\zeta} $ and $ \rootgr{-\zeta} $ and hence by $ w $. The remaining cases are covered by \cite[5.4.13,~5.4.15]{torben}.
\end{proof}

It remains to study the more complicated case of $ \zeta $-Weyl elements acting on $ \rootgr{\xi} $ where $ \subsys{\zeta, \xi} $ is of type $ H_2 $. We begin with collecting some \enquote{rank-2 results}, that is, results which hold in arbitrary $ H_2 $-graded groups (and not only those which are embeddable into an $ H_3 $-graded group).

\begin{lemma}
\label{lem:i25additivity}
	Let $x_\alpha,x'_\alpha\in U_\alpha$ and $y_\ep,y'_\ep\in U_\ep$. Then the following hold:
	\begin{enumerate}[(i)]
		\item \label{lem:i25additivity:right}$[x_\alpha,y_\ep y'_\ep]_\beta = [x_\alpha,y'_\ep]_\beta [x_\alpha,y_\ep]_\beta = [x_\alpha,y_\ep]_\beta [x_\alpha,y'_\ep]_\beta$.
		\item \label{lem:i25additivity:right:switch}$[x_\alpha x'_\alpha,y_\ep]_\delta = [x_\alpha,y_\ep]_\delta [x'_\alpha,y_\ep]_\delta$.
	\end{enumerate}
	In particular, $[x_\alpha,y_\ep\inv]_\beta=[x_\alpha,y_\ep]_\beta\inv$ and $[x_\alpha\inv,y_\ep]_\delta=[x_\alpha,y_\ep]_\delta\inv$.
\end{lemma}
\begin{proof}
	Using~\ref{rem:commrel}~\ref{rem:commrel:add-right}, we compute
	\begin{align*}
		[x_\alpha,y_\ep y'_\ep]&=[x_\alpha, y'_\ep][x_\alpha,y_\ep]^{y'_\ep}\\
		&=[x_\alpha, y'_\ep]_\beta [x_\alpha, y'_\ep]_\gamma [x_\alpha, y'_\ep]_\delta \big([x_\alpha,y_\ep]_\beta [x_\alpha,y_\ep]_\gamma [x_\alpha,y_\ep]_\delta \big)^{y'_\ep}.
	\end{align*}
	By~\ref{rem:commrel}~\ref{rem:commrel:conj}, $([x_\alpha,y_\ep]_\beta [x_\alpha,y_\ep]_\gamma [x_\alpha,y_\ep]_\delta)^{y'_\ep}$ is equal to
	\[\left[y'_\ep,[x_\alpha,y_\ep]_\beta\inv\right]_\delta \left[y'_\ep,[x_\alpha,y_\ep]_\beta\inv\right]_\gamma [x_\alpha,y_\ep]_\beta \left[y'_\ep,[x_\alpha,y_\ep]_\gamma\inv\right][x_\alpha,y_\ep]_\gamma [x_\alpha,y_\ep]_\delta.\]
	Note that the only factor in this product that lies in $U_\beta$ is $[x_\alpha,y_\ep]_\beta$. Using relation~\ref{rem:commrel}~\ref{rem:commrel:comm}, we can reorder the factors to obtain that
	\[ [x_\alpha,y_\ep y'_\ep] = [x_\alpha, y'_\ep]_\beta [x_\alpha,y_\ep]_\beta z_\gamma z_\delta \]
	for some $ z_\gamma \in \rootgr{\gamma} $ and $ z_\delta \in \rootgr{\delta} $. By \ref{prop:prodmapbij}, the first equation in~\ref{lem:i25additivity:right} follows. Since $ \rootgr{\beta} $ is abelian, the second equation in~\ref{lem:i25additivity:right} holds as well. Finally, assertion~\ref{lem:i25additivity:right:switch} follows from~\ref{lem:i25additivity:right} (applied to the $ H_2 $-quintuple $ (\epsilon, \delta, \gamma, \beta, \alpha) $) and~\ref{lem:commpart}:
	\begin{align*}
		\commpart{x_\alpha x_\alpha'}{y_\epsilon}{\delta} &= \commpart{y_\epsilon}{x_\alpha x_\alpha'}{\delta}^{-1} = \brackets[\big]{\commpart{y_\epsilon}{x_\alpha'}{\delta} \commpart{y_\epsilon}{x_\alpha}{\delta}}^{-1} = \commpart{x_\alpha}{y_\epsilon}{\delta} \commpart{x_\alpha'}{y_\epsilon}{\delta},
	\end{align*}
	as desired.
\end{proof}

\begin{lemma}\label{lem:int2add}
	Let $ x_\alpha, x_\alpha' \in \rootgr{\alpha} $ and $ y_\delta, y_\delta' \in \rootgr{\delta} $. Then the following hold:
	\begin{enumerate}[(i)]
		\item $ \commpart{x_\alpha}{y_\delta y_\delta'}{\beta} = \commpart{x_\alpha}{y_\delta}{\beta} \commpart{x_\alpha}{y_\delta'}{\beta} $.
		
		\item $ \commpart{x_\alpha x_\alpha'}{y_\delta}{\beta} = \commpart{x_\alpha}{y_\delta}{\beta} \commpart{x_\alpha'}{y_\delta}{\beta} \commutator[\big]{\commpart{x_\alpha}{y_\delta}{\gamma}}{x_\alpha'} $.
		
		\item $ \commpart{x_\alpha}{y_\delta y_\delta}{\beta} = \commpart{x_\alpha}{y_\delta}{\beta} \commpart{x_\alpha}{y_\delta'}{\beta} \commutator[\big]{\commpart{x_\alpha}{y_\delta}{\beta}}{y_\delta'} $.
		
		\item $ \commpart{x_\alpha x_\alpha'}{y_\delta}{\gamma} = \commpart{x_\alpha}{y_\delta}{\gamma} \commpart{x_\alpha'}{y_\delta}{\gamma} $.
	\end{enumerate}
\end{lemma}
\begin{proof}
	This is a similar computation as in~\ref{lem:i25additivity}. See \cite[7.4.4]{torben} for details.
\end{proof}

\begin{lemma}
\label{lem:weylactioni25pair}
	 Let $x_\ep\in U_\ep$ and let $w_\alpha=a_{-\alpha}b_\alpha a_{-\alpha}$ be an $\alpha$-Weyl element for $b_\alpha\in U_\alpha^\sharp$ and $a_{-\alpha}\in \invset{-\alpha}$. Then the following hold:
	 	\begin{enumerate}[label=(\roman*)]
	 		\item\label{lem:weylactioni25pair:basic} $x_\ep^{w_\alpha}=[b_\alpha,x_\ep\inv]_\beta$.
	 		\item\label{lem:weylactioni25pair:square} $x_\ep^{w_\alpha^2}=\left[a_{-\alpha},[b_\alpha,x_\ep\inv]_\beta\inv\right]_\ep=\left[a_{-\alpha},[b_\alpha,x_\ep]_\beta\right]_\ep$.
	 		\item\label{lem:weylactioni25pair:basic-inv} $x_\ep^{w_\alpha\inv}=[b_\alpha\inv,x_\ep\inv]_\beta$.
	 		\item\label{lem:weylactioni25pair:cancel} $x_\ep=\left[a_{-\alpha}\inv,[b_\alpha,x_\ep]_\beta\right]_\ep$.
		\end{enumerate}	 	 
\end{lemma}
\begin{proof}
	Let $b_\alpha\in U_\alpha^\sharp$ and $a_{-\alpha}\in U_{-\alpha}$ such that $w_\alpha=a_{-\alpha}b_\alpha a_{-\alpha}$ is an $\alpha$-Weyl element and let $x_\ep\in U_\ep$. Using~\ref{rem:commrel}~\ref{rem:commrel:conj}, we compute
	\begin{align*}
		x_\ep^{w_\alpha}&=x_\ep^{a_{-\alpha}b_\alpha a_{-\alpha}}=x_\ep^{b_\alpha a_{-\alpha}}=([b_\alpha,x_\ep\inv]x_\ep)^{a_{-\alpha}} \\
		&=[b_\alpha,x_\ep\inv]_\beta^{a_{-\alpha}} [b_\alpha,x_\ep\inv]_\gamma^{a_{-\alpha}} [b_\alpha,x_\ep\inv]_\delta^{a_{-\alpha}} x_\ep \\
		&= \left[a_{-\alpha},[b_\alpha,x_\ep\inv]_\beta\inv\right][b_\alpha,x_\ep\inv]_\beta\left[a_{-\alpha},[b_\alpha,x_\ep\inv]_\gamma\inv\right][b_\alpha,x_\ep\inv]_\gamma \\
		& \qquad \mathord{} \cdot \big[a_{-\alpha},[b_\alpha,x_\ep\inv]_\delta\inv\big][b_\alpha,x_\ep\inv]_\delta x_\ep.
	\end{align*}
	Just as in the proof of \ref{lem:i25additivity}, reordering the elements in the above product does not produce any new factors in $U_\beta$ other than $[b_\alpha,x_\ep\inv]_\beta$. Since $x_\ep^{w_\alpha}\in U_{\ep^{\reflbr{\alpha}}}=U_\beta$, we conclude that~\ref{lem:weylactioni25pair:basic} holds. By applying~\ref{lem:weylactioni25pair:basic} to the $ \alpha $-Weyl element $w_\alpha\inv=a_{-\alpha}\inv b_\alpha\inv a_{-\alpha}\inv$, we obtain~\ref{lem:weylactioni25pair:basic-inv}.
	
	For~\ref{lem:weylactioni25pair:square}, recall from~\ref{prop:weylelements} that $w_\alpha$ is also a ($-\alpha$)-Weyl element with factorisation $w_\alpha=b_\alpha a_{-\alpha}a_{-\alpha}^{w_\alpha}$. Hence by two applications of~\ref{lem:weylactioni25pair:basic},
	\begin{align*}
		x_\ep^{w_\alpha^2}&= [b_\alpha,x_\ep\inv]_\beta^{w_\alpha} =\left[a_{-\alpha},[b_\alpha,x_\ep\inv]_\beta\inv\right]_\ep
	\end{align*}
	The second equality in~\ref{lem:weylactioni25pair:square} now follows by \ref{lem:i25additivity}.
	
	Finally, we have
	\[x_\ep=x_\ep^{w_\alpha w_\alpha\inv}=[b_\alpha,x_\ep\inv]_\beta^{w_\alpha\inv}.\]
	Using the factorisation $w_\alpha\inv=(a_{-\alpha}\inv)^{w_\alpha} a_{-\alpha}\inv b_\alpha\inv$ of $ w_\alpha^{-1} $ as a $ (-\alpha) $-Weyl element (see~\ref{prop:weylelements}), an application of~\ref{lem:weylactioni25pair:basic} to the $ H_2 $-quintuple $ (-\alpha, \ep, \delta, \gamma, \beta) $ yields $x_\ep=\left[a_{-\alpha}\inv,[b_\alpha,x_\ep]_\beta\right]_\ep$. This finishes the proof.
\end{proof}

\begin{prop}
\label{prop:betaepbijection}
Let $b_\alpha\in U_\alpha^\sharp$. Then the map $U_\ep\rightarrow U_\beta,~x_\ep\mapsto[b_\alpha,x_\ep]_\beta$ is an isomorphism of groups. In particular, for each $x_\beta\in U_\beta$ there exist $x_\alpha\in U_\alpha$ and $x_\ep\in U_\ep$ such that $x_\beta=[x_\alpha,x_\ep]_\beta$.
\end{prop}
\begin{proof}
Let $ w_\alpha $ be an $ \alpha $-Weyl element in $ \rootgr{-\alpha} b_\alpha \rootgr{-\alpha} $. By \ref{lem:i25additivity} and \ref{lem:weylactioni25pair}, the map in the assertion is $x_\ep\mapsto (x_\ep^{w_\alpha})\inv$, which is an isomorphism by the definition of Weyl elements and because $ \rootgr{\beta} $ is abelian.
\end{proof}

From now on, we will use in a crucial way that the grading of $ G $ has rank (at least)~3. At first, we show that the action of the square of a $ \zeta $-Weyl element (for a fixed root $ \zeta $) does not depend on the choice of the Weyl element.

\begin{prop}
\label{thm:weyli25}
Let $w_\alpha$, $ v_\alpha $ be two $\alpha$-Weyl elements and let $\zeta\in \Psi \setminus \Set{\alpha, -\alpha}$. Then $ w_\alpha^2 $ and $ v_\alpha^2 $ act identically on $ \rootgr{\zeta} $ and $w_\alpha^4$ acts trivially on $U_\zeta$.
\end{prop}
\begin{proof}
	Since $ w_\alpha $ is also a $ (-\alpha) $-Weyl element by~\ref{prop:weylelements}~\ref{prop:weylelements:minus}, we may assume that $\zeta\in\{\beta,\gamma,\delta,\ep\}$. In each case, we will derive a formula for the action of $ w_\alpha^2 $ on $ \rootgr{\zeta} $ which does not depend on the choice of $ w_\alpha $, which proves the first assertion. It then follows that $ w_\alpha^2 $ and $ (w_\alpha^{-1})^2 $ act identically on $ \rootgr{\zeta} $, which says precisely that $ w_\alpha^4 $ acts trivially on $ \rootgr{\zeta} $.
	
	Recall from \ref{prop:222} that $\zeta$ lies in a unique $H_2$-subsystem $\tilde{\Psi}$ different from $\Psi$. We choose an $H_2$-quintuple $(\tilde{\alpha},\tilde{\beta},\tilde{\gamma},\tilde{\delta},\tilde{\ep})$ in $\tilde{\Psi}$ such that $\zeta=\tilde{\beta}$. Since any rank-2 subsystem is uniquely determined by two non-proportional roots in this subsystem, we have $ \Psi \cap \tilde{\Psi} = \Set{\zeta, -\zeta} = \Set{\tilde{\beta}, -\tilde{\beta}} $ and
	\begin{equation}\label{eq:weyli25:distinct}
		\Psi, \subsys{\alpha, \tilde{\alpha}}, \subsys{\alpha, \tilde{\ep}}, \subsys{\alpha, \tilde{\gamma}}, \subsys{\alpha, -\tilde{\delta}} \text{ are pairwise distinct.}
	\end{equation}
	
	Now let $ x_\zeta \in \rootgr{\zeta} $ and put $ w \defl w_\alpha^2 $. By \ref{prop:betaepbijection} we can write $x_\zeta=[x_{\tilde{\alpha}},x_{\tilde{\ep}}]_\zeta$ for some $x_{\tilde{\alpha}}\in U_{\tilde{\alpha}}$ and $x_{\tilde{\ep}}\in U_{\tilde{\ep}}$. Hence by~\ref{lem:commpart},
	\begin{equation*}
		x_\zeta^w = \commpart{x_{\tilde{\alpha}}}{x_{\tilde{\ep}}}{\zeta}^w = \commpart{x_{\tilde{\alpha}}^w}{x_{\tilde{\ep}}^w}{\zeta}.
	\end{equation*}
	If $ \subsys{\alpha, \tilde{\alpha}} $ and $ \subsys{\alpha, \tilde{\ep}} $ are of type $ A_1 \times A_1 $ or $ A_2 $, then~\ref{thm:weyla2} provides formulas for the action of $ w $ on $ \rootgr{\tilde{\alpha}} $ and $ \rootgr{\tilde{\ep}} $ which do not depend on the choice of $ w_\alpha $. Thus the first assertion holds in this case. Further, since $ \alpha $ lies in exactly two $ H_2 $-subsystems by~\ref{prop:222}, it follows from~\eqref{eq:weyli25:distinct} that at most one of the subsystems $ \subsys{\alpha, \tilde{\alpha}} $ and $ \subsys{\alpha, \tilde{\ep}} $ can be of type $ H_2 $.
	
	We remain with the case that one of the subsystems $ \subsys{\alpha, \tilde{\alpha}} $ and $ \subsys{\alpha, \tilde{\ep}} $ is of type $ H_2 $. Again by applying~\ref{prop:betaepbijection}, but to the $ H_2 $-pair $ (\tilde{\gamma}, -\tilde{\delta}) $, we can write $x_\zeta=[x_{\tilde{\gamma}},x_{-\tilde{\delta}}]_\zeta$ for some $x_{\tilde{\gamma}}\in U_{\tilde{\gamma}}$ and $x_{-\tilde{\delta}}\in U_{-\tilde{\delta}}$. As above, we then have $ x_\zeta^w = \commpart{x_{\tilde{\gamma}}^w}{x_{-\tilde{\delta}}^w}{\zeta} $.
	By our case assumption,~\ref{prop:222} and~\eqref{eq:weyli25:distinct}, both $ \subsys{\alpha, \tilde{\gamma}} $ and $ \subsys{\alpha, -\tilde{\delta}} $ are of type $ A_1 \times A_1 $ or $ A_2 $. By~\ref{thm:weyla2}, we conclude that the first assertion holds in this case as well, which finishes the proof.
\end{proof}

We may already suspect from the proof of~\ref{thm:weyli25} that the action of $ w_\alpha^2 $ on $ \rootgr{\zeta} $ depends on the choice of $ \zeta \in \Psi \setminus \{\alpha, -\alpha \} $. However, we will see in~\ref{thm:possametrivial} that only the two cases in the following definition are distinct, and in~\ref{thm:posdiffinv} that the actions in these two cases differ by inversion.

\begin{definition}
\label{def:i25pos}
Let $\zeta,\xi \in H_3$. We say that $\xi$ is \emph{in involution position with respect to $\zeta$} if $(\zeta,\xi)$ or $(\zeta,-\xi)$ is an $H_2$-pair. We say that $\xi$ is \emph{in inverted involution position with respect to~$\zeta$} if $ \subsys{\xi, \zeta} $ is of type $ H_2 $ but $ \xi $ is not in involution position with respect to $ \zeta $. For any $ \rho \in H_3 $, we say that $ \xi $, $ \zeta $ are \emph{in the same position with respect to $ \rho $} if both are in involution position with respect to $ \rho $ or both are in inverted involution position with respect to $ \rho $.
\end{definition}

\begin{rem}\label{rem:position}
	In the notation of this section, $\beta, -\beta, \ep $ and $-\ep$ are in involution position with respect to $\alpha$ while $\gamma, -\gamma,\delta$ and $ -\delta $ are in inverted involution position with respect to $ \alpha $. Further, for any $ w \in \Weyl(H_3) $, the root $ \xi $ is in (inverted) involution position with respect to $ \zeta $ if and only if $ \xi^w $ is in (inverted) involution position with respect to $ \zeta^w $. Moreover, $ \xi $ is in (inverted) involution position with respect to $ \zeta $ if and only if $ \zeta $ is in (inverted) involution position with respect to $ \xi $.
\end{rem}

\begin{notation}
\label{conv:weylg2i25}
For the rest of this section, we denote by $\tilde{\Psi}$ the unique $H_2$-subsystem different from $ \Psi $ which contains $\alpha$. Further, we choose an $H_2$-quintuple $(\alpha,\tilde{\beta},\tilde{\gamma},\tilde{\delta},\tilde{\ep})$ in $ \tilde{\Psi} $ starting with $ \alpha $. Note that $ \Psi \cap \tilde{\Psi} = \Set{\alpha, -\alpha} $. Finally, we fix a positive system $ \Pi $ in $ H_3 $ containing $ \{\alpha, \beta, \gamma, \delta, \ep\} $, whose elements will henceforth be called \enquote{positive roots}.
\end{notation}

Before we can prove~\ref{thm:possametrivial}, we need the following lemma about~$ H_3 $.

\begin{lemma}\label{rem:orthogonalroots}
	There exist unique positive roots $\eta,\vartheta\in H_3 \cap \alpha^\perp$ such that
	\begin{align*}
		\beta^{\reflbr{\eta}} &= \tilde{\beta}, & \gamma^{\reflbr{\eta}} &= \tilde{\gamma}, & \delta^{\reflbr{\eta}} &= \tilde{\delta}, & \ep^{\reflbr{\eta}} &= \tilde{\ep}, \\
		\beta^{\reflbr{\vartheta}} &= -\tilde{\ep}, & \gamma^{\reflbr{\vartheta}} &= -\tilde{\delta}, & \delta^{\reflbr{\vartheta}} & =-\tilde{\gamma}, & \ep^{\reflbr{\vartheta}} &= -\tilde{\beta}.
	\end{align*}
	Further, $ \eta $ and $ \vartheta $ are the only positive roots in $ H_3 \cap \alpha^\perp $.
\end{lemma}
\begin{proof}
	By \ref{prop:222}~\ref{prop:222:A1xA1}, there are exactly two positive roots in $ H_3 \cap \alpha^\perp $, so the uniqueness assertion is trivial. Let $ \rho $ denote any root in $ H_3 \cap \alpha^\perp $. Then $ \rho \notin \Psi $ and hence $ \rho \notin \langle \Psi \rangle_\R $. By the defining formula for reflections, it follows that for all $ \xi \in \Psi $, we have either $ \xi^{\reflbr{\rho}} = \xi $ or $ \xi^{\reflbr{\rho}} \notin \Psi $. At the same time, by~\ref{prop:222}~\ref{prop:222:A1xA1}, $ \refl{\rho} $ cannot fix all roots in $ \Psi $. Since $ \alpha^{\reflbr{\rho}} = \alpha $, we conclude that $ \Psi^{\reflbr{\rho}} $ is the unique $ H_2 $-subsystem of $ H_3 $ containing $ \alpha $ which is different from $ \Psi $. That is, $ \Psi^{\reflbr{\rho}} = \tilde{\Psi} $. Since $ \refl{\rho} $ is orthogonal and $ \alpha^{\reflbr{\rho}} = \alpha $, we must have $ \beta^{\reflbr{\rho}} \in \{ \tilde{\beta}, -\tilde{\ep} \} $.
	
	Now denote by $ \eta $ and $ \vartheta $ the two positive roots in $ H_3 \cap \alpha^\perp $. For a contradiction, suppose that $ \beta^{\reflbr{\eta}} = \beta^{\reflbr{\vartheta}} $. Then $ \refl{\eta}^{-1} \refl{\vartheta} $ is an orthogonal map of determinant $ 1 $ which fixes the hyperplane $ \gen{\Psi}_\R $ in $ \gen{H_3}_\R $ pointwise. It follows that $ \refl{\eta} = \refl{\vartheta} $, so $ \eta = \vartheta $ because both roots are positive. This is a contradiction, so $ \beta^{\reflbr{\eta}} \ne \beta^{\reflbr{\vartheta}} $. By interchanging the roles of $ \eta $ and $ \vartheta $ if necessary, we conclude that $ \beta^{\reflbr{\eta}}= \tilde{\beta} $ and $ \beta^{\reflbr{\vartheta}} = -\tilde{\ep} $. The remaining values of $ \refl{\eta} $ and $ \refl{\vartheta} $ on $ \Psi $ are uniquely determined by their values on $ \alpha $ and $ \beta $, so the assertion follows.
\end{proof}

\begin{prop}
\label{thm:possametrivial}
	Let $\xi,\zeta\in\Psi\cup\tilde{\Psi} \setminus \{\alpha, -\alpha\}$ be non-proportional and let $w_\xi\in M_\xi,~w_\zeta\in M_\zeta$ be Weyl elements. Assume that $\xi$ and $\zeta$ are in the same position with respect to $\alpha$. Then the action of $w_\xi^2w_\zeta^2$ on $U_\alpha$ is trivial. Equivalently, the actions of $w_\xi^2$ and $w_\zeta^2$ on $U_\alpha$ are identical.
\end{prop}
\begin{proof}
Let $\eta,\vartheta\in H_3 \in \alpha^\perp$ be as in \ref{rem:orthogonalroots} and choose Weyl elements $w_\eta\in M_\eta$ and $w_\vartheta\in M_\vartheta$.
At first, we consider the case that $\xi$ and $\zeta$ lie in a common $H_2$-subsystem. Without loss of generality, we assume that both lie in $ \Psi $. Further, by \ref{prop:weylelements}~\ref{prop:weylelements:minus} and \ref{thm:weyli25}, we may assume that $\zeta \neq \xi$ and that both $ \xi, \zeta $ are positive roots. Thus either $ \{\xi, \zeta\} = \{\beta, \ep\} $ or $ \{\xi, \zeta\} = \{\gamma, \delta\} $.
In any case, $\xi^{\reflbr{\eta}\reflbr{\vartheta}}=-\zeta$ by \ref{rem:orthogonalroots}, so it follows from \ref{prop:weylelements}~\ref{prop:weylelements:conj} that
	\[w_\xi^{w_\eta w_\vartheta}\in M_{\xi^{\reflbr{\eta}\reflbr{\vartheta}}}=M_{-\zeta}=M_\zeta.\]
Thus by~\ref{thm:weyli25}, $ w_\zeta^2 $ and $ (w_\xi^{w_\vartheta w_\eta})^2 $ act identically on $ \rootgr{\alpha} $. Here $ w_\vartheta $ and $ w_\eta $ act trivial on $ \rootgr{\alpha} $ because $ \vartheta, -\vartheta, \eta $ and $ -\eta $ are orthogonal and hence adjacent to $ \alpha $. We conclude that $ w_\zeta^2 $ and $ w_\xi^2 $ act identically on $ \rootgr{\alpha} $, as desired.

Now we consider the case that $\xi$ and $\zeta$ are not contained in a common $H_2$-sub\-sys\-tem. Then $ \zeta $ and $ \xi^{\refl{\eta}} $ lie in a common $ H_2 $-subsystem and, by~\ref{rem:position}, are in the same position with respect to $ \alpha $. Therefore, by the previous case, the actions of $w_\zeta^2$ and $(w_\xi^{w_\eta})^2$ on $U_\alpha$ are identical. Since $ w_\eta $ acts trivially on $ \rootgr{\alpha} $, we infer that $w_\zeta^2$ and $w_\xi^2$ act identically on $U_\alpha$, as desired.
\end{proof}

Again, we need an elementary lemma about $ H_3 $ before we can prove~\ref{thm:posdiffinv}.

\begin{lemma}\label{lem:H3-ex-A2sub}
	There exists $ \rho \in \{\tilde{\gamma}, \tilde{\delta}\} $ such that both $\Psi_{\rho,\beta}$ and $\Psi_{\rho,\gamma}$ are of type~$A_2$.
\end{lemma}
\begin{proof}
	Note that if $ \rho $ has the desired properties, then $ -\rho $ has the desired properties for $(\alpha,\tilde{\beta},\tilde{\gamma},\tilde{\delta},\tilde{\ep})$ in~\ref{conv:weylg2i25} replaced by $ (\alpha, -\tilde{\ep}, -\tilde{\delta}, -\tilde{\gamma}, -\tilde{\beta}) $. Thus for each fixed root $ \alpha $, there are only four cases to consider: Namely, there are two ways to choose $ \Psi $ and two ways to choose an $ H_2 $-quintuple starting from $ \alpha $ in $ \Psi $. It is a straightforward computation (see~\ref{rem:gap}) to prove the assertion in each of these cases for one fixed choice of $ \alpha $. Then it follows from~\ref{lem:WH3transitive} that it holds for all possible choices of~$ \alpha $.
\end{proof}

\begin{prop}
\label{thm:posdiffinv}
	Let $\xi,\zeta\in\Psi\cup\tilde{\Psi} \setminus \{\alpha, -\alpha\} $ and let $w_\xi\in M_\xi,~w_\zeta\in M_\zeta$ be Weyl elements. Assume that $\xi$ and $\zeta$ are not in the same position with respect to $\alpha$. Then $w_\xi^2w_\zeta^2$ acts on $U_\alpha$ by inversion. Equivalently, the actions of $w_\xi^2$ and $w_\zeta^2$ on $U_\alpha$ differ by inversion.
\end{prop}
\begin{proof}
If $w_\xi^2w_\zeta^2$ acts on $U_\alpha$ by inversion, then so does $ w_\zeta^{-2} w_\xi^{-2} $. Hence we may assume that $\xi$ is in involution position with respect to $\alpha$ and that $\zeta$ is in inverted involution position. By \ref{thm:possametrivial}, we can thus assume that $ \xi = \beta $ and $ \zeta = \gamma $.

Let $x_\alpha\in U_\alpha$ be arbitrary. Using \ref{prop:betaepbijection}, we can write $x_\alpha=[x_\beta,x_{-\gamma}]_\alpha$ for some $x_\beta\in U_\beta,~x_{-\gamma}\in U_{-\gamma}$. Choose $a_{-\beta}\in U_{-\beta}$ and $b_\beta\in U_\beta$ such that $w_\beta=a_{-\beta}b_\beta a_{-\beta}$. Applying~\ref{lem:weylactioni25pair}~\ref{lem:weylactioni25pair:square} to the $ (-\beta) $-Weyl element $ w_\beta = b_\beta a_{-\beta} b_\beta^{w_\beta} $ from~\ref{prop:weylelements}~\ref{prop:weylelements:minus}, we see that
\[ x_\alpha^{w_\beta^2}=\big[b_\beta,[a_{-\beta},x_\alpha]_{-\gamma}\big]_\alpha. \]
Now choose $ \rho \in \{\tilde{\gamma}, \tilde{\delta}\} $ as in~\ref{lem:H3-ex-A2sub} and choose some $ \rho $-Weyl element $ w_\rho $. Then, by~\ref{thm:weyla2}, $ w_\rho^2 $ acts on $U_\beta$ and $U_\gamma$ by inversion. Since $ w_\gamma^2 $ and $ w_\rho^2 $ act identically on $ \rootgr{\alpha} $ by~\ref{thm:possametrivial}, we conclude that
\begin{align*}
	x_\alpha^{w_\beta^2w_\gamma^2} &= \big[b_\beta,[a_{-\beta},x_\alpha]_{-\gamma}\big]_\alpha^{w_\rho^2} = \big[b_\beta^{w_\rho^2},[a_{-\beta},x_\alpha]_{-\gamma}^{w_\rho^2}\big]_\alpha = \left[b_\beta\inv,[a_{-\beta},x_\alpha]_{-\gamma}\inv\right]_\alpha.
\end{align*}
Combining \ref{lem:i25additivity} and~\ref{lem:weylactioni25pair}~\ref{lem:weylactioni25pair:cancel}, we infer that
\[ x_\alpha^{w_\beta^2w_\gamma^2} = \left[b_\beta\inv,[a_{-\beta},x_\alpha]_{-\gamma}\right]_\alpha^{-1} = x_\alpha^{-1}. \]
This shows that $ w_\beta^2 w_\gamma^2 $ acts on $ \rootgr{\alpha} $ by inversion, as desired.
\end{proof}

Before we can summarise the results of this section in~\ref{prop:h3-weyl-summary}, we need to introduce the following notation.

\begin{definition}
\label{def:invorg}
	Let $ \rho \in H_3 $. The \emph{twisting involution on $ \rootgr{\rho} $} is the map $ \map{}{\rootgr{\rho}}{\rootgr{\rho}}{x_\rho}{\invoroot{x_\rho}{\rho}} $ defined by $ \invoroot{x_\rho}{\rho} \defl x_\rho^w $ where $ \zeta $ is any root in involution position with respect to $ \rho $, $ w_\zeta $ is any $ \zeta $-Weyl element and $ w \defl w_\zeta^2 $. By~\ref{thm:possametrivial}, the twisting involution does not depend on the choices of $ \zeta $ and $ w_\zeta $. We will usually write $ \invo{x_\rho} $ in place of $ \invoroot{x_\rho}{\rho} $.
\end{definition}

\begin{lemma}\label{lem:invo-twist-weyl}
	Let $ \rho, \xi \in H_3 $, let $ w_\xi $ be a $ \xi $-Weyl element and put $ \bar{\rho} \defl \rho^{\reflbr{\xi}} $. Then $ (\invoroot{x_\rho}{\rho})^{w_\xi} = \invoroot{(x_\rho^{w_\xi})}{\bar{\rho}} $ and $ (\invoroot{x_\rho}{\rho})^{-1} = \invoroot{(x_\rho^{-1})}{\rho} $ for all $ x_\rho \in \rootgr{\rho} $.
\end{lemma}
\begin{proof}
	Choose a root $ \zeta $ in involution position with respect to $ \rho $ and put $ w \defl w_\zeta^2 $ for some $ \zeta $-Weyl element $ w_\zeta $. Then $ v \defl w^{w_\xi} = (w_\zeta^{w_\xi})^2 $ is the square of a $ \zeta^{\reflbr{\xi}} $-Weyl element. Since $ \zeta^{\reflbr{\xi}} $ is in involution position with respect to $ \rho^{\reflbr{\xi}} = \bar{\rho} $, we infer that
	\[ (\invoroot{x_\rho}{\rho})^{w_\xi} = (x_\rho^w)^{w_\xi} = (x_\rho^{w_\xi})^{v} = \invoroot{(x_\rho^{w_\xi})}{\bar{\rho}} \]
	This proves the first assertion. The second assertion is clear.
\end{proof}

\begin{prop}\label{prop:h3-weyl-summary}
	Let $ \xi, \zeta \in H_3 $, let $ x_\xi \in \rootgr{\xi} $, let $ w_\zeta $ be a $ \zeta $-Weyl element and put $ w \defl w_\zeta^2 $. Then for all $ x_\xi \in \rootgr{\xi} $, we have
	\[ x_\xi^w = \begin{cases}
		x_\xi & \text{if } \subsys{\xi, \zeta} = A_1 \text{ or } \subsys{\xi, \zeta} = A_1 \times A_1, \\
		x_\xi^{-1} & \text{if } \subsys{\xi, \zeta} = A_2, \\
		\invo{x_\xi} & \text{if } \zeta \text{ is in involution position with respect to } \xi, \\
		(\invo{x_\xi})^{-1} & \text{if } \zeta \text{ is in inverted involution position with respect to } \xi.
	\end{cases} \]
\end{prop}
\begin{proof}
	This follows from~\ref{thm:weyla2},~\ref{thm:posdiffinv} and~\ref{def:invorg}.
\end{proof}

\begin{rem}\label{rem:weyl-h4}
	We show how the previous assertions can be extended to $H_4$-graded groups $(G', (\rootgr{\alpha}')_{\alpha \in H_4})$. First, we can define the notions of \emph{involution position} and \emph{inverted involution position} for roots in $H_4$ in the same way as in~\ref{def:i25pos}. Second, we have to verify that analogues of~\ref{thm:possametrivial} and~\ref{thm:posdiffinv} hold. For this, let $\alpha, \xi, \zeta \in H_4$ such that $\subsys{\alpha,\xi}$ and $\subsys{\alpha, \zeta}$ are of type $H_2$ and let $w_\xi$, $w_\zeta$ be Weyl elements of type $\xi$ and $\zeta$, respectively. If $\subsys{\alpha, \xi, \zeta}$ is of rank~3, then it must be of type $H_3$ because this is the only type of subsystems of $H_3$ that contains more than one $H_2$-subsystem. If $\subsys{\alpha,\xi, \zeta}$ is of rank~2, then it is of type $H_2$. Hence in any case, $\alpha, \xi, \zeta$ are contained in some $H_3$-subsystem of $\roots$. Thus it follows from~\ref{thm:possametrivial} and~\ref{thm:posdiffinv} that $w_\xi^2$ and $w_\zeta^2$ act identically on $\rootgrX{\alpha}$ if they are in the same position with respect to $\alpha$, and that their actions differ by inversion if one is in involution position and one is inverted involution position with respect to $\alpha$. This observation allows us to define the twisting involution in $H_4$-graded groups as in~\ref{def:invorg}. The assertions of~\ref{lem:invo-twist-weyl} and~\ref{prop:h3-weyl-summary} then remain valid as well.
\end{rem}

\begin{example}\label{ex:twist-is-switch}
	Let $\ring$ be a commutative ring, let $X \in \Set{D_6, E_8}$, put $\ell \defl \rank(X)/2$, let $ G = \Chev(X, \ring) $ be the Chevalley group from~\ref{ex:chev-D6} or~\ref{ex:chev-E8} and denote by $(\risomHex{\alpha})_{\alpha \in H_\ell}$ the family of root isomorphisms from~\ref{rem:ex-twist-signs}. Then a straightforward matrix computation (see~\ref{rem:gap}) shows that for all $ \alpha \in H_\ell $ and all $ x,y \in \ring $, the twisting involution is given by $ \invo{\risomHex{\alpha}(x,y)} = \risomHex{\alpha}(-x,y) $.
\end{example}

\section{The parametrisation theorem}
\label{sec:gpt}

\begin{notation}\label{conv:gpt}
	In this section, we fix a root system $\roots$ with the property that any pair of orthogonal roots in $\roots$ is adjacent and that the Weyl group acts transitively on $\roots$. (In other words, $\roots$ is irreducible of type $A$, $D$, $E$ or $H$.) We denote by $(G, (\rootgr{\alpha})_{\alpha \in \roots})$ a $\roots$-graded group and by $\rootbase$ a root base of $\roots$ . For each $ \delta \in \rootbase $, we fix a $ \delta $-Weyl element $ w_\delta $ and we put $ w_{-\delta} \defl w_\delta^{-1} $. For any word $ \word{\delta} = (\delta_1, \ldots, \delta_k) $ over $ \rootbase \cup (-\rootbase) $, we put $ w_{\word{\delta}} \defl w_{\delta_1} \cdots w_{\delta_k} $ and $ \refl{\word{\delta}} \defl \reflbr{\word{\delta}} \defl \refl{\delta_1} \cdots \refl{\delta_k}$.
\end{notation}

The parametrisation theorem is a general tool which yields the existence of a parametrisation of $G$, that is, of a family of root isomorphisms which is \enquote{compatible} with the fixed family $(w_\delta)_{\delta \in \rootbase}$ of Weyl elements (see~\ref{def:para} for the precise definition).
In this section, we state this theorem for $ \roots $-graded groups and give a self-contained proof (see~\ref{thm:para}). The restrictions on $\roots$ in~\ref{conv:gpt} simplify this task drastically. We will only be interested in the case $\roots \in \Set{H_3, H_4}$, but there are no additional complications by allowing the other cases in~\ref{conv:gpt}. More details, motivation and a complete proof of the parametrisation theorem in full generality can be found in \cite[Chapter~4]{torben}. Combining this theorem with the results from the previous section, we obtain a parametrisation of $ (G, (\rootgr{\alpha})_{\alpha \in \roots}) $ for $\roots \in \Set{H_3, H_4}$, see~\ref{conclusion:para}. This parametrisation will be refined to a standard coordinatisation (in the sense of~\ref{def:stand-coord}) in the following sections using the blueprint technique.

\begin{rem}\label{rem:H3-ortho-triv}
	We will mostly use our assumption on orthogonal roots (see~\ref{conv:gpt}) in the following way: If $\alpha, \beta$ are orthogonal, then any $ \beta $-Weyl element centralises~$ \rootgr{\alpha} $.
\end{rem}

\begin{definition}
	For any word $ \word{\alpha} = (\alpha_1, \ldots, \alpha_m) $ over $ \rootbase \union (-\rootbase) $, we define the corresponding \emph{inverse word} to be $ \word{\alpha}^{-1} \defl (-\alpha_m, \ldots, -\alpha_1) $. For any root $ \alpha $, a \emph{$ \rootbase $-expression of $ \alpha $} is a word over $ \rootbase \union (-\rootbase) $ of the form $ (\word{\delta}^{-1}, \delta, \word{\delta}) $ where $ \delta \in \rootbase $ and $ \word{\delta} $ is a word over $ \rootbase \union (-\rootbase) $ such that $ \delta^{\reflbr{\word{\delta}}} = \alpha $.
\end{definition}

\begin{rem}\label{rem:expr-is-weyl}
	Let $\beta$ be a root. By the transitivity properties of the Weyl group, there exist $\rootbase$-expressions of $\beta$. If $ \word{\beta} $ is a $ \rootbase $-expression of a root $ \beta $, then $ w_{\word{\beta}} $ is a $ \beta $-Weyl element by~\ref{prop:weylelements}~\ref{prop:weylelements:conj}. In contrast, if $\word{\beta}$ is any word over $\rootbase \union (-\rootbase)$ with $\refl{\word{\beta}} = \refl{\beta}$, then $w_{\word{\beta}}$ need not be a $\beta$-Weyl element because it is not clear that $w_{\word{\beta}} \in \rootgr{-\beta} \rootgr{\beta} \rootgr{-\beta}$. This complication is the reason why $\rootbase$-expressions are more practical to work with.
\end{rem}

In order to define parametrisations, we need to first introduce parity maps and parameter systems.

\begin{notation}\label{notation:braidword}
	Let $ A $ be any set. For all $s,t\in A$ and $m\in\N_0$, we denote by $ \braidword{m}(s,t) $ the word $(u_1,\ldots,u_m)$ with $u_k=s$ if $k$ is odd and $u_k=t$ if $k$ is even. We call $\braidword{m}(s,t)$ the \emph{braid word of $(s,t)$ of length $m$}.
\end{notation}

\begin{definition}\label{param:parmap-def}
	Let $ A $ be an abelian group. A \emph{$ \rootbase $-parity map with values in $ A $} is a map $ \map{\inverparsym}{\roots \times \rootbase}{\twistgroup}{(\alpha, \delta)}{\inverpar{\alpha}{\delta} = \inverparsym(\alpha, \delta)} $. Given any $ \rootbase $-parity map $ \inverparsym $, we define $ \inverparbr{\alpha}{-\delta} \defl \inverparbr{\alpha^{\reflbr{\delta}}}{\delta}^{-1} $ for all $ \alpha \in \roots $, $ \delta \in \rootbase $ and for any word $ \word{\delta} = (\delta_1, \ldots, \delta_m) $ over $ \rootbase \union (-\rootbase) $, we define $ \inverpar{\alpha}{\word{\delta}} = \inverparbr{\alpha}{\word{\delta}} $ to be
	\begin{align*}
		\prod_{i=1}^m \inverparbr{\alpha^{\reflbr{\delta_1, \dots, \delta_{i-1}}}}{\delta_i} = \inverparbr{\alpha}{\delta_1} \cdot \inverparbr{\alpha^{\reflbr{\delta_1}}}{\delta_2} \cdots \inverparbr{\alpha^{\reflbr{\delta_1, \dots, \delta_{m-1}}}}{\delta_m},
	\end{align*}
	with the convention that $ \inverpar{\alpha}{\emptyset} = 1_A $ for the empty word $ \emptyset $.
	When $ \rootbase $ and $ A $ are clear from the context, we will simply say that $ \inverparsym $ is a \emph{parity map}. If $ \inverparbr{\alpha}{\braidword{m}(\beta, \gamma)} = \inverparbr{\alpha}{\braidword{m}(\gamma, \beta)} $ for all $ \alpha \in \roots $ and $ \beta, \gamma \in \rootbase $ where $ m = m_{\alpha,\beta} $ denotes the group-theoretic order of $ \refl{\beta} \refl{\gamma} $ in the Weyl group, then $ \inverparsym $ is called \emph{braid-invariant}. If for all roots $ \alpha, \beta \in \roots $ such that $ \alpha $ is adjacent to $ \beta $ and $ -\beta $ and for all $ \rootbase $-expressions $ \word{\beta} $ of $ \beta $, we have $ \inverparbr{\alpha}{\word{\beta}} = 1_\twistgroup $, then $ \inverparsym $ is called \emph{adjacency-trivial}. If for all $ \alpha \in \roots $ and all $ b \in \twistgroup $, there exists a word $ \word{\delta} $ over $ \rootbase $ such that $ \refl{\word{\delta}} $ stabilises $ \alpha $ and $ \inverparbr{\alpha}{\word{\delta}} = b$, then $ \inverparsym $ is called \emph{complete}.
\end{definition}

\begin{example}\label{rem:stand-parmap}
	In~\ref{rem:fold-ex-parmap}, we have computed $\rootbase$-parity maps for $\roots = H_3$ and $\roots=H_4$ with values in $\Set{\pm 1} \times \Set{\pm 1}$. We call them the \emph{standard parity maps for $H_3$ and $H_4$}. The fact that they are \enquote{induced} by certain $\roots$-graded groups ensures that they have many desirable properties (see~\ref{lem:standard-parmap-properties}). Still, there is nothing canonical about them: A twisting procedure as in~\ref{def:param-twist} would allow us to obtain many more parity maps.
	
	Note that the standard parity maps depend on the choice of a root base in $H_3$ or $H_4$. This is clear from the tables giving the values of the parity map, where each root $\alpha$ is expressed as a linear combination of a fixed root base in standard order. In the following, we will always tacitly assume that we choose the standard parity maps with respect to the root base $\rootbase$ fixed in~\ref{conv:gpt}.
\end{example}

\begin{rem}\label{rem:complete}
	By our assumption that the Weyl group acts transitively on $\roots$, a parity map $\inverparsym$ is complete if and only if for \emph{some} $\alpha \in \roots$ and all $ b \in \twistgroup $, there exists a word $ \word{\delta} $ over $ \rootbase $ such that $ \refl{\word{\delta}} $ stabilises $ \alpha $ and $ \inverparbr{\alpha}{\word{\delta}} = b$.
\end{rem}

\begin{rem}\label{rem:adj-ortho-eq}
	Let $ \alpha, \beta \in \roots $. Then $ \alpha $ is adjacent to $ \beta $ and $ -\beta $ if and only if $ \alpha $ is orthogonal to $ \beta $. Here the backwards implication holds by the assumptions in~\ref{conv:gpt}. This observation will be relevant for verifying that a parity map is adjacency-trivial.
\end{rem}

\begin{definition}
	A \emph{(faithful) parameter system} is a pair $ (\twistgroup, M) $ consisting of abelian groups $ \twistgroup, M $ such that $ \twistgroup $ acts (faithfully) on the underlying set of $ M $.
\end{definition}

\begin{definition}
\label{def:para}
Let $(A,M)$ be a parameter system and let $\eta$ be a parity map with values in $A$. A \emph{parametrisation of $(G, (\rootgr{\alpha})_{\alpha \in \roots})$ by $(A,M)$ with respect to $\eta$ and $(w_\delta)_{\delta\in\rootbase}$} is a family $(\risom{\alpha}:M\rightarrow U_\alpha)_{\alpha\in\Phi}$ of group isomorphisms such that for all $\alpha\in\Phi$, $\delta\in\Delta$ and $x\in M$, we have
	\[\theta_\alpha(x)^{w_\delta}=\theta_{\alpha^{\reflbr{\delta}}}(\eta_{\alpha,\delta}.x).\]
The family $(\theta_\alpha)_{\alpha\in\Phi}$ will also be called a \emph{parametrisation of $G$} and the maps in this family are called \emph{root isomorphisms}. We will also say that \emph{$G$ is parametrised by $(A,M)$ with respect to $\eta$, $(w_\delta)_{\delta\in\Delta}$ and with root isomorphisms $(\theta_\alpha)_{\alpha\in\Phi}$}. 
\end{definition}

\begin{example}\label{ex:paramsys}
	Let $ \twistgroup \defl \Set{\pm 1}^2 $, let $ \ring $ be a commutative ring and put $ M \defl (\ring^2, +) $. Let $ \twistgroup $ act on $ M $ via $ (\ep,\ep').(x,y) \defl (\ep x,\ep'y) $ for all $ x,y \in \ring $ and $ \ep,\ep' \in \{\pm 1\} $. Then $ (\twistgroup, M) $ is a parameter system, called the \emph{standard parameter system for $ \ring $}. It is faithful if and only if $ 2_\ring \ne 0_\ring $. We will show in~\ref{thm:prod-param} that every root graded group of type $H_3$ or $H_4$ is parametrised by a parameter system of this form.
\end{example}

The final missing ingredient to formulate the parametrisation theorem is the notion of twisting groups.

\begin{definition}\label{def:twistgrp}
	Let $ \twistgroup $ be an abelian group equipped with an action
	\[ \map{}{\twistgroup \times \rootgr{\alpha}}{\rootgr{\alpha}}{(b,x_\alpha)}{b.x_\alpha} \]
	on the underlying set of $ \rootgr{\alpha} $ for each $ \alpha \in \roots $. Then $ \twistgroup $ is called a \emph{twisting group for $ (G, (w_\delta)_{\delta \in \rootbase}) $} if for all $ \alpha \in \roots $, $ \delta \in \rootbase \union (-\rootbase) $ and $ b \in \twistgroup $, we have $ (b.x_\alpha)^{w_\delta} = b.(x_\alpha^{w_\delta}) $ for all $ x_\alpha \in \rootgr{\alpha} $. We will refer to the actions of $ \twistgroup $ on the root groups as \emph{twisting actions (on the root groups)}.
\end{definition}

\begin{example}
\label{ex:twistgrouph3}
Assume that $\roots \in \Set{H_3, H_4}$ and put $ \twistgroup \defl \Set{\pm 1} \times \Set{\pm 1} $. For all $ \alpha \in \roots $, we declare that $ \twistgroup $ acts on $ \rootgr{\alpha} $ via $ (-1,-1).x_\alpha \defl x_\alpha^{-1} $ and $ (-1,1).x_\alpha \defl \invoroot{x_\alpha}{\alpha} $ for all $ x_\alpha \in \rootgr{\alpha} $ where $\invosymroot{\alpha}$ denotes the twisting involution from~\ref{def:invorg} and~\ref{rem:weyl-h4}. Then by~\ref{lem:invo-twist-weyl} and~\ref{rem:weyl-h4}, $ \twistgroup $ is a twisting group for $ (G, (w_\delta)_{\delta \in \rootbase}) $, called the \emph{standard twisting group for $ (G, (w_\delta)_{\delta \in \rootbase}) $}. We will use this twisting group to construct a parametrisation of $ G $.

Now let, in addition, $\ring$ be a commutative ring, let $X \in \Set{D_6, E_8}$, put $\ell \defl \rank(X)/2$, let $ G = \Chev(X, \ring) $ be the Chevalley group from~\ref{ex:chev-D6} or~\ref{ex:chev-E8} and denote by $(\risomHex{\alpha})_{\alpha \in H_\ell}$ the family of root isomorphisms from~\ref{rem:ex-twist-signs}. Then by~\ref{ex:twist-is-switch} and~\ref{rem:weyl-h4}, $ (\ep, \ep').\risomHex{\alpha}(x,y) = \risomHex{\alpha}\brackets[\big]{(\ep, \ep').(x,y)} $ for all $ x,y \in \ring $ and all $ (\ep, \ep') \in \twistgroup $.
\end{example}

\begin{definition}
	Let $ \inverparsym $ be a $ \rootbase $-parity map with values in a twisting group $ \twistgroup $ for $ (G, (w_\delta)_{\delta \in \rootbase}) $. We say that $ G $ is \emph{square-compatible with respect to $ \inverparsym $ (and $ (w_\delta)_{\delta \in \rootbase} $)} if for all $ \alpha \in H_3 $ and all $ \delta \in \rootbase $, we have $ x_\alpha^w = \inverparbr{\alpha}{(\delta, \delta)}.x_\alpha $ for all $ x_\alpha \in \rootgr{\alpha} $ where $ w \defl w_\delta^2 $.
\end{definition}

We now show that the standard parity maps have all the desirable properties that we have defined in this section.

\begin{prop}\label{lem:standard-parmap-properties}
	Let $X \in \Set{D_6, E_8}$, put $\ell \defl \rank(X)/2$ and let $\map{\inverparsym}{H_\ell \times \rootbase}{\Set{\pm 1} \times \Set{\pm 1}}{}{}$ be the standard parity map for $H_\ell$. Then the following hold:
	\begin{enumerate}[(a)]
		\item \label{lem:standard-parmap-properties:param}Let $\ring$ be a commutative ring, consider the Chevalley group $\Chev(X,\ring)$ with root groups $ (\rootgrHex{\alpha})_{\alpha \in H_\ell} $ from~\ref{ex:chev-D6} or~\ref{ex:chev-E8} and let $ (w_\delta')_{\delta \in \rootbase} $ denote the family of Weyl elements from~\ref{rem:fold-ex-parmap}. Then the family $ (\risomHex{\alpha})_{\alpha \in H_\ell} $ of root isomorphisms from~\ref{rem:ex-twist-signs} is a parametrisation of $ (\Chev(X, \ring), (\rootgrHex{\alpha})_{\alpha \in H_\ell}) $ by the standard parameter system for $ \ring $ (see~\ref{ex:paramsys}) with respect to $ \inverparsym $ and $ (w_\delta')_{\delta \in \rootbase} $.
		
		\item \label{lem:standard-parmap-properties:props}$\inverparsym$ is braid-invariant, adjacency-trivial, complete and satisfies $\inverparbr{\alpha}{\delta} = \inverparbr{-\alpha}{\delta}$ for all $\alpha \in H_\ell$, $\delta \in \rootbase$.
		
		\item \label{lem:standard-parmap-properties:square}Regard $\inverparsym$ as having values in the standard twisting group $\twistgroup=\Set{\pm 1} \times \Set{\pm 1}$ for $(G, (w_\delta)_{\delta \in \rootbase})$ from~\ref{ex:twistgrouph3}. Then $G$ is square-compatible with respect to $\inverparsym$ and $(w_\delta)_{\delta \in \rootbase}$.
	\end{enumerate}
\end{prop}
\begin{proof}
	Assertion~\ref{lem:standard-parmap-properties:param} holds by~\eqref{eq:fold-ex-parmap}, which is the defining property of $ \inverparsym $.
	
	In the remaining part of the proof, we choose any commutative ring $ \ring $ with $ 2_\ring \ne 0_\ring $ and we fix the notation from~\ref{lem:standard-parmap-properties:param}. To show that $ \inverparsym $ is adjacency-trivial, let $ \alpha, \beta \in H_\ell $ such that $ \alpha $ is adjacent to $ \beta $ and $ -\beta $ and let $ \word{\beta} $ be a $ \rootbase $-expression of $ \beta $. Then $ \alpha^{\reflbr{\word{\beta}}} = \alpha^{\reflbr{\beta}} = \alpha $ by~\ref{rem:adj-ortho-eq} and $ w_{\word{\beta}} $ acts trivially on $ \rootgr{\alpha} $ by~\ref{rem:H3-ortho-triv}. Thus by~\ref{lem:standard-parmap-properties:param},
	\begin{align*}
		\risomHex{\alpha}(r,s) = \risomHex{\alpha}(r,s)^{w_{\word{\beta}}} = \risomHex{\alpha}\brackets[\big]{\inverpar{\alpha}{\word{\beta}}.(r,s)}
	\end{align*}
	for all $ r,s \in \ring $. Since the standard parameter system is faithful by~\ref{ex:paramsys} and our assumption that $ 2_\ring \ne 0_\ring $, it follows that $\inverpar{\alpha}{\word{\beta}} = (1,1) $. Hence $ \inverparsym $ is adjacency-trivial. In a similar way, it follows from the braid relations (\ref{braidrel}) that $ \inverparsym $ is braid-invariant. 
	
	Observe that $\inverparsym(\rho_2,(\rho_1,\rho_1))=(-1,-1)$ and $\inverparsym(\rho_2,(\rho_3,\rho_3))=(-1,1)$ by Figure~\ref{fig:parmapex-H3} and that these values generate $\twistgroup$. Hence $\inverparsym$ is complete by~\ref{rem:complete}. Further, we have already observed in~\ref{rem:fold-ex-parmap} that $\inverparbr{\alpha}{\delta} = \inverparbr{-\alpha}{\delta}$ for all $\alpha \in H_\ell$, $\delta \in \rootbase$. This finishes the proof of~\ref{lem:standard-parmap-properties:props}.
	
	In addition to the action of $\twistgroup$ on the root groups of $G$, we let $ \twistgroup $ act on $ \ring^2 $ as in~\ref{ex:paramsys}. Since the group $(\Chev(X, \ring), (\rootgrHex{\alpha})_{\alpha \in H_\ell})$ is parametrised with respect to $\eta$, we have
	\[ \risomHex{\alpha}(x,y)^{w_\delta' w_\delta'}= \risomHex{\alpha}\brackets[\big]{\inverparbr{\alpha}{(\delta, \delta)}.(x,y)} \]
	for all $x,y\in\ring$. At the same time, the action of $ (w_\delta')^2 $ on $ \risomHex{\alpha}(x,y) $ is completely described by~\ref{prop:h3-weyl-summary} and~\ref{ex:twist-is-switch}. Since the parameter system $ (\twistgroup, \ring^2) $ is faithful, we conclude:
	$ \inverparbr{\alpha}{(\delta, \delta)} = (1, 1) $ if $ \subsys{\alpha, \delta} $ is of type $ A_1 $ or $ A_1 \times A_1 $; $ \inverparbr{\alpha}{(\delta, \delta)} = (-1, -1) $ if $ \subsys{\alpha, \delta} $ is of type $ A_2 $; $ \inverparbr{\alpha}{(\delta, \delta)} = (-1, 1) $ if $ \alpha $ is in involution position with respect to $ \delta $ and $ \inverparbr{\alpha}{(\delta, \delta)} = (1,-1) $ if $ \alpha $ is in inverted involution position with respect to $ \alpha $.
	Since~\ref{prop:h3-weyl-summary} applies to $ (G, (\rootgr{\alpha})_{\alpha \in H_\ell}) $ as well, the previous statement and the definition of the action of $ \twistgroup $ on $ (\rootgr{\alpha})_{\alpha \in H_\ell} $ imply that $ G $ is square-compatible with respect to~$\inverparsym$.
\end{proof}

We now state and prove the parametrisation theorem. At first, we need to introduce some terminology which captures the transformations occurring in the solution of the word problem for Coxeter groups.

\begin{definition}
\label{def:homotopy}
	Let $f$, $g$ be two words over $ \rootbase $.
	\begin{enumerate}[(a)]
		\item We say that $ f $ is \emph{elementary braid-homotopic to $ g $} if there exist $ \alpha, \beta \in \rootbase $ and (possibly empty) words $ f_1, f_2 $ over $ \rootbase $ such that $ f = f_1 \braidword{m}(\alpha, \beta) f_2 $ and $ g = f_1 \braidword{m}(\beta, \alpha) f_2 $ where $ m $ denotes the group-theoretic order of $ \refl{\alpha} \refl{\beta} $ in the Weyl group and $ \braidword{m} $ is defined as in~\ref{notation:braidword}. 
		
		\item We say that $ f $ is \emph{elementary square-homotopic to $ g $} if there exist $ \alpha \in \rootbase $ and (possibly empty) words $ f_1, f_2 $ over $ \rootbase $ such that $ f = f_1 (\alpha, \alpha) f_2 $ and $ g = f_1 f_2 $, or such that $f = f_1 f_2$ and $g = f_1 (\alpha, \alpha) f_2$. We say that $ f $ and $ g $ are \emph{square-homotopic} if they lie in the equivalence relation generated by the relation \enquote{elementary braid-homotopic}.
		
		\item We say that $f$ and $g$ are \emph{homotopic} if they lie in the equivalence relation generated by the relation \enquote{braid-homotopic or square-homotopic}.
	\end{enumerate}
\end{definition}

\begin{lemma}[{\cite[4.5.3, 4.5.4]{torben}}]\label{lem:homotopic-weyl}
	Let $\inverparsym$ be a braid-invariant $ \rootbase $-parity map with values in a twisting group $ \twistgroup $ for $ (G, (w_\delta)_{\delta \in \rootbase}) $ such that $ G $ is square-compatible with respect to $ \inverparsym $. Let $\word{\beta}, \word{\beta}'$ be two words over $\rootbase \cup (-\rootbase)$ such that $\reflbr{\word{\beta}} = \reflbr{\word{\beta}'}$. Then $w_{\word{\beta}}$ acts on $\rootgr{\alpha}$ by $\inverparbr{\alpha}{\word{\beta}}$ if and only if $w_{\word{\beta}'}$ acts on $\rootgr{\alpha}$ by $\inverparbr{\alpha}{\word{\beta}'}$.
\end{lemma}
\begin{proof}
	Since $ G $ is square-compatible with respect to $ \inverparsym $, we may replace every letter of $\word{\beta}$ and $\word{\beta}'$ in $-\rootbase$ by its negative. Hence we may assume that $\word{\beta}, \word{\beta'}$ are words over $\rootbase$. Since $\reflbr{\word{\beta}} = \reflbr{\word{\beta}'}$, the words $\word{\beta}$ and $\word{\beta'}$ are then homotopic. If they are square-homotopic, then the assertion holds because $ G $ is square-compatible with respect to $ \inverparsym $. If they are braid-homotopic, the assertion holds by the braid relations for Weyl elements (see~\ref{braidrel}) and because $\inverparsym$ is braid-invariant. It follows that the assertion holds in this special case. See \cite[4.5.3, 4.5.4]{torben} for a more detailed proof.
\end{proof}

\begin{lemma}[{\cite[4.5.5]{torben}}]\label{lem:weyl-homotopic}
	Let $ \inverparsym $ be an adjacency-trivial, braid-invariant $ \rootbase $-parity map with values in a twisting group $ \twistgroup $ for $ (G, (w_\delta)_{\delta \in \rootbase}) $ such that $ G $ is square-compatible with respect to $ \inverparsym $. Let $ \alpha \in \roots $ and let $ \word{\beta} $ be a word over $\rootbase \cup (-\rootbase)$ such that $\alpha^{\reflbr{\word{\beta}}} = \alpha$. Then $ x_\alpha^{w_{\word{\beta}}} = \inverpar{\alpha}{\word{\beta}}.x_\alpha $ for all $ x_\alpha \in \rootgr{\alpha} $.
\end{lemma}
\begin{proof}
	Since $\alpha^{\reflbr{\word{\beta}}} = \alpha$, there exist $ \rho_1, \ldots, \rho_n \in \roots \cap \alpha^\perp $ with $ \refl{\word{\beta}} = \refl{\rho_1} \cdots \refl{\rho_n} $. By~\ref{lem:homotopic-weyl}, we may replace $\word{\beta}$ by any word $\word{\gamma}$ with $\reflbr{\word{\beta}} = \reflbr{\word{\gamma}}$, so we may assume that $\word{\beta} = (\word{\gamma}^1, \ldots, \word{\gamma}^n)$ where $\word{\gamma}^1, \ldots, \word{\gamma}^n$ are $\rootbase$-expressions of $\rho_1, \ldots, \rho_n$, respectively. For all $i \in \Set{1,\ldots, n}$, $\inverparbr{\alpha}{\word{\gamma}^i} = 1_\twistgroup$ because $\inverparsym$ is adjacency-trivial and further, $w_{\word{\gamma}^i}$ centralises $\rootgr{\alpha}$ by~\ref{rem:expr-is-weyl} and~\ref{rem:H3-ortho-triv}. Hence $w_{\word{\beta}}$ centralises $\rootgr{\alpha}$ and $\inverparbr{\alpha}{\word{\beta}} = 1_\twistgroup$, so the assertion holds.
\end{proof}

\begin{note}
	The roots $\rho_1, \ldots, \rho_n$ in the proof of~\ref{lem:weyl-homotopic} need not be contained in $\rootbase$, so we may not write $w_{\rho_i}$ and $\inverparbr{\alpha}{\rho_i}$.
\end{note}

\begin{prop}[{\cite[4.5.6]{torben}}]\label{prop:param-prop}
	Let $ \inverparsym $ be an adjacency-trivial, braid-in\-vari\-ant $ \rootbase $-parity map with values in a twisting group $ \twistgroup $ for $ (G, (w_\delta)_{\delta \in \rootbase}) $ such that $ G $ is square-compatible with respect to $ \inverparsym $. Let $ \alpha \in \roots $ and let $ \word{\beta} $, $ \word{\gamma} $ be two words over $ \rootbase \union (-\rootbase) $ such that $ \alpha^{\reflbr{\word{\beta}}} = \alpha^{\reflbr{\word{\gamma}}} $.
	Put $ b \defl \inverpar{\alpha}{\word{\beta}} \inverpar{\alpha}{\word{\gamma}}^{-1} $. Then $ x_\alpha^{w_{\word{\beta}}} = (b.x_\alpha)^{w_{\word{\gamma}}} $ for all $ x_\alpha \in \rootgr{\alpha} $.
\end{prop}
\begin{proof}
	This follows from~\ref{lem:weyl-homotopic}. See \cite[4.5.6]{torben} for details.
\end{proof}

The following result is precisely the content of the parametrisation theorem \cite[4.5.16]{torben} for root systems $\roots$ satisfying the restrictions in~\ref{conv:gpt}.

\begin{theorem}[{\cite[4.5.16]{torben}}]
\label{thm:para}
	Let $ \inverparsym $ be a complete, adjacency-trivial, braid-invariant $ \rootbase $-parity map with values in a twisting group $ \twistgroup $ for $ (G, (w_\delta)_{\delta \in \rootbase}) $ such that $ G $ is square-compatible with respect to $ \inverparsym $. Then there exist an abelian group $ (\prodring, +) $, an action $ \map{\omega}{\twistgroup \times \prodring}{\prodring}{(b,r)}{b.r} $ of $ \twistgroup $ on the underlying set of $ \prodring $ and a parametrisation $ (\map{\risom{\alpha}}{\prodring}{\rootgr{\alpha}}{}{})_{\alpha \in H_3} $ of $ G $ by $ (\twistgroup, \prodring) $ with respect to $ \inverparsym $ and $ (w_\delta)_{\delta \in \Delta} $ such that the twisting actions of $ \twistgroup $ on the root groups are induced by $ \omega $. That is, we have $ b.\risom{\alpha}(s) = \risom{\alpha}(b.s) $ for all $ \alpha \in H_3 $, $ s \in \prodring $ and $ b \in \twistgroup $.
\end{theorem}
\begin{proof}[Sketch of the proof]
	Fix a root $ \alpha_0 \in H_3 $. Choose a group $ \prodring $ which is isomorphic to $ \rootgr{\alpha_0} $ and an isomorphism $ \map{\risom{\alpha_0}}{\prodring}{\rootgr{\alpha_0}}{}{} $. We define an action of $ \twistgroup $ on $ \prodring $ by $ b.s \defl \risom{\alpha_0}^{-1}(b.\risom{\alpha_0}(s)) $ for all $ b \in \twistgroup $ and $ s \in \prodring $. For any root $ \beta $, we define $ \map{\risom{\beta}}{\prodring}{\rootgr{\beta}}{s}{\risom{\alpha_0}(s)^{w_{\word{\gamma}}}} $ where $ \word{\gamma} $ is any word over $ \rootbase \union (-\rootbase) $ with $ \alpha_0^{\reflbr{\word{\gamma}}} = \beta $ and $ \inverparbr{\alpha}{\word{\gamma}} = 1_\twistgroup $. Such a word $ \word{\gamma} $ exists by the completeness of $ \inverparsym $, and $ \risom{\beta} $ does not depend on the choice of $ \word{\gamma} $ by~\ref{prop:param-prop}. One can verify that this construction has the desired properties.
\end{proof}

\begin{notation}
	For the rest of this section, we assume that $ \roots = H_\ell $ for some $ \ell \in \Set{3,4} $, that $ X $ is defined as in~\ref{not:root-base-order} and that $ \rootbase $ is in standard order. We write $ \rootbase = (\rho_1, \rho_2, \rho_3) $ if $ \ell = 3 $ and $ \rootbase = (\rho_0, \rho_1, \rho_2, \rho_3) $ if $ \ell=4 $.
\end{notation}

\begin{conclusion}
\label{conclusion:para}
	Denote by $ \twistgroup \defl \{\pm 1\}^2 $ the standard twisting group for $(G, (w_\delta)_{\delta \in \rootbase})$ from~\ref{ex:twistgrouph3} and by $ \map{\inverparsym}{H_\ell}{\twistgroup}{}{} $ the standard parity map from~\ref{rem:stand-parmap}. By~\ref{lem:standard-parmap-properties}, $(G, (\rootgr{\alpha})_{\alpha \in H_\ell})$ and $\inverparsym$ satisfy the assumptions in~\ref{thm:para}. Hence there exist an abelian group $ (\prodring, +) $ on which $ \twistgroup $ acts and a parametrisation $ (\risom{\alpha})_{\alpha \in H_\ell} $ of $ G $ by $ (\twistgroup, \prodring) $ with respect to $ \inverparsym $ and $ (w_\delta)_{\delta \in \rootbase} $. Putting $ \invo{s} \defl (-1,1).s $ for all $ s \in \prodring $, we can summarise the property of being a parametrisation as follows:
	\begin{equation}\label{eq:H3-param}
		\risom{\alpha}(s)^{w_\delta} = \begin{cases}
			\risom{\alpha^{\reflbr{\delta}}}(s) & \text{if } \inverpar{\alpha}{\delta} = (1,1), \\
			\risom{\alpha^{\reflbr{\delta}}}(-s) & \text{if } \inverpar{\alpha}{\delta} = (-1,-1), \\
			\risom{\alpha^{\reflbr{\delta}}}(\invo{s}) & \text{if } \inverpar{\alpha}{\delta} = (-1,1), \\
			\risom{\alpha^{\reflbr{\delta}}}(-\invo{s}) & \text{if } \inverpar{\alpha}{\delta} = (1,-1)
		\end{cases}
	\end{equation}
	for all $ \alpha \in H_\ell $, $ \delta \in \rootbase $ and $ s \in \prodring $. Note that $ \invo{(\invo{s})} = s $ for all $ s \in \prodring $ because $ (-1,1)^2 = (1,1) $. Further, we define a multiplication on $ \prodring $ by
	\begin{equation}\label{eq:def-mult}
		a\cdot b\coloneqq \theta_{\rho_1+\rho_2}\inv([\theta_{\rho_1}(a),\theta_{\rho_2}(b)]).
	\end{equation}
	It satisfies the distributive laws by~\ref{lem:oneadditivity} because $ (\rho_1, \rho_2) $ is an $ A_2 $-pair. Finally, we define $ 1_\prodring $ to be the element of $ \prodring $ which satisfies $ w_{\rho_2} \in \rootgr{-\rho_2} \risom{\rho_2}(1_\prodring) \rootgr{-\rho_2} $, which is uniquely determined by~\ref{prop:a2gradprop}~\ref{prop:a2gradprop:factor-unique}. By~\ref{prop:a2gradprop}~\ref{prop:a2gradprop:conj-formula}, we have
	\begin{align*}
		\risom{\rho_1 + \rho_2}(s \cdot 1_\prodring) &= \commutator{\risom{\rho_1}(s)}{\risom{\rho_2}(1_\prodring)} = \risom{\rho_1}(s)^{w_{\rho_2}} = \risom{\rho_1 + \rho_2}(\inverpar{\rho_1}{\rho_2}.s) = \risom{\rho_1 + \rho_2}(s)
	\end{align*}
	for all $ s \in \prodring $. A similar computation (see~\cite[5.7.12]{torben}) shows that $ 1_\prodring \cdot r = r $ for all $ r \in \prodring $ as well, so that $ (\prodring, +, \cdot) $ is a nonassociative ring with unit $ 1_\prodring $.
\end{conclusion}

\begin{goal}\label{goal:ring-decomp}
	We will show in Section~\ref{sec:bph3} that the multiplication on $ \prodring $ from~\ref{conclusion:para} is associative (\ref{thm:associative}) and commutative (\ref{prop:commutative}), that $ \prodring $ decomposes as $ \prodring = \prodring_1 \oplus \prodring_2 $ for a pair of isomorphic ideals $ \prodring_1 $, $ \prodring_2 $ of $ \prodring $ (\ref{prop:dirsum}, \ref{prop:R1R2iso}) and that $ \invo{(r,s)} = (-r,s) $ for all $ (r,s) \in \prodring_1 \oplus \prodring_2 $ (\ref{thm:invocomp}). Further, we will show that the commutator relations in $ G $ are precisely those given in Figure~\ref{fig:excommrel} (\ref{prop:comm-formula}, \ref{prop:comm-formula-H4}). This says precisely that the parametrisation in~\ref{conclusion:para} is a standard coordinatisation in the following sense.
\end{goal}

\begin{definition}\label{def:stand-coord}
	Let $\ring$ be a commutative ring. A \emph{standard coordinatisation of $ (G, (\rootgr{\alpha})_{\alpha \in H_\ell}) $ by $\ring \times \ring$} (or, in the terminology of \cite[4.7.5]{torben}, a \emph{coordinatisation with standard signs}) is a family $(\map{\risom{\alpha}}{(\ring \times \ring,+)}{\rootgr{\alpha}}{}{})_{\alpha \in H_\ell}$ of root isomorphisms with the following properties.
	\begin{enumerate}[(i)]
		\item \label{def:stand-coord:weyl}$ w_\delta' \defl \risom{-\delta}(-1_\ring, -1_\ring) \risom{\delta}(1_\ring, 1_\ring) \risom{-\delta}(-1_\ring, -1_\ring) $ is a $ \delta $-Weyl element for all $ \delta \in \rootbase $.
		
		\item \label{def:stand-coord:param}$(\risom{\alpha})_{\alpha \in H_\ell}$ is a parametrisation of $(G, (\rootgr{\alpha})_{\alpha \in H_\ell})$ (in the sense of~\ref{def:para}) by the standard parameter system $(\ring \times \ring, \Set{\pm 1} \times \Set{\pm 1})$ from~\ref{ex:paramsys} with respect to the standard parity map from~\ref{rem:stand-parmap} and $(w_\delta')_{\delta \in \rootbase}$.
		
		\item \label{def:stand-coord:comm}The same commutator relations as in~\ref{lem:excommrel} are satisfied.
	\end{enumerate}
	We also say that $ (\risom{\alpha})_{\alpha \in H_\ell} $ is a \emph{standard coordinatisation with respect to $ (w_\delta')_{\delta \in \rootbase} $}.
\end{definition}

\begin{example}
	The families $(\risomH{\alpha})_{\alpha \in H_3}$ and $(\risomH{\alpha})_{\alpha \in H_\ell}$ from~\ref{rem:ex-twist-signs} are standard coordinatisations of $\Chev(D_6, \ring)$ and $\Chev(E_8, \ring)$, respectively. Here $\ring$ is any commutative ring.
\end{example}

\section{Standard coordinatisations and the Steinberg presentation}
\label{sec:stein}


\begin{notation}
	In this section, we fix $ \ell \in \{3,4\} $ and a commutative ring $ \ring $. We put $ X \defl D_6 $ if $ \ell=3 $ and $ X \defl E_8 $ if $ \ell=4 $.
\end{notation}

The purpose of this section is to collect some properties of $ H_\ell $-graded groups with a standard coordinatisation. They are not needed for the proof of our main result~\ref{thm:prod-param}, by which every $ H_\ell $-graded group has a standard coordinatisation.

The following remark is analogous to \cite[10.4.16]{torben}, which is about $F_4$-gradings.

\begin{rem}\label{rem:conjugate-formulas}
	Let $ (G, (\rootgr{\alpha})_{\alpha \in H_\ell}) $ be an $ H_\ell $-graded group with a standard coordinatisation $ (\risom{\alpha})_{\alpha \in H_\ell} $ by $ \ring \times \ring $.
	Recall from~\ref{def:commmap} that for all non-proportional $\alpha, \gamma \in H_\ell$ and all $\beta \in \oprootint{\alpha,\gamma}$, we have a commutation map $\map{\commmap{\alpha,\gamma}{\beta}}{M \times M}{M}{}{}$ such that
	\begin{equation*}
		\commutator{\risom{\alpha}(r)}{\risom{\gamma}(s)} = \prod_{\beta \in \oprootint{\alpha,\beta}} \risom{\beta}\brackets[\big]{\commmap{\alpha,\gamma}{\beta}(r,s)}
	\end{equation*}
	for all $r,s \in M$ where the product is taken in the (unique) interval ordering on $\oprootint{\alpha,\gamma}$ starting from $\alpha$. Axiom~\ref{def:stand-coord}~\ref{def:stand-coord:comm} says precisely that the commutation map $\commmap{\alpha,\gamma}{\beta}$ is prescribed for all $\alpha,\gamma \in \paratwopos(\rootbase)$ and $\beta \in \oprootint{\alpha,\gamma}$, but no explicit statement is made about the case $ \alpha,\gamma \notin \paratwopos(\rootbase) $. However, let $ (\xi, \zeta) $ be an arbitrary pair of non-proportional roots in $ H_\ell $. Then by~\ref{rem:paratwopos}, there exists $ u \in \Weyl(H_\ell) $ such that $ \alpha \defl \xi^u$ and $\gamma \defl \zeta^u $ lie in $\paratwopos(\rootbaseH)$. Choose $\delta_1, \ldots, \delta_k \in \rootbase$ such that $u^{-1} = \reflbr{\delta_1} \cdots \reflbr{\delta_k}$ and put $w \defl w_{\delta_1}' \cdots w_{\delta_k}'$ where $w_{\delta_1}', \ldots, w_{\delta_k}'$ are the Weyl elements from Axiom~\ref{def:stand-coord}~\ref{def:stand-coord:weyl}. Then it follows from the equation above that
	\[ \commutator{\risom{\alpha}(r)^w}{\risom{\gamma}(s)^w} = \prod_{\beta \in \oprootint{\alpha,\beta}} \risom{\beta}\brackets[\big]{\commmap{\alpha,\gamma}{\beta}(r,s)}^w \]
	is the commutator formula for the pair $(\rootgr{\xi}, \rootgr{\zeta})$. By Axiom~\ref{def:stand-coord}~\ref{def:stand-coord:param}, it looks just like the commutator formula for $(\rootgr{\alpha}, \rootgr{\gamma})$, except that the action of the twisting group may introduce or change some signs. The precise signs that appear in the new formula can be computed using the values of the standard parity map, which are given in Figures~\ref{fig:parmapex-H3},~\ref{fig:parmapex-H4-1} and~\ref{fig:parmapex-H4-2}. (We will see examples of such computations in~\ref{prop:comm-formula-H4} and~\ref{lem:com-form-neg}.) It follows that the commutation maps $\commmap{\xi, \zeta}{\beta}$ with respect to a standard coordinatisation are prescribed for all non-proportional $\xi, \zeta \in H_\ell$ and all $\beta \in \oprootint{\xi,\zeta}$ and that they may be computed explicitly. They are called the \emph{standard commutation maps (with respect to $M=\ring \times \ring$)}. Note that they are maps $\map{}{M \times M}{M}{}{}$ which depend only on $\ring$ but not on $G$.
\end{rem}

We now construct a universal $ H_\ell $-graded group for every commutative ring. Every $ H_\ell $-graded group is the quotient of one such group by~\ref{lem:stein-epi} and our main result~\ref{thm:prod-param}.

\begin{definition}\label{def:stein}
	The \emph{Steinberg group of type $H_\ell$ over $\ring$} is the abstract group $\Stein_{H_\ell}(\ring)$ which is generated by elements $\risomst{\alpha}(r,s)$ for $\alpha \in H_\ell$ and $r,s \in \ring$ with respect to the following relations:
	\begin{enumerate}[(a)]
		\item $\risomst{\alpha}(r+r',s+s') = \risomst{\alpha}(r,s) \risomst{\alpha}(r',s')$ for all $\alpha \in H_\ell$ and all $r,r',s,s' \in \ring$.
		
		\item For all non-proportional $\alpha, \gamma \in H_\ell$, we require the commutator relation
		\[ \commutator{\risomst{\alpha}(r)}{\risomst{\gamma}(s)} = \prod_{\beta \in \oprootint{\alpha,\beta}} \risomst{\beta}\brackets[\big]{\commmap{\alpha,\gamma}{\beta}(r,s)} \]
		for all $r,s \in \ring \times \ring$ where the product is taken in the (unique) interval ordering on $\oprootint{\alpha,\gamma}$ starting from $\alpha$ and where
		\[ \map{\commmap{\alpha,\gamma}{\beta}}{(\ring \times \ring) \times (\ring \times \ring)}{\ring \times \ring}{}{} \]
		denotes the standard commutation map with respect to $\ring \times \ring$ (see~\ref{rem:conjugate-formulas}).
	\end{enumerate}
	For all $ \alpha \in H_\ell $, we define a subgroup $\rootgrst{\alpha} \defl \Set{\risomst{\alpha}(r) \mid r \in \ring \times \ring}$ and we consider $ \risomst{\alpha} $ as a surjective homomorphism from $ \ring \times \ring $ to $ \rootgrst{\alpha} $.
\end{definition}

\begin{rem}
	The Steinberg group is functorial in $ \ring $: If $ \map{\phi}{\ring}{\ring'}{}{} $ is a homomorphism of commutative rings, then there is a unique homomorphism $ \map{\Stein(\varphi)}{\Stein_{H_\ell}(\ring)}{\Stein_{H_\ell}(\ring')}{}{} $ with $ \Stein(\varphi) \circ \risomst{\alpha,\ring} = \risomst{\alpha,\ring'} $ for all $ \alpha \in H_\ell $.
\end{rem}

\begin{lemma}\label{lem:stein-epi}
	Let $ (G, (\rootgr{\alpha})_{\alpha \in H_\ell}) $ be an $ H_\ell $-graded group with a standard coordinatisation $ (\risom{\alpha})_{\alpha \in H_\ell} $ by $ \ring \times \ring $.
	Then there exists a unique surjective homomorphism $ \map{\pi}{\Stein_{H_\ell}(\ring)}{G}{}{} $ such that $ \pi(\risomst{\alpha}(r,s)) = \risom{\alpha}(r,s) $ for all $ \alpha \in H_\ell $ and all $ r,s \in \ring $. Further, for any positive system $ \possys $ in $ H_\ell $, the restriction $ \map{\pi}{\rootgrst{\possys}}{\rootgr{\possys}}{}{} $ is bijective.
\end{lemma}
\begin{proof}
	The existence and uniqueness of $ \pi $ follows from the universal property of groups defined by generators and relations. Now let $ \possys $ be a positive system in $ H_\ell $ and let $ \alpha_1, \ldots, \alpha_k $ be an ordering of the roots in $ \possys $ with the properties in~\ref{prop:prodmapbij:possys}. Let $ x \in \rootgrst{\possys} $ such that $ \pi(x) = 1 $. Since $ \Stein_{H_\ell}(\ring) $ has $ H_\ell $-commutator relations by definition, it follows from~\ref{prop:prodmapbij:possys}~\ref{prop:prodmapbij:possys:surj} that there exist $ r_1, \ldots, r_k \in \ring \times \ring $ with $ x = \prod_{i=1}^k \risomst{\alpha_i}(r_i) $. Then $ 1 = \pi(x) = \prod_{i=1}^k \risom{\alpha_i}(r_i) $, which by~\ref{prop:prodmapbij:possys}~\ref{prop:prodmapbij:possys:bij} implies that $ r_i = 0 $ for all $ i \in \{1, \ldots, k\} $. Hence $ x=1 $, so $ \map{\pi}{\rootgrst{\possys}}{\rootgr{\possys}}{}{} $ is bijective.
\end{proof}

\begin{lemma}\label{lem:stein-possys}
	The maps $(\risomst{\alpha})_{\alpha \in H_\ell}$ are injective and $(\Stein_{H_\ell}(\ring), (\rootgrst{\alpha})_{\alpha \in H_\ell})$ satisfies Axiom~\ref{def:rgg}~\ref{def:rgg:pos}.
\end{lemma}
\begin{proof}
	Choose an $ H_\ell $-graded group $ (G, (\rootgr{\alpha})_{\alpha \in H_\ell}) $ with a standard coordinatisation $ (\risom{\alpha})_{\alpha \in H_\ell} $ by $ \ring \times \ring $. (For example, one of the Chevalley groups from Section~\ref{sec:fold}.) Let $ \map{\pi}{\Stein_{H_\ell}(\ring)}{G}{}{} $ be the homomorphism from~\ref{lem:stein-epi}. Let $ \possys $ be a positive system in $ H_\ell $, let $ \alpha \in -\possys $ and let $ x \in \rootgrst{\alpha} \cap \rootgrst{\possys} $. Then $ \pi(x) \in \rootgr{\alpha} \cap \rootgr{\possys} = \{1_G\} $. As $ \pi $ is injective on $ \rootgrst{\possys} $ by~\ref{lem:stein-epi}, it follows that $ x=1 $. Hence $ \Stein_{H_\ell}(\ring) $ satisfies Axiom~\ref{def:rgg}~\ref{def:rgg:pos}. Further, for each $ \alpha \in H_\ell $, the injectivity of $ \risom{\alpha} $ implies the injectivity of $ \risomst{\alpha} $ because $ \risom{\alpha} = \pi \circ \risomst{\alpha} $.
\end{proof}

The existence of Weyl elements in $ \Stein_{H_\ell}(\ring) $ essentially follows from the existence of Weyl elements in some $ H_\ell $-graded group with a standard coordinatisation.

\begin{lemma}\label{lem:stein-weyl}
	Let $ \alpha \in H_\ell $ and let $ r,s \in \ring $ be invertible. Then
	\[ w_\alpha^{\text{St}}(r,s) \defl \risomst{-\alpha}(-r^{-1}, -s^{-1}) \risomst{\alpha}(r,s) \risomst{-\alpha}(-r^{-1}, -s^{-1}) \]
	is an $ \alpha $-Weyl element in $ \Stein_{H_\ell}(\ring) $.
\end{lemma}
\begin{proof}
	Let $ \beta \in H_\ell $ and $ w_\alpha \defl w_\alpha^{\text{St}}(r,s) $. We have to show that $ (\rootgrst{\beta})^{w_\alpha} = \rootgrst{\beta^{\reflbr{\beta}}} $. Since $ r,s $ are arbitrary and $ w_\alpha^{-1} = w_\alpha(-r,-s) $, it suffices to show that $ (\rootgrst{\beta})^{w_\alpha} \subseteq \rootgrst{\beta^{\reflbr{\alpha}}} $. First, assume that $ \beta \notin \Set{\alpha, -\alpha} $. Let $ x_\beta \in \rootgrst{\beta} $. Let $ \possys $ be the unique positive system in the rank-2 subsystem $ \subsys{\alpha, \beta} $ such that $ \alpha,\beta \in \possys $ and such that the root base of $ \subsys{\alpha,\beta} $ associated to $ \possys $ contains $ \alpha $.
	Put $ P \defl \possys \setminus \{\alpha\} $. Then for all $ g \in \rootgrst{\alpha} \cup \rootgrst{-\alpha} $ and all $ \gamma \in P $, we have $ (\rootgrst{\gamma})^g \subseteq \commutator{g}{\rootgrst{\gamma}} \rootgrst{\gamma} \subseteq \rootgrst{P} $ by~\ref{rem:commrel}~\ref{rem:commrel:conj}. Hence
	\[ (\rootgrst{P})^{w_\alpha} \subseteq \rootgrst{P} \subseteq \rootgrst{\possys}. \]
	Now consider the Chevalley group $ (\Chev(X, \ring), (\rootgrHex{\alpha})_{\alpha \in X}) $ from Section~\ref{sec:fold} and the surjective homomorphism $ \map{\pi}{\Stein_{H_\ell}(\ring)}{\Chev(X,\ring)}{}{} $ from~\ref{lem:stein-epi}. Since $ \pi(w_\alpha) $ is an $ \alpha $-Weyl element in $ \Chev(X, \ring) $ by \ref{rem:fold-ex-parmap}, we have $ \pi(x_\beta^{w_\alpha}) \in \rootgrHex{\beta^{\reflbr{\alpha}}} $. Since the restriction $ \map{\pi}{\rootgrst{\possys}}{\rootgrHex{\possys}}{}{} $ is an isomorphism and $ x_\beta^{w_\alpha} \in (\rootgrst{P})^{w_\alpha} \subseteq \rootgrst{\possys} $, we infer that $ x_\beta^{w_\alpha} \in \rootgrst{\beta^{\reflbr{\alpha}}} $, as desired.
	
	Now assume that $ \beta = \alpha $. Choose an $ A_2 $-pair $ (\gamma, \gamma') $ with $ \gamma+\gamma'=\alpha $. Recall from~\ref{rem:conjugate-formulas} that up to signs, the commutation map $ \commmap{\gamma,\gamma'}{\alpha} $ looks like the map $ \commmap{\rho_0,\rho_1}{} $ in Figure~\ref{fig:excommrel}. Hence $ \commutator{\rootgrst{\gamma}}{\rootgrst{\gamma'}} = \rootgrst{\alpha} $. It follows that
	\[ (\rootgrst{\alpha})^{w_\alpha} = \commutator{(\rootgrst{\gamma})^{w_\alpha}}{(\rootgrst{\gamma'})^{w_\alpha}} \subseteq \commutator{\rootgrst{-\gamma'}}{\rootgrst{-\gamma}} = \rootgrst{-\alpha}. \]
	The case $\beta = -\alpha$ can be treated similarly. Thus $w_\alpha^{\text{St}}(r,s)$ is a Weyl element.
\end{proof}

\begin{prop}\label{prop:stein-is-rgg}
	$ (\Stein_{H_\ell}(\ring), (\rootgrst{\alpha})_{\alpha \in H_\ell}) $ is an $ H_\ell $-graded group with standard coordinatisation $ (\risomst{\alpha})_{\alpha \in H_\ell} $.
\end{prop}
\begin{proof}
	By~\ref{def:stein}, \ref{lem:stein-possys} and~\ref{lem:stein-weyl}, $(\rootgrst{\alpha})_{\alpha \in H_\ell}$ is an $H_\ell$-grading of $\Stein_{H_\ell}(\ring)$ with root isomorphisms $(\risomst{\alpha})_{\alpha \in H_\ell}$. By similar arguments as in~\ref{lem:stein-weyl} (that is, projecting into the appropriate Chevalley group and using that the projection is injective on positive systems), the action of the Weyl elements on the root groups in $\Stein_{H_\ell}(\ring)$ satisfies the same formulas as the corresponding Weyl elements in Chevalley groups. Hence Axiom~\ref{def:stand-coord}~\ref{def:stand-coord:param} is satisfied. Thus $ (\risomst{\alpha})_{\alpha \in H_\ell} $ is a standard coordinatisation.
\end{proof}

\begin{prop}\label{pro:weyl-inv}
	Let $ (G, (\rootgr{\alpha})_{\alpha \in H_\ell}) $ be an $ H_\ell $-graded group with a standard coordinatisation $ (\risom{\beta})_{\beta \in H_\ell} $ by $ \ring \times \ring $ and let $\alpha \in H_\ell$. Then the map
	\[ \map{}{\invring{\ring} \times \invring{\ring}}{\weylset{\alpha}}{(r,s)}{w_\alpha(r,s) \defl \risom{-\alpha}(-r^{-1}, -s^{-1}) \risom{\alpha}(r,s) \risom{-\alpha}(-r^{-1}, -s^{-1})} \]
	is a bijection. Here $ \invring{\ring} $ and $ \weylset{\alpha} $ denote the sets of invertible elements in $ R $ and of $ \alpha $-Weyl elements in $ G $, respectively.
\end{prop}
\begin{proof}
	A similar assertion is proven in~\cite[5.6.6]{torben} for $\roots$-graded groups where $\roots$ is of type $A$, $D$ or $E$, and the same kinds of arguments apply in our situation. For the reader's convenience, we repeat the main ideas. Note that for all $ r,s \in \invring{\ring} $, $ w_\alpha^{\text{St}}(r,s) $ is a Weyl element in $ \Stein_{H_\ell}(\ring) $ by~\ref{lem:stein-weyl} and hence its image $ w_\alpha(r,s) $ in $ G $ is a Weyl element in $ G $.
	
	
	At first, let $\alpha \in H_\ell$ and let $w_\alpha$ be any $\alpha$-Weyl element in $G$. By~\ref{prop:a2gradprop}~\ref{prop:a2gradprop:factor-unique} there exist unique $a,b,c,d \in \ring$ such that
	\[ w_\alpha = \risom{-\alpha}(a,b) \risom{\alpha}(c,d) \risom{-\alpha}(a,b). \]
	Choose a root $\gamma \in H_\ell$ such that $(\alpha, \gamma)$ is an $A_2$-pair. Then by the commutator relations in $\Stein_{H_\ell}(\ring)$, there exist signs $\ep_1, \ep_2 \in \Set{\pm 1}$ such that
	\[ \commutator{\risom{\gamma}(r_1, r_2)}{\risom{\alpha}(s_1, s_2)} = \risom{\alpha+\gamma}(\ep_1 r_1 s_1, \ep_2 r_2 s_2) \]
	for all $r_1, r_2, s_1, s_2 \in \ring$. By an application of~\ref{prop:a2gradprop}~\ref{prop:a2gradprop:conj-formula}, we infer that
	\[ \risom{\gamma}(r,s)^{w_\alpha} = \commutator{\risom{\gamma}(r,s)}{\risom{\alpha}(c,d)} = \risom{\alpha+\gamma}(\ep_1 rc, \ep_2 sd) \]
	for all $r,s \in \ring$. Since the map $\map{}{\ring \times \ring}{\rootgr{\alpha+\gamma}}{(r,s)}{\risom{\gamma}(r,s)^{w_\alpha}}$ is bijective, it follows that $c$ and $d$ are invertible. Similar arguments applied to the formula
	\[ \risom{\gamma}(1_\ring,1_\ring) = \commutator[\big]{\commutator{\risom{\gamma}(1_\ring,1_\ring)}{\risom{\alpha}(c,d)}}{\risom{-\alpha}(a,b)}^{-1} \]
	from~\ref{prop:a2gradprop}~\ref{prop:a2gradprop:conj-formula} show that $c^{-1} = -a$ and $d^{-1}=-b$. (The precise computation, which we skip, relies on the knowledge of the signs $\ep_1$, $\ep_2$.)
	This shows that the map in the assertion is surjective. By~\ref{prop:a2gradprop}~\ref{prop:a2gradprop:factor-unique}, it is also injective.
\end{proof}

\begin{rem}\label{rem:rgd}
	Let $ \roots $ be any root system. An \emph{RGD-system of type $ \roots $} is (essentially) a $ \roots $-graded group $ (G, (\rootgr{\alpha})_{\alpha \in \roots}) $ such that $ \invset{\alpha} = \rootgr{\alpha} \setminus \Set{1_G} $. These objects were introduced by Tits in~\cite{Tits-TwinBuildingsKacMoody}, and they are (essentially) in bijective correspondence with thick Moufang buildings of type $ \roots $ (see \cite[Section~7.8]{AbramenkoBrown} for details). It follows from~\ref{pro:weyl-inv} that an $ H_\ell $-graded group with a standard coordinatisation cannot be an RGD-system (because $ \ring \times \ring $ contains non-invertible elements, such as $ (x,0) $). Since every $ H_\ell $-grading has a standard coordinatisation by our main result~\ref{thm:prod-param}, it follows that there exist no RGD-systems (and hence no thick Moufang buildings) of types $ H_3 $ and $ H_4 $. This is a classical result of Tits \cite[3.5]{Tits-EndSpiegWeyl}, which he proved by showing the analogous result for type $ H_2 $ \cite[Theorème~1]{Tits-NonEx}.
\end{rem}

\begin{rem}\label{rem:not-chev}
	Let $ (G, (\rootgr{\alpha})_{\alpha \in H_\ell}) $ be an $ H_\ell $-graded group with a standard coordinatisation $ (\risom{\alpha})_{\alpha \in H_\ell} $ by $ \ring \times \ring $.
	It follows from~\ref{lem:stein-epi} that for any positive system $ \possys $ in $ H_\ell $, the subgroup $ \rootgr{\possys} $ of $ G $ is isomorphic to the subgroup $ \rootgrHex{\possys} $ of $ \Chev(X, \ring) $. Thus $ G $ cannot be \enquote{substantially different} from $\Chev(X, \ring)$, but it need not be isomorphic to $ \Chev(X,\ring) $. For example, $ G $ could be the adjoint Chevalley group of type $ D_6 $ as well as the simply-connected Chevalley group of type $ D_6 $, which are not isomorphic. We do not prove that $G$ is (isomorphic to) a Chevalley group, but we are also not aware of an example of an $ H_\ell $-graded group which is provably not a Chevalley group.
\end{rem}

\section{The blueprint technique}
\label{sec:bp}

In \cite{titsronan}, Ronan-Tits introduced the notion of blueprints as a tool for constructing (and hence proving the existence of) buildings of certain types. More precisely, they show that if certain computations involving rewriting rules and a homotopy cycle of the longest element in the Weyl group lead to a certain result, then a blueprint can be used to construct a building. What we refer to as the \enquote{blueprint technique} is the reversal of this idea: Starting from a root graded group $G$ (which should be thought of as a group acting on a \enquote{generalised building}) with a parametrisation (in the sense of~\ref{def:para}), one can perform the same kind of computations to obtain identities that are satisfied by the commutation maps (defined as in~\ref{def:commmap}) in $G$. This application of the ideas in~\cite{titsronan} was first described in \cite{torben}.

In this section, we give a brief summary of the general technique by the example of $ H_\ell $-gradings. A more detailed exposition can be found in \cite[6.1]{torben}. In Section~\ref{sec:bph3}, we will perform the actual blueprint computations for $ H_\ell $-graded groups in order to achieve the goal in~\ref{goal:ring-decomp}.

We begin by introducing some notation and terminology.

\begin{notation}
\label{conv:ideablueprint}
	In this section, we fix $\ell \in \Set{3,4}$ and an $H_\ell$-graded group $(G, (\rootgr{\alpha})_{\alpha \in H_\ell})$. We also choose a root base $\rootbase$ of $H_\ell$ in standard order, which we denote by $(\rho_0, \rho_1, \rho_2, \rho_3)$ if $\ell = 4$ and by $(\rho_1, \rho_2, \rho_3)$ if $\ell = 3$. We will always consider $H_3$ as the subsystem $\subsys{\rho_1, \rho_2, \rho_3}$ of $H_\ell$. Finally, we fix a family $(w_\delta)_{\delta \in \rootbase}$ of Weyl elements and define the nonassociative ring $ \prodring $ and the root isomorphisms $ (\risom{\alpha})_{\alpha \in H_3} $ as in~\ref{conclusion:para}.
\end{notation}

\begin{rem}
	Nearly all computations in this section and the following one take place in the subgroup of $ G $ generated by $ (\rootgr{\alpha})_{\alpha \in H_3} $. Only in~\ref{prop:comm-formula-H4} and~\ref{prop:weyl-1} will we have to leave the $ H_3 $-subsystem to ensure that the commutator relation for $ (\rho_0, \rho_1) $ and the Weyl element $ w_{\rho_0} $ behave as they should.
\end{rem}

\begin{definition}\label{def:red-word}
	A word $f$ over $\rootbase$ is \emph{reduced} if it is not homotopic (in the sense of~\ref{def:homotopy}) to shorter word.
\end{definition}

\begin{notation}
\label{nota:words}
For any word $\bar{\delta}=(\delta_1,\ldots,\delta_k)$ over $\Delta$, we put
\[ U_{\bar{\delta}}\coloneqq U_{\delta_1}\times\cdots\times U_{\delta_k}. \]
\end{notation}

\begin{definition}
\label{def:rr}
	For any word $\word{\delta}$ over $\rootbase$, the \emph{blueprint invariant for $ \word{\delta} $ (with respect to $ (w_\delta)_{\delta \in \rootbase} $)} is the map
	\[ \map{\blumap{\word{\delta}}}{\rootgr{\word{\delta}}}{G}{(g_1, \ldots, g_m)}{\prod_{i=1}^m{w_{\delta_i}g_i}}.\]
	For any two words $ \word{\alpha}, \word{\beta} $ over $ \rootbase $, a \emph{blueprint rewriting rule of type $(\bar{\alpha},\bar{\beta})$ (with respect to $ (w_\delta)_{\delta \in \rootbase} $)} is a map $ \map{\rwrule{}}{\rootgr{\word{\alpha}}}{\rootgr{\word{\beta}}}{}{} $ which satisfies $\gamma_{\bar{\beta}}\circ\rwrule{}=\gamma_{\bar{\alpha}}$.
\end{definition}

\begin{lemma}[{\cite[6.2.11]{torben}}]\label{lem:blue-conc}
	Concatenations of blueprint rewriting rules and inverses of bijective blueprint rewriting rules are blueprint rewriting rules.
\end{lemma}

\begin{lemma}[{\cite[6.2.14]{torben}}]\label{lem:blue-times}
	Let $ \word{\alpha}, \word{\alpha}', \word{\alpha}'', \word{\beta} $ be words over $ \rootbase $ and let $ \rwrule{} $ be a blueprint rewriting rule of type $ (\word{\alpha}, \word{\beta}) $. Then $ \map{\id \times \rwrule{} \times \id}{\rootgr{\word{\alpha}'} \times \rootgr{\word{\alpha}} \times \rootgr{\word{\alpha}''}}{\rootgr{\word{\alpha}'} \times \rootgr{\word{\beta}} \times \rootgr{\word{\alpha}''}}{}{} $ is a blueprint rewriting rule.
\end{lemma}

The heart of the blueprint technique is the following result. It follows from the fact that the blueprint invariant is injective for any reduced word, which is a consequence of~\ref{prop:prodmapbij}.

\begin{theorem}[{\cite[6.2.8]{torben}}]\label{thm:blue}
\label{thm:rr}
	Let $\bar{\alpha},\bar{\beta}$ be two reduced words over $\Delta$ such that $\sigma_{\bar{\alpha}}=\sigma_{\bar{\beta}}$, let $x\in U_{\bar{\alpha}}$ and let $ \map{\rwrule{1}, \rwrule{2}}{\rootgr{\word{\alpha}}}{\rootgr{\word{\beta}}}{}{} $ be two blueprint rewriting rules. Then $\rwrule{1}(x)=\rwrule{2}(x)$.
\end{theorem}

\begin{algorithm}
\label{algo:blueprint}
Put $\theta_i\coloneqq\theta_{\rho_i}$ and $ \rootgr{i} \defl \rootgr{\rho_i} $ for all $ i \in \{1,2,3\} $ and $ \rootgr{i_1 \cdots i_k} \defl \rootgr{i_1} \times \cdots \times \rootgr{i_k} $ for all $ i_1, \ldots, i_k \in \{1,2,3\} $. Further, we denote by $ \longel $ the longest element in the Weyl group of $ H_3 $ and choose an arbitrary reduced expression of $ \longel $, say
	\[f_1=(\rho_3,\rho_2,\rho_3,\rho_2,\rho_3,\rho_1,\rho_2,\rho_3,\rho_2,\rho_3,\rho_1,\rho_2,\rho_3,\rho_2,\rho_1).\]
That is, $ \reflbr{f_1} = \longel $. To simplify notation, we will write $f_1=323231232312321$. Further, we define $ f_2, \ldots, f_{63} $ to be the \emph{homotopy cycle} for $ \longel $ given by Figure~\ref{fig:homcycle}. Observe that $ f_i $ is elementary braid-homotopic to $ f_{i+1} $ for all $ i \in \{1, \ldots, 62\} $ and that $ f_1 = f_{63} $.

Recall that the Weyl group of $ H_3 $ satisfies the braid relations
\[ \reflbr{\rho_1 \rho_2 \rho_1} = \reflbr{\rho_2 \rho_1 \rho_2}, \; \reflbr{\rho_2 \rho_3 \rho_2 \rho_3 \rho_2} = \reflbr{\rho_3 \rho_2 \rho_3 \rho_2 \rho_3}, \; \reflbr{\rho_1 \rho_3} = \reflbr{\rho_3 \rho_1}. \]
We assume that for all distinct $ i,j \in \{1,2,3\} $, we have a blueprint rewriting rule $ \rwrule{ij} $ of type $ \brackets{\braidword{m_{ij}}(\rho_i, \rho_j), \braidword{m_{ji}}(\rho_j, \rho_i)} $. For example, we have blueprint rewriting rules $ \map{\rwrule{12}}{\rootgr{121}}{\rootgr{212}}{}{} $ and $ \map{\rwrule{21}}{\rootgr{212}}{\rootgr{121}}{}{} $. Using these \enquote{elementary} blueprint rewriting rules and~\ref{lem:blue-times}, we can define blueprint rewriting rules $ \rwrule{i} $ of type $ (f_i, f_{i+1}) $ for all $ i \in \{1, \ldots, 62\} $. Similarly, we can define blueprint rewriting rules $ \rwrule{i}' $ of type $ (f_{i+1}, f_i) $ for all $ i \in \{2, \ldots, 63\} $.

Now let $y_1,y_2,\ldots,y_{15}\in\prodring$ be arbitrary ring elements and consider the tuple 
	\[x_1 \coloneqq \big(\theta_3(y_1),\theta_2(y_2),\theta_3(y_3),\ldots,\theta_1(y_{15})\big) \in \rootgr{f_1} = \rootgr{f_{63}}.\]
For all $ i \in \{2, \ldots, 63\} $, we can successively compute tuples $ x_i \defl \rwrule{i-1}(x_{i-1}) $. Since the identity map on $ \rootgr{f_1} $ and $ \rwrule{62} \circ \cdots \circ \rwrule{1} $ are blueprint rewriting rules, it follows from~\ref{thm:blue} that $ x_1 = x_{63} $. Writing
\[ x_{63} = \big(\theta_3(y_1^{(63)}),\theta_2(y_2^{(63)}),\theta_3(y_3^{(63)}),\ldots,\theta_1(y_{15}^{(63)})\big), \]
we obtain 15 identities of the form $ y_i = y_i^{(63)} $ which hold for arbitrary ring elements $ y_1, \ldots, y_{15} $.

In practice, we will only compute $ x_2, \ldots, x_{32} $ and, in a separate computation, tuples $ x_{63}', \ldots, x_{32}' $ defined by $ x_{63}' \defl x_1 $ and $ x_i' \defl \rwrule{i+1}'(x_{i+1}') $ for all $ i \in \{32, \ldots, 62\} $. By the same argument as above, we then have $ x_{32} = x_{32}' $, which yields 15 identities in $ \prodring $ (which are shorter than the 15 identities in the previous paragraph). We will refer to these identities as the \emph{blueprint identities}.
\end{algorithm}

\begin{figure}[htb]
	\begin{multicols}{3}
	\begin{center}
		\begin{enumerate}[(1)]
			\item \underline{32323}1232312321
			\item 2323\underline{212}32312321
			\item 232312\underline{13}2312321
			\item 232312312\underline{31}2321
			\item 2323123\underline{121}32321
			\item 232312321\underline{23232}1
			\item 232\underline{31}2321323231
			\item 23213232\underline{13}23231
			\item 2321\underline{32323}123231
			\item 23\underline{212}3232123231
			\item 23121323\underline{212}3231
			\item 2312\underline{13}231213231
			\item 2312312312\underline{13}231
			\item 2312312\underline{31}231231
			\item 2\underline{31}231213231231
			\item 21323\underline{121}3231231
			\item 2132321232312\underline{31}
			\item 21323212323\underline{121}3
			\item 2132321\underline{23232}123
			\item 213232\underline{13}2323123
			\item 21\underline{32323}12323123
			\item \underline{212}323212323123
			\item 121323\underline{212}323123
			\item 12\underline{13}23121323123
			\item 12312\underline{31}21323123
			\item 123\underline{121}321323123
			\item 12321232\underline{13}23123
			\item 12321232312\underline{31}23
			\item 123212323\underline{121}323
			\item 12321\underline{23232}12323
			\item 1232\underline{13}232312323
			\item 123231232\underline{31}2323
			\item 1232312321\underline{32323}
			\item 12323123\underline{212}3232
			\item 1232312\underline{31}213232
			\item 12323\underline{121}3213232
			\item 1\underline{23232}123213232
			\item \underline{13}2323123213232
			\item 31232\underline{31}23213232
			\item 3123213232\underline{13}232
			\item 312321\underline{32323}1232
			\item 3123\underline{212}32321232
			\item 312\underline{31}2132321232
			\item 3\underline{121}32132321232
			\item 321232\underline{13}2321232
			\item 3212323123\underline{212}32
			\item 321232312\underline{31}2132
			\item 3212323\underline{121}32132
			\item 321\underline{23232}1232132
			\item 32\underline{13}23231232132
			\item 3231232\underline{31}232132
			\item 323123213232\underline{13}2
			\item 32312321\underline{32323}12
			\item 323123\underline{212}323212
			\item 32312\underline{31}21323212
			\item 323\underline{121}321323212
			\item 32321232\underline{13}23212
			\item 323212323123\underline{212}
			\item 32321232312\underline{31}21
			\item 323212323\underline{121}321
			\item 32321\underline{23232}12321
			\item 3232\underline{13}232312321
			\item 323231232312321
		\end{enumerate}
		\end{center}
	\end{multicols}
\caption{A homotopy cycle for the longest word in $H_3$. The underlines indicate the elementary braid homotopy moves which are applied at each step.}
\label{fig:homcycle}
\end{figure}

The computation in~\ref{algo:blueprint} will be referred to as the \emph{blueprint computation}. In order to actually perform this computation, it remains to compute the (elementary) blueprint rewriting rules.

\begin{notation}
\label{conv:bph3i25quin}
	From now on, we put $\alpha\coloneqq\rho_2$, $\ep\coloneqq\rho_3$ and we denote the $H_2$-quintuple corresponding to the $ H_2 $-pair $ (\alpha, \ep) $ by $(\alpha,\beta,\gamma,\delta,\ep)$. We will also make frequent use of Notation~\ref{def:commmap} for commutation maps (with respect to $ (\risom{\zeta})_{\zeta \in H_\ell} $).
\end{notation}

\begin{definition}
\label{def:rrh3}
We define the following maps:
		\begin{align*}
			\rwrule{12}:~ &U_{\rho_1}\times U_{\rho_2}\times U_{\rho_1}\rightarrow U_{\rho_2}\times U_{\rho_1}\times U_{\rho_2},\\*
			&\left(\theta_{\rho_1}(a),\theta_{\rho_2}(b),\theta_{\rho_1}(c)\right)\mapsto\left(\theta_{\rho_2}(c),\theta_{\rho_1}(-b-ca),\theta_{\rho_2}(a)\right),\\
			\rwrule{21}:~ &U_{\rho_2}\times U_{\rho_1}\times U_{\rho_2}\rightarrow U_{\rho_1}\times U_{\rho_2}\times U_{\rho_1},\\*
		&\left(\theta_{\rho_2}(a),\theta_{\rho_1}(b),\theta_{\rho_2}(c)\right)\mapsto\left(\theta_{\rho_1}(c),\theta_{\rho_2}(-b-ac),\theta_{\rho_1}(a)\right),\\
			\rwrule{13}:~ &U_{\rho_1}\times U_{\rho_3}\rightarrow U_{\rho_3}\times U_{\rho_1},\\*
			&\left(\theta_{\rho_1}(a),\theta_{\rho_3}(b)\right)\mapsto\left(\theta_{\rho_3}(b),\theta_{\rho_1}(a)\right),\\
		\rwrule{31}:~ &U_{\rho_3}\times U_{\rho_1}\rightarrow U_{\rho_1}\times U_{\rho_3},\\*
			&\left(\theta_{\rho_3}(a),\theta_{\rho_1}(b)\right)\mapsto\left(\theta_{\rho_1}(b),\theta_{\rho_3}(a)\right).
	\end{align*}
Furthermore, we define
	\begin{align*}
		\rwrule{23}:~ &U_{\rho_2}\times U_{\rho_3}\times U_{\rho_2}\times U_{\rho_3}\times U_{\rho_2}\rightarrow U_{\rho_3}\times U_{\rho_2}\times U_{\rho_3}\times U_{\rho_2}\times U_{\rho_3},\\
		&\left(\theta_{\rho_2}(a),\theta_{\rho_3}(b),\theta_{\rho_2}(c),\theta_{\rho_3}(d),\theta_{\rho_2}(e)\right)\\
		&\quad \mapsto \left(\theta_{\rho_3}(e),\theta_{\rho_2}(B),\theta_{\rho_3}(C),\theta_{\rho_2}(D),\theta_{\rho_3}(a)\right)
	\end{align*}
where
	\begin{align*}
		B &\coloneqq -\Big(\psi_{\delta,\alpha}^\beta(-b,e)+\psi_{\ep,\alpha}^\beta(a,e)+\psi_{\gamma,\alpha}(c,e)-d\Big)^*,\\
		C & \coloneqq -\Big(\psi_{\delta,\alpha}^\gamma(-b,e)+\psi_{\delta,\beta}\brackets[\big]{-b, \psi_{\delta, \alpha}^{\beta}(-b,e)} \\*
		&\qquad\quad \mathord{} +\psi_{\ep,\alpha}^\gamma(a,e)+\psi_{\ep,\beta}^\gamma\brackets[\big]{a,\psi_{\ep,\alpha}^\beta(a,e)+\psi_{\gamma,\alpha}(c,e)-d} \\*
		&\qquad\quad \mathord{}+\psi_{\delta,\beta}\brackets[\big]{-b+\psi_{\ep,\alpha}^\delta(a,e),\psi_{\ep,\alpha}^\beta(a,e) + \psi_{\gamma,\alpha}(c,e)-d}+c\Big),\\		
		D &\coloneqq -\bigg(-b+\psi_{\ep,\alpha}^\delta(a,e)+\psi_{\ep,\beta}^\delta\brackets[\big]{a,\psi_{\ep,\alpha}^\beta(a,e)+\psi_{\gamma,\alpha}(c,e)-d} \\*
		&\qquad\quad \mathord{}+\psi_{\ep,\gamma}\brackets[\Big]{a,\psi_{\ep,\beta}^\gamma\brackets[\big]{a,\psi_{\ep,\alpha}^\beta(a,e)+\psi_{\gamma,\alpha}(c,e)-d}+\psi_{\ep,\alpha}^\gamma(a,e)+c}\bigg)^*.		
	\end{align*}
Similarly, we define
	\begin{align*}
		\rwrule{32}:~ &U_{\rho_3}\times U_{\rho_2}\times U_{\rho_3}\times U_{\rho_2}\times U_{\rho_3}\rightarrow U_{\rho_2}\times U_{\rho_3}\times U_{\rho_2}\times U_{\rho_3}\times U_{\rho_2},\\
		&\left(\theta_{\rho_3}(a),\theta_{\rho_2}(b),\theta_{\rho_3}(c),\theta_{\rho_2}(d),\theta_{\rho_3}(e)\right)\\
		&\quad \mapsto \left(\theta_{\rho_2}(e),\theta_{\rho_3}(B),\theta_{\rho_2}(C),\theta_{\rho_3}(D),\theta_{\rho_2}(a)\right)
	\end{align*}
where
	\begin{align*}		
		B &\coloneqq -\Big(\psi_{\beta,\ep}^\delta(-b^*,e)+\psi_{\alpha,\ep}^\delta(a,e)+\psi_{\gamma,\ep}(-c,e)-d^*\Big),\\		
		C &\coloneqq \psi_{\beta,\ep}^\gamma(-b^*,e)+\psi_{\beta,\delta}\brackets[\big]{-b^*,\psi_{\beta,\ep}^\delta(-b^*,e)}\\*
		&\qquad\quad \mathord{}+\psi_{\beta,\delta}\brackets[\Big]{-b^*+\psi_{\alpha,\ep}^\beta(a,e),\psi_{\alpha,\ep}^\delta(a,e)+\psi_{\gamma,\ep}(-c,e)-d^*}\\*
		&\qquad\quad \mathord{} +\psi_{\alpha,\delta}^\gamma\brackets[\big]{a,\psi_{\alpha,\ep}^\delta(a,e)+\psi_{\gamma,\ep}(-c,e)-d^*}+\psi_{\alpha,\ep}^\gamma(a,e)-c,\\
		D&\coloneqq -\bigg(-b^*+\psi_{\alpha,\ep}^\beta(a,e)+\psi_{\alpha,\delta}^\beta\brackets[\big]{a,\psi_{\alpha,\ep}^\delta(a,e)+\psi_{\gamma,\ep}(-c,e)-d^*}\\*
		& \qquad \quad \mathord{}+\psi_{\alpha,\gamma}\brackets[\Big]{a,\psi_{\alpha,\delta}^\gamma\brackets[\big]{a,\psi_{\alpha,\ep}^\delta(a,e)+\psi_{\gamma,\ep}(-c,e)-d^*}+\psi_{\alpha,\ep}^\gamma(a,e)-c}\bigg).
	\end{align*}
\end{definition}

The following computation shows how to arrive at the formulas in~\ref{def:rrh3}.

\begin{lemma}[{\cite[6.3.4]{torben}}]
\label{lem:a1a2gammacomp}
The maps $\rwrule{12}$, $\rwrule{21}$, $\rwrule{13}$ and $\rwrule{31}$ are blueprint rewriting rules.
\end{lemma}
\begin{proof}
Let $a,b,c,A,B,C\in\mathcal{R}$ such that
	\[\rwrule{12}\brackets[\big]{\theta_{\rho_1}(a),\theta_{\rho_2}(b),\theta_{\rho_1}(c)} = \brackets[\big]{\theta_{\rho_2}(A),\theta_{\rho_1}(B),\theta_{\rho_2}(C)}\]
and put $f\coloneqq(\rho_1,\rho_2,\rho_1)$. Then, by \ref{def:rr}, Figure~\ref{fig:parmapex-H3} and~\eqref{eq:H3-param}, we have
	\begin{align*}
		\gamma_f(\theta_{\rho_1}(a),\theta_{\rho_2}(b),\theta_{\rho_1}(c)) &= w_{\rho_1}w_{\rho_2}w_{\rho_1}\theta_{\rho_1}(a)^{w_{\rho_2}w_{\rho_1}}\theta_{\rho_2}(b)^{w_{\rho_1}}\theta_{\rho_1}(c)\\
		&= w_{\rho_1}w_{\rho_2}w_{\rho_1}\theta_{\rho_1+\rho_2}(a)^{w_{\rho_1}}\theta_{\rho_1+\rho_2}(-b)\theta_{\rho_1}(c) \\
		&= w_{\rho_1}w_{\rho_2}w_{\rho_1}\theta_{\rho_2}(a)\theta_{\rho_1+\rho_2}(-b)\theta_{\rho_1}(c).
	\end{align*}
Similarly, with $g\coloneqq(\rho_2,\rho_1,\rho_2)$, we compute
	\begin{align*}
		\gamma_g(\theta_{\rho_2}(A),\theta_{\rho_1}(B),\theta_{\rho_2}(C)) &= w_{\rho_2}w_{\rho_1}w_{\rho_2}\theta_{\rho_2}(A)^{w_{\rho_1}w_{\rho_2}}\theta_{\rho_1}(B)^{w_{\rho_2}}\theta_{\rho_2}(C)\\
		&= w_{\rho_2}w_{\rho_1}w_{\rho_2}\theta_{\rho_1+\rho_2}(-A)^{w_{\rho_2}}\theta_{\rho_1+\rho_2}(B)\theta_{\rho_2}(C) \\
		&= w_{\rho_2}w_{\rho_1}w_{\rho_2}\theta_{\rho_1}(A)\theta_{\rho_1+\rho_2}(B)\theta_{\rho_2}(C).
	\end{align*}
Here $w_{\rho_1}w_{\rho_2}w_{\rho_1}=w_{\rho_2}w_{\rho_1}w_{\rho_2}$ by \ref{braidrel} and
	\begin{align*}
		\theta_{\rho_2}(a)\theta_{\rho_1+\rho_2}(-b)\theta_{\rho_1}(c) &= \theta_{\rho_1+\rho_2}(-b)\theta_{\rho_1}(c)\theta_{\rho_2}(a)\left[\theta_{\rho_1}(c),\theta_{\rho_2}(a)\right]\inv\\
		&= \theta_{\rho_1+\rho_2}(-b)\theta_{\rho_1}(c)\theta_{\rho_2}(a)\theta_{\rho_1+\rho_2}(-ca)\\
		&= \theta_{\rho_1}(c)\theta_{\rho_1+\rho_2}(-b-ca)\theta_{\rho_2}(a).
	\end{align*}
Since $A=c$, $B=-b-ca$ and $C=a$ by~\ref{def:rrh3}, we infer that $\rwrule{12}$ is a blueprint rewriting rule. By a similar calculation, the same holds for $\rwrule{13}$. Lastly, note that $\rwrule{21}$ and $\rwrule{31}$ are the inverses of $\rwrule{12}$ and $\rwrule{13}$, respectively. Hence they are blueprint rewriting rules by~\ref{lem:blue-conc}.
\end{proof}

\begin{lemma}
\label{lem:i25gammacomp}
The maps $\rwrule{23}$ and $\rwrule{32}$ are blueprint rewriting rules.
\end{lemma}
\begin{proof}
The proof proceeds just as the one of~\ref{lem:a1a2gammacomp}. We omit the details.
\end{proof}

\section[The blueprint computation for \texorpdfstring{$ H_3 $}{H3}-graded groups]{The blueprint computation for \texorpdfstring{$ H_3 $}{H3}-graded\\ groups}
\label{sec:bph3}

\begin{notation}\label{notation:blueh3}
	The notation introduced in~\ref{conv:ideablueprint} and~\ref{conv:bph3i25quin} continues to hold.
\end{notation}

Using the blueprint rewriting rules from~\ref{def:rrh3} and the homotopy cycle from Figure~\ref{fig:homcycle}, we can perform the blueprint computation for $ (G, (\rootgr{\alpha})_{\alpha \in H_3}) $, following the algorithm in~\ref{algo:blueprint}. This yields 15 blueprint identities involving indeterminates $ y_1, \ldots, y_{15} \in \prodring $. Using these identities, we will show in this section that $ \prodring $ is an associative commutative ring (see~\ref{prop:commutative}, \ref{thm:associative}) which decomposes as a direct sum $ \prodring = \prodring_1 \oplus \prodring_2 $ of two isomorphic ideals $ \prodring_1 $, $ \prodring_2 $ (see~\ref{prop:dirsum}, \ref{prop:R1R2iso}), and we will also compute explicit formulas for the commutation maps $ \commmap{\xi,\zeta}{\rho} $ from~\ref{conclusion:para}.

\begin{rem}
\label{rem:bpevagen}
	Unfortunately, the majority of the 15 blueprint identities are rather long and complicated: They fill a total of 20 A4 pages. Consequently, they are not accessible for manipulation by hand, so we have decided not to include them in this paper. However, since these identities hold for all possible choices of variables $ y_1, \ldots, y_{15} \in \prodring $, we can infer more succinct \emph{evaluated blueprint identities} by replacing most of the indeterminates with $ 0_\prodring $ (or, in some cases, $ 1_\prodring $). The evaluated blueprint identities can be simplified by hand, which yields identities labelled (1), \dots, (35). We will not present how to deduce (1), \dots, (35) from the evaluated blueprint identities because this is a straightforward procedure, but we will give some remarks in the proofs of~\ref{prop:blue-eval1}, \ref{prop:blue-eval2} and~\ref{prop:blue-eval3}. The interested reader can find the original blueprint identities, the evaluated blueprint identities and the code that produced them in~\cite{repo}.
\end{rem}

Before we begin with the derivation of the identities (1), \dots, (35), we make some preliminary observations.

\begin{prop}
\label{prop:commutative}
The multiplication on $\prodring$ is commutative.
\end{prop}
\begin{proof}
This follows from the second blueprint identity, which reads $-y_{14}-y_5y_{15}=-y_{14}-y_{15}y_5$.
\end{proof}

\begin{rem}
\label{rem:commmapinv}
	Let $\zeta, \xi \in H_3$ be non-proportional and let $\rho \in \oprootint{\zeta, \xi}$. Then
	\begin{equation*}
		\commmap{\zeta, \xi}{\rho}(x,y) = -\commmap{\xi, \zeta}{\rho}(y,x)
	\end{equation*}
	for all $x, y \in \prodring$. This holds by \ref{lem:commpart}.
\end{rem}

\begin{rem}\label{lem:psiadditive-one}
	For all $\zeta,\xi\in H_3$ such that $\oprootint{\zeta,\xi}$ contains exactly one element, the map $ \psi_{\zeta,\xi} $ is bi-additive by~\ref{lem:oneadditivity}.
\end{rem}

\begin{rem}
\label{lem:psisadditive-weak}
	Let $ (\bar{\alpha}, \bar{\beta}, \bar{\gamma}, \bar{\delta}, \bar{\epsilon}) $ be an $ H_2 $-quintuple and let $x,x',y,y'\in\prodring$. Then the following hold by~\ref{lem:i25additivity} and~\ref{lem:int2add}:
	\begin{enumerate}[(i)]
		\item $ \psi_{\bar{\alpha}, \bar{\delta}}^{\bar{\beta}}(x,y+y') = \psi_{\bar{\alpha}, \bar{\delta}}^{\bar{\beta}}(x,y)+\psi_{\bar{\alpha}, \bar{\delta}}^{\bar{\beta}}(x,y') $.
		
		\item $ \psi_{\bar{\alpha}, \bar{\delta}}^{\bar{\beta}}(x+x',y) = \psi_{\bar{\alpha}, \bar{\delta}}^{\bar{\beta}}(x,y) + \psi_{\bar{\alpha}, \bar{\delta}}^{\bar{\beta}}(x',y) + \psi_{\bar{\gamma}, \bar{\alpha}}\brackets[\big]{\psi_{\bar{\alpha}, \bar{\delta}}^{\bar{\gamma}}(x,y), y'} $.
		
		\item $ \psi_{\bar{\alpha}, \bar{\delta}}^{\bar{\gamma}}(x,y+y') = \psi_{\bar{\alpha}, \bar{\delta}}^{\bar{\gamma}}(x,y) + \psi_{\bar{\alpha}, \bar{\delta}}^{\bar{\gamma}}(x,y') + \psi_{\bar{\beta}, \bar{\delta}}\brackets[\big]{\psi_{\bar{\alpha}, \bar{\delta}}^{\bar{\beta}}(x,y), y'} $.
		
		\item $ \psi_{\bar{\alpha}, \bar{\delta}}^{\bar{\gamma}}(x+x',y) = \psi_{\bar{\alpha}, \bar{\delta}}^{\bar{\gamma}}(x,y) + \psi_{\bar{\alpha}, \bar{\delta}}^{\bar{\gamma}}(x',y) $.
		
		\item $ \psi_{\bar{\alpha}, \bar{\epsilon}}^{\bar{\beta}}(x,y+y') = \psi_{\bar{\alpha}, \bar{\epsilon}}^{\bar{\beta}}(x,y) + \psi_{\bar{\alpha}, \bar{\epsilon}}^{\bar{\beta}}(x,y') $.
	\end{enumerate}
\end{rem}

\begin{prop}\label{prop:blue-eval1}
	All identities in the column \enquote{Result} of Figure~\ref{fig:eva1} hold for all $ y_1, \ldots, y_{15} \in \prodring $.
\end{prop}
\begin{proof}
	This is a straightforward computation, as explained in~\ref{rem:bpevagen}. For example, in order to derive (4), we put $ y_{12} \defl 1_\prodring $ and $ y_i \defl 0_\prodring $ for all $ i \in \{1, \ldots, 15\} \setminus \{4,5,12\} $ in blueprint identity~9. This yields the evaluated blueprint identity
	\[ \invo{\commmap{\gamma,\epsilon}{}\brackets[\big]{-\invo{y_4}(-\invo{1_\prodring}),y_5}} = -\invo{\commmap{\gamma,\epsilon}{}(-y_5,y_4)}. \]
	Using (2), (3) and~\ref{lem:psiadditive-one}, this identity is seen to be equivalent to (4). The identities have to be derived in the given order, so that each identity may be used in the proof of all subsequent identities.
\end{proof}

\premidfigure
\begin{figure}[p]
\centering
{\renewcommand{\arraystretch}{1.7}
\begin{tabular}{cccc}
	\toprule
	No. & \makecell{Blueprint\\identity} & \makecell{Non-zero\\variables} & Result \\
	\midrule
	(1) & $12$ & $y_{4}$ & $(-y_4)^*=-y_4^*$ \\
	(2) & $5$ & $y_{2},y_{15}$ & $y_{15}y_{2}^*=(y_{15}y_{2})^*$ \\
	(3) & $9$ & $y_{5},y_{8}$ & $\psi_{\ep,\gamma}(y_5,y_8)=-\psi_{\alpha,\gamma}(y_5,y_8)$ \\
	(4) & $9$ & $y_{4},y_{5}$, $y_{12}=1$ & $\psi_{\alpha,\gamma}(y_4,y_5)=\psi_{\alpha,\gamma}(y_5,y_4)$ \\
	(5) & $14$ & $y_{1},y_{3}$ & $\psi_{\alpha,\gamma}(y_1,y_3)^*=\psi_{\alpha,\gamma}(y_1,y_3)$ \\
	(6) & $5$ & $y_{1},y_{12}$ & $\psi_{\alpha,\gamma}(y_1,y_{12}^*)=-\psi_{\alpha,\gamma}(y_{1},y_{12})$ \\
	(7) & $7$ & $y_{4}=1,y_{10},y_{12}$ & $\psi_{\beta,\delta}(y_{10},y_{12})=\psi_{\ep,\gamma}(y_{10},y_{12})$ \\
	(8) & $7$ & $y_{4},y_{10},y_{12}$ & $ \psi_{\alpha,\gamma}(y_4y_{10},y_{12}) =\psi_{\alpha,\gamma}(y_4y_{12},y_{10}) $ \\
	(9) & $13$ & $y_{1},y_{4}$ & $\psi_{\ep,\beta}^\gamma(y_1,y_4)=-\psi_{\alpha,\delta}^\gamma(y_1,y_4^*)$ \\
	(10) & $3$ & $y_{7},y_{15}$ & $\psi_{\beta,\ep}^\delta(y_7,y_{15})=-\psi_{\delta,\alpha}^\beta(y_7,y_{15})$ \\
	(11) & $3$ & $y_{10},y_{14}$ & $\psi_{\beta,\ep}^\delta(y_{10}^*,y_{14})=\psi_{\beta,\ep}^\delta(y_{10},-y_{14})$ \\
	(12) & $12$ & $y_{2},y_{5}$ & $\psi_{\beta,\ep}^\delta(y_2,y_5)^*=-\psi_{\beta,\ep}^\delta(y_2,-y_5)$ \\
	(13) & $7$ & $y_{7},y_{10}$ & $\psi_{\beta,\ep}^\delta(y_7,y_{10}^*)=\psi_{\beta,\ep}^\delta(y_7,y_{10})$ \\
	(14) & $11$ & $y_{1},y_{3},y_{8}$ & $\psi_{\alpha,\gamma}(\psi_{\alpha,\gamma}(y_1,y_3),y_8)=0$ \\
	(15) & $11$ & $y_{1},y_{4},y_{9}$ & $\psi_{\alpha,\gamma}(\psi_{\alpha,\delta}^\gamma(y_1,y_4),y_9)=0$ \\
	(16) & $6$ & $y_{3},y_{5},y_{13}$ & $\psi_{\alpha,\delta}^\gamma(\psi_{\alpha,\gamma}(y_3,y_5),y_{13})=0$ \\
	(17) & $7$ & $y_{1},y_{4},y_{12}$ & $\psi_{\alpha,\delta}^\gamma(\psi_{\alpha,\delta}^\gamma(y_1,y_4),y_{12})=0$ \\
	(18) & $8$ & $y_{4},y_{6},y_{15}$ & $ \psi_{\alpha,\gamma}(y_6,\psi_{\alpha,\ep}^\gamma(y_{15},y_4)) =y_4\psi_{\beta,\ep}^\delta(y_6,y_{15}) $ \\
	(19) & $6$ & $y_{7},y_{10},y_{14}$ & $\psi_{\alpha,\gamma}(\psi_{\beta,\ep}^\delta(y_{10},y_{14}),y_7)=0$ \\
	(20) & $9$ & $y_{2},y_{5},y_{14}=1$ & $\psi_{\beta,\ep}^\delta(y_5,y_2)=-\psi_{\ep,\beta}^\gamma(y_5,-y_2)$ \\
	(21) & $3$ & $y_{1},y_{14}$ & $\psi_{\ep,\alpha}^\beta(y_1,y_{14})=\psi_{\alpha,\ep}^\delta(y_1,-y_{14})$ \\
	\bottomrule
\end{tabular}
}
\caption{Evaluation of the blueprint identities, part~1. See~\ref{prop:blue-eval1}.}
\label{fig:eva1}
\end{figure}
\postmidfigure

\begin{rem}\label{rem:invo-one}
	Let $ x,y \in \prodring $. Since $ \prodring $ is commutative, it follows from~(2) that
	\begin{align*}
		\invo{y} &= \invo{(y1)} = y \invo{1} = \invo{1}y \quad \text{and} \quad xy = \invo{(\invo{(xy)})} = \invo{(x\invo{y})} = \invo{x} \invo{y}.
	\end{align*}
\end{rem}

\begin{notation}
	From now on, we put $ \comone \defl \commmap{\alpha,\gamma}{} $, $ \comtwo \defl \commmap{\alpha,\delta}{\beta} $, $ \comthrone \defl \commmap{\alpha, \ep}{\beta} $ and $ \comthrtwo \defl \commmap{\alpha, \ep}{\gamma} $.
\end{notation}

The following lemma summarises the identities in Figure~\ref{fig:eva1}. In particular, we see that all commutation maps $ \commmap{\xi,\zeta}{\rho} $ for distinct $ \xi,\zeta \in \{\alpha, \beta, \gamma, \delta, \ep\} $ and $ \rho \in \oprootint{\xi, \zeta} $ can be expressed in terms of the four maps $ \comone $, $ \comtwo $, $ \comthrone $ and $ \comthrtwo $.

\begin{lemma}\label{lem:psi-id}
	The following hold (for all $ x,y,z \in \prodring $), where $ \switch $ denotes the map $ \map{}{}{}{(x,y)}{(y,x)} $:
	\begin{enumerate}[(i)]
		\item \label{lem:psi-id:invo}$ \invo{(xy)} = x\invo{y} $.
		\item \label{lem:psi-id:leave-out-1}$ \comone = \psi_{\alpha, \gamma} = \psi_{\gamma, \epsilon} = -\psi_{\beta, \delta} = -\psi_{\gamma, \alpha} = -\psi_{\epsilon, \gamma} = \psi_{\delta, \beta} $.
		
		\item \label{lem:psi-id:invo-1}$ \comone(x,y) = \comone(y,x) = \invo{\comone(x,y)} $ and $ \comone(\invo{x},y) = \comone(x,\invo{y}) = -\comone(x,y) $.

		\item \label{lem:psi-id:assoc}$ \comone(x,y) $ depends only on the image of $ xy $ in the largest associative quotient of $ \prodring $.

		\item \label{lem:psi-id:leave-out-2}
		$ \comtwo = \psi_{\alpha, \delta}^\beta = -\psi_{\ep, \beta}^\delta = -\psi_{\beta, \ep}^\gamma = \psi_{\delta,\alpha}^\gamma = -\psi_{\delta, \alpha}^\beta \circ \switch = \psi_{\beta, \ep}^\delta \circ \switch = \psi_{\ep, \beta}^\gamma \circ \switch = -\psi_{\alpha, \delta}^\gamma \circ \switch $.
		
		\item \label{lem:psi-id:biadd}$ \comone $ and $ \comtwo $ are biadditive and $ \comthrone $ is additive in the second component.
		
		\item \label{lem:psi-id:invo-2}
		$ \comtwo(\invo{x},y) = \comtwo(x,y) = \invo{\comtwo(x,y)} $ and $ \comtwo(x, \invo{y}) = -\comtwo(x,y) $.

		\item \label{lem:psi-id:leave-out-31}$ \comthrone(x,y) = \commmap{\alpha,\ep}{\beta}(x,y) = -\commmap{\alpha,\ep}{\delta}(y,-x) = -\commmap{\ep, \alpha}{\beta}(y,x) = \commmap{\ep,\alpha}{\delta}(-x,y) $.
		
		\item \label{lem:psi-id:one-three}$ \comone\brackets[\big]{x, \comthrtwo(y,z)} = \comone\brackets[\big]{\comthrtwo(y,z), x} = z \comtwo(y,x) $.
		
		\item \label{lem:psi-id:zero}$ \comone\brackets[\big]{\comone(x,y),z} = \comone\brackets[\big]{\comtwo(x,y),z} = \comtwo\brackets[\big]{x, \comone(y,z)} = \comtwo\brackets[\big]{x, \comtwo(y,z)} = 0 $.
	\end{enumerate}
\end{lemma}
\begin{proof}
	Assertion~\ref{lem:psi-id:invo} is~(2). Assertion~\ref{lem:psi-id:invo-1} follows from (4), (5) and~(6) and~\ref{lem:psi-id:assoc} follows from~(4) and~(8). Assertion~\ref{lem:psi-id:leave-out-1} follows from~(3) and~(7) for the first three equalities and from~(4) and~\ref{rem:commmapinv} for the remaining equalities. By~(10), (20) and~(9), we have
	\begin{equation}\label{eq:psi-id}
		\comtwo(x,y) = \psi_{\alpha, \delta}^\beta(x,y) = -\psi_{\epsilon, \beta}^\delta(x,y) = \psi_{\beta, \epsilon}^\gamma(-x,y) = -\psi_{\delta, \alpha}^\gamma(-\invo{x}, y).
	\end{equation}
	In particular, the images of the maps $ \psi_{\alpha, \delta}^\beta $, $ \psi_{\epsilon, \beta}^\delta $, $ \psi_{\beta, \epsilon}^\gamma $, $ \psi_{\delta, \alpha}^\gamma $ are identical, which allows us to apply (15) and (19) in a more general context. Assertion~\ref{lem:psi-id:biadd} follows from~\ref{lem:psiadditive-one}, \ref{lem:psisadditive-weak}, (15), (19) and the previous observation. Assertion~\ref{lem:psi-id:invo-2} follows from (11), (13), ~\eqref{eq:psi-id} and~\ref{lem:psi-id:biadd}. Assertion~\ref{lem:psi-id:leave-out-2} follows from~\eqref{eq:psi-id}, \ref{lem:psi-id:biadd} and~\ref{lem:psi-id:invo-2}. Finally, assertion~\ref{lem:psi-id:leave-out-31} is (21), assertion~\ref{lem:psi-id:one-three} is~(18) and assertion~\ref{lem:psi-id:zero} is the culmination of~(14) to~(17).
\end{proof}

\begin{prop}\label{prop:blue-eval2}
	All identities in the column \enquote{Result} of Figure~\ref{fig:eva2} hold for all $ y_1, \ldots, y_{15} \in \prodring $.
\end{prop}
\begin{proof}
	We proceed as in~\ref{prop:blue-eval1}. For most identities, the first step is to apply~\ref{lem:psi-id}~\ref{lem:psi-id:zero} in order to remove zero summands. For (25), we use (24). For (30), we use (22) and (23). For (28), we use (24), (23) and (21). For (31), we use (21) and the substitutions $ y_5 \mapsto -y_5 $ and $ y_{15} \mapsto -y_{15} $. For (32), we use that the images of $ \comone $ and $ \comtwo $ are fixed by $ \invosym $.
\end{proof}

\premidfigure
\begin{figure}[htb]
	\centering{\renewcommand{\arraystretch}{1.7}
	\begin{tabular}{cccc}
		\toprule
		No. & \makecell{Blueprint\\identity} & \makecell{Non-zero\\variables} & Result \\
		\midrule
		(22) & $3$ & $y_{4},y_{15}$ & $ \comthrone(-y_{15}, y_4) = \comthrone(y_{15}, \invo{y_4}) $ \\
		(23) & $12$ & $y_{1},y_{5}$ & $ \comthrone(-y_{5}, y_1) = \invo{\comthrone(y_5, y_1)} $ \\
		(24) & $11$ & $y_{1},y_{5},y_{8}$ & $ \comone\brackets[\big]{\comthrone(y_1, y_5), y_8} = \comone(y_5, y_8)y_1 $ \\
		(25) & $4$ & $y_{4},y_{10},y_{15}$ & $ \comone(y_4, y_{10}) y_{15} = -\comtwo(y_{15}, y_4 y_{10}) $ \\
		(26) & $5$ & $y_{1},y_{4},y_{15}$ & $ \comtwo\brackets[\big]{y_1, \comthrone(y_{15}, y_4)} = \comtwo(y_1, y_4) y_{15} $ \\
		(27) & $3$ & $y_{5},y_{10},y_{15}$ & $ \comtwo(y_5 y_{15}, y_{10}) = \comtwo(y_{15}, y_{10}) y_5 $ \\
		(28) & $13$ & $y_{1},y_{5}$ & $ \comthrtwo(y_1, y_5) = \comthrtwo(y_5, y_1)$ \\
		(29) & $7$ & $y_{1},y_{5},y_{12}$ & $\begin{aligned}
			y_1(y_5y_{12}) &=-\comthrone\brackets[\big]{-\comone(y_1, y_{12}), y_5} \\[-0.1cm]
			& \qquad \mathord{}+\comtwo\brackets[\big]{y_{12}, \comthrtwo(y_1, y_5)}
		\end{aligned}  $ \\
		(30) & $3$ & $y_{3},y_{5},y_{15}$ & $ \comone\brackets[\big]{\comthrone(y_{15}, y_3), y_{15} y_5} = \comthrone\brackets[\big]{y_{15}, \comone(y_3, y_5)} $ \\
		(31) & $3$ & $y_{1},y_{5},y_{15}$ & $ \comthrone\brackets[\big]{y_{15}, \comthrone(y_5, y_1)} = \comthrone(y_{15}y_5, y_1) $ \\
		(32) & $4$ & $y_{1},y_{7},y_{15}$ & $ \comtwo\brackets[\big]{\comtwo(y_{15}, y_7), y_1} =-\comtwo(y_{15}, y_1 y_7) $ \\
		\bottomrule
	\end{tabular}
	}
	\caption{Evaluation of the blueprint identities, part~2. See~\ref{prop:blue-eval2}.}
	\label{fig:eva2}
\end{figure}
\postmidfigure

Before we can prove the remaining identities (33), (34), (35), we have to collect some consequences of (22), \dots, (32).

\begin{prop}
\label{prop:ringhom1}
The map $ \map{-\comtwo(1,\mapdot)}{\prodring}{\prodring}{}{} $ preserves the addition and multiplication of $ \prodring $. In other words, $ \comtwo(1,x+y) = \comtwo(1,x)+\comtwo(1,y) $ and $ \comtwo(1,xy) = -\comtwo(1,x)\comtwo(1,y) $ for all $ x,y \in \prodring $.
\end{prop}
\begin{proof}
	The map is additive by~\ref{lem:psi-id}~\ref{lem:psi-id:biadd} and we have
	\[ \comtwo(1,xy) = -\comtwo\brackets[\big]{\comtwo(1, y), x} = -\comtwo\brackets[\big]{\comtwo(1,y) 1, x} = -\comtwo(1,x) \comtwo(1,y)\]
	for all $ x,y \in \prodring $ by (32) and (27).
\end{proof}

\begin{lemma}\label{lem:g11}
	The following hold:
	\begin{enumerate}[(i)]
		\item $ \assoc{x}{y}{\comtwo(1,1)} = \assoc{x}{\comtwo(1,1)}{y} = \assoc{\comtwo(1,1)}{x}{y} = 0 $ for all $ x,y \in \prodring $ where $ \assoc{x}{y}{z} \defl (xy)z - x(yz) $.
		
		\item $ \comtwo(x,y) \comtwo(1,1) = -\comtwo(x,y) = \comtwo(1,1)\comtwo(x,y) $ for all $ x,y \in \prodring $.
	\end{enumerate}
\end{lemma}
\begin{proof}
	Let $ x,y \in \prodring $. By~(27) and the commutativity of $ \prodring $, we have
	\begin{align*}
		x \brackets[\big]{\comtwo(1,1)y} &= x \comtwo(y,1) = \comtwo(xy,1) = \comtwo(x,1)y = \brackets[\big]{x\comtwo(1,1)}y,
	\end{align*}
	so that $ \assoc{x}{\comtwo(1,1)}{y} = 0 $. It follows in a similar way that $ \assoc{x}{y}{\comtwo(1,1)} = \assoc{\comtwo(1,1)}{x}{y} = 0 $ as well. Further, by two applications of~(27) and one application of (32),
	\begin{align*}
		\comtwo(x,y) \comtwo(1,1) &= \comtwo\brackets[\big]{\comtwo(1,1)x,y} = \comtwo\brackets[\big]{\comtwo(x,1),y} = -\comtwo(x,y).
	\end{align*}
	This finishes the proof.
\end{proof}

\begin{prop}
\label{prop:ringhom2}
The map $ \map{-\comtwo(\mapdot,1)}{\prodring}{\prodring}{}{} $ preserves the addition and multiplication of $ \prodring $. In other words, $ \comtwo(x+y,1) = \comtwo(x,1)+\comtwo(y,1) $ and $ \comtwo(xy,1) = -\comtwo(x,1)\comtwo(y,1) $ for all $ x,y \in \prodring $.
\end{prop}
\begin{proof}
	Again, the map is additive by~\ref{lem:psi-id}~\ref{lem:psi-id:biadd}. Further,
	\begin{align*}
		\comtwo(xy,1) &= \comtwo(y,1)x = -\brackets[\big]{\comtwo(y,1) \comtwo(1,1)} x = -\comtwo(y,1)\brackets[\big]{\comtwo(1,1) x} = -\comtwo(y,1) \comtwo(x,1) 
	\end{align*}
	for all $ x,y \in \prodring $ by~(27) and~\ref{lem:g11}.
\end{proof}

\begin{notation}
	We put $ \prodring_1 \defl \ker(g(\mapdot, 1)) = \{x \in \prodring \mid \comtwo(x,1) = 0\} $ and $ \prodring_2 \defl \ker(\comtwo(1,\mapdot)) = \{ x \in \prodring \mid \comtwo(1,x) = 0\} $. Further, we define maps
	\[ \map{\ringproj_1}{\prodring}{\prodring}{x}{-\comthrone\brackets[\big]{-\comone(x,1), 1}} \quad \text{and} \quad \map{\ringproj_2}{\prodring}{\prodring}{x}{-\comtwo(x,1)}. \]
\end{notation}

\begin{rem}\label{rem:R2-vanish}
	By~\ref{prop:ringhom1} and~\ref{prop:ringhom2}, $ \prodring_1 $ and $ \prodring_2 $ are ideals of $ \prodring $. Further, for all $ x \in \prodring $ and $ y \in \prodring_2 $, we have $ \comtwo(x,y) = \comtwo(1,y)x = 0x = 0 $ by~(27).
\end{rem}

\begin{rem}\label{rem:comthrone-R2}
	It follows from~(26) that $ \comthrone(y_{15}, y_4) \in \prodring_2 $ for all $ y_4 \in \prodring_2 $ and $ y_{15} \in \prodring $.
\end{rem}

\begin{rem}\label{rem:f-is-g}
	By~(25), we have $ \comone(1,x) = -\comtwo(1,x) $ for all $ x \in \prodring $. In particular, $ \prodring_2 = \ker(f(1,\mapdot)) = \ker(f(\mapdot,1)) $.
\end{rem}

\begin{prop}\label{prop:blue-eval3}
	All identities in the column \enquote{Result} of Figure~\ref{fig:eva3} hold for all $ y_1, \ldots, y_{15} \in \prodring $.
\end{prop}
\begin{proof}
	We proceed as in~\ref{prop:blue-eval1} and~\ref{prop:blue-eval2}. For (33), we use in particular~\ref{lem:psi-id}~\ref{lem:psi-id:one-three}, \ref{rem:R2-vanish} and~\ref{rem:comthrone-R2}. For (34), we use~(24), (32) and~\ref{rem:f-is-g}. For~(35), we use~\ref{lem:psi-id}~\ref{lem:psi-id:assoc} and~\ref{lem:psi-id}~\ref{lem:psi-id:zero}.
\end{proof}

\premidfigure
\begin{figure}[htb]
\centering
{\renewcommand{\arraystretch}{2.0}
\begin{tabular}{ccccc}
	\toprule
	No. & \makecell{Blueprint\\identity} & \makecell{Non-zero\\variables} & Result & Remarks \\
	\midrule
	(33) & $6$ & $y_{5},y_{10},y_{15}$ & $ \comtwo\brackets[\big]{\comthrone(y_5, y_{10}), y_{15}} = 0 $ & if $y_5\in\prodring_2$ \\
	(34) & $11$ & $y_{1},y_{5},y_{12}=1$ & \makecell{$ \comtwo\brackets[\big]{\comthrtwo(y_1, y_5), 1} $\\$= \comtwo\brackets[\big]{\brackets[\big]{\comone(y_1, y_5) y_5} y_1, 1} $} & \\
	(35) & $6$ & $y_{1},y_{4},y_{14}$ & \makecell{$ \comtwo\brackets[\big]{\comthrone(y_{14}, y_1), y_4} $\\$= \comthrone(-y_{14}, \comtwo(y_1, y_4)) $} & \\
	\bottomrule
\end{tabular}
}
\caption{Evaluation of the blueprint identities, part~3. See~\ref{prop:blue-eval3}.}
\label{fig:eva3}
\end{figure}
\postmidfigure

It remains to establish that $ \prodring $ decomposes as $ \prodring_1 \oplus \prodring_2 $, that it is associative and to compute explicit formulas for the commutation maps and for the involution~$ \invosym $.

\begin{prop}\label{prop:dirsum}
	We have $ \prodring = \prodring_1 \oplus \prodring_2 $, that is, the ring $ \prodring $ is the direct sum of its ideals $ \prodring_1 $ and $ \prodring_2 $. Further, the projection maps on $ \prodring_1 $ and $ \prodring_2 $ are given by $ \ringproj_1 $ and $ \ringproj_2 $, respectively.
\end{prop}
\begin{proof}
	Let $ x \in \prodring $.
	By~\ref{lem:psi-id}~\ref{lem:psi-id:zero}, $ f(x,1) $ lies in $ \prodring_2 $, so~(33) yields that
	\[ \comtwo\brackets[\big]{\ringproj_1(x), 1} = -\comtwo\brackets[\big]{\comthrone(-f(x,1), 1), 1} = 0. \]
	Hence $ \ringproj_1(x) $ lies in $ \prodring_1 $. Further, $ \ringproj_2(x) = -\comtwo(x,1) $ lies in $ \prodring_2 $ by~\ref{lem:psi-id}~\ref{lem:psi-id:zero}.
	
	Identity~(29) with $ y_{12} = y_5 = 1 $ and $ y_1 = x $ says that
	\[ x = - \comthrone\brackets[\big]{-\comone(x,1), 1} + \comtwo\brackets[\big]{1, \comthrtwo(x,1)} = \ringproj_1(x) + \comtwo\brackets[\big]{1, \comthrtwo(x,1)}. \]
	Since $ \comtwo\brackets[\big]{1, \comthrtwo(x,1)} = -\comone\brackets[\big]{1, \comthrtwo(x,1)} = -\comtwo(x,1) = \ringproj_2(x) $ by~\ref{rem:f-is-g} and~\ref{lem:psi-id}~\ref{lem:psi-id:one-three}, we conclude that $ x = \ringproj_1(x) + \ringproj_2(x) $. In particular, $ \prodring = \prodring_1 + \prodring_2 $.
	
	Finally, let $ x \in \prodring_1 \cap \prodring_2 $. Then $ \ringproj_2(x) = -\comtwo(x,1) = 0 $ by the definition of $ \prodring_1 $. Further,
	\begin{align*}
		\ringproj_1(x) &= - \comthrone\brackets[\big]{-\comone(x,1), 1} = - \comthrone\brackets[\big]{\comtwo(1,x), 1} = -\comthrone(0,1) = 0
	\end{align*}
	by~\ref{rem:f-is-g}. Thus $ x = \ringproj_1(x) + \ringproj_2(x) = 0 $. This finishes the proof.
\end{proof}

\begin{rem}\label{rem:subrings}
	It follows from~\ref{prop:dirsum} that $ \prodring_1 $, $ \prodring_2 $ are nonassociative rings with identity elements $ \ringproj_1(1_\prodring) = -\comthrone(-\comone(1,1), 1) $ and $ \ringproj_2(1_\prodring) = -\comtwo(1,1) = \comone(1,1) $, respectively, and that $ \map{\ringproj_1}{\prodring}{\prodring_1}{}{} $, $ \map{\ringproj_2}{\prodring}{\prodring_2}{}{} $ are homomorphisms.
\end{rem}

\begin{notation}
	We define $ \map{\subisom}{\prodring_1}{\prodring_2}{x}{\comone(1,x) = -\comtwo(1,x)} $.
\end{notation}

\begin{prop}\label{prop:R1R2iso}
	The map $ \subisom $
	is a well-defined isomorphism of nonassociative rings with inverse $ \map{\subisom^{-1}}{\prodring_2}{\prodring_1}{x}{\comthrone(x,1)} $.
\end{prop}
\begin{proof}
	The image of $ \subisom $ lies in $ \prodring_2 $ by~\ref{lem:psi-id}~\ref{lem:psi-id:zero} and it preserves the addition and multiplication by~\ref{prop:ringhom1} and~\ref{rem:subrings}. Further, it follows from~(24) that
	\begin{align*}
		\subisom(1_{\prodring_1}) &= \comone\brackets[\big]{1, -\comthrone\brackets[\big]{-\comone(1,1), 1}} = \comone(1,1) \comone(1,1) = 1_{\prodring_2} 1_{\prodring_2} = 1_{\prodring_2},
	\end{align*}
	so $ \subisom $ is a homomorphism. The kernel of $ \subisom $ is $ \prodring_1 \cap \prodring_2 $, so $ \subisom $ is injective. Further,
	\begin{align*}
		\subisom\brackets[\big]{\comthrone(x,1)} &= \comone\brackets[\big]{\comthrone(x,1), 1} = \comone(1,1) x = 1_{\prodring_2} x = x
	\end{align*}
	for all $ x \in \prodring_2 $ by~(24). Hence $ \subisom $ is an isomorphism with the desired inverse.
\end{proof}

\begin{rem}
	It is \emph{not} true that the map $ \map{}{\prodring_2}{\prodring_1}{x}{\comone(1,x)} $ is an isomorphism. In fact, it is the zero map by~\ref{rem:f-is-g}.
\end{rem}

\begin{prop}
\label{thm:associative}
The multiplication on $\prodring$ is associative.
\end{prop}
\begin{proof}
	For all $ x,y,z \in \prodring_1 $, we have
	\begin{align*}
		\subisom\brackets[\big]{x(yz)} &= \comone\brackets[\big]{1, x(yz)} = \comone\brackets[\big]{1, (xy)z} = \subisom\brackets[\big]{(xy)z}
	\end{align*}
	by~\ref{lem:psi-id}~\ref{lem:psi-id:assoc}. Since $ \map{\subisom}{\prodring_1}{\prodring_2}{}{} $ is an isomorphism by~\ref{prop:R1R2iso}, it follows that $ \prodring_1 $ and $ \prodring_2 $ are associative. Hence $ \prodring = \prodring_1 \oplus \prodring_2 $ is associative as well.
\end{proof}
\begin{prop}
\label{thm:invocomp}
For all $x\in\prodring$, we have $\invo{x}=-\ringproj_1(x)+\ringproj_2(x)$.
\end{prop}
\begin{proof}
	Let $ x \in \prodring_1 $. By~\ref{lem:psi-id}~\ref{lem:psi-id:invo-1}, we have $ \subisom(\invo{x}) = \comone(1,\invo{x}) = -\comone(1,x) = \subisom(-x) $ and $ \invo{\subisom(x)} = \invo{\comone(1,x)} = \comone(1,x) = \subisom(x) $. Since $ \map{\subisom}{\prodring_1}{\prodring_2}{}{} $ is an isomorphism and $ \invosym $ is additive, the assertion follows.
\end{proof}

\begin{notation}\label{not:prod}
	For all $ x,y \in \prodring_1 $, we put $ \ringcoord{x}{y} \defl x + \subisom(y) \in \prodring $. This defines an isomorphism between the (abstract) direct sum $ \prodring_1 \oplus \prodring_1 $ and $ \prodring $ with the property that $ \subisom(\ringcoord{x}{0}) = \ringcoord{0}{x} $ and $ \subisom^{-1}(\ringcoord{0}{x}) = \ringcoord{x}{0} $ for all $ x \in \prodring_1 $.
\end{notation}

\begin{prop}\label{prop:comm-formula}
	Let $ x_1, x_2, y_1, y_2 \in \prodring_1 $. Then the following hold:
	\begin{enumerate}[(i)]
		\item \label{prop:comm-formula:comone}$ \comone(\ringcoord{x_1}{x_2}, \ringcoord{y_1}{y_2}) = \ringcoord{0}{x_1 y_1} $.
		
		\item \label{prop:comm-formula:comtwo}$ \comtwo(\ringcoord{x_1}{x_2}, \ringcoord{y_1}{y_2}) = \ringcoord{0}{-y_1 x_2} $.
		
		\item \label{prop:comm-formula:comthrone}$ \comthrone(\ringcoord{x_1}{x_2}, \ringcoord{y_1}{y_2}) = \ringcoord{x_2 y_1}{x_1 x_2 y_2} $.
		
		\item \label{prop:comm-formula:comthrtwo}$ \comthrtwo(\ringcoord{x_1}{x_2}, \ringcoord{y_1}{y_2}) = \ringcoord{-x_2 y_2}{x_1 x_2 y_1 y_2} $.
	\end{enumerate}
\end{prop}
\begin{proof}
	Put $ x \defl \ringcoord{x_1}{x_2}, y \defl \ringcoord{y_1}{y_2} \in \prodring $. By~\ref{lem:psi-id}~\ref{lem:psi-id:zero} and~\ref{lem:psi-id}~\ref{lem:psi-id:assoc}, we have
	\begin{align*}
		\comone(\ringcoord{x_1}{x_2}, \ringcoord{y_1}{y_2}) &= \comone\brackets[\big]{x_1 + \comone(x_2,1), y_1 + \comone(y_2,1)} = \comone(x_1, y_1) = \comone(1,x_1 y_1) \\
		&= \subisom(x_1 y_1) = \ringcoord{0}{x_1 y_1},
	\end{align*}
	which proves~\ref{prop:comm-formula:comone}. It follows from~(32) and~(25) that
	\begin{align*}
		\comtwo(x,y) &= -\comtwo\brackets[\big]{\comtwo(x,1), y \cdot 1} = \comone(y,1) \comtwo(x,1).
	\end{align*}
	Here $ \comtwo(x,1) = -\ringproj_2(x) = \ringcoord{0}{-x_2} $ and, by~\ref{prop:comm-formula:comone}, $ \comone(y,1) = \ringcoord{0}{y_1} $. Hence~\ref{prop:comm-formula:comtwo} holds. Further, by an application of~(24),
	\begin{align*}
		\ringproj_1\brackets[\big]{\comthrone(x,y)} &= -\comthrone\brackets[\big]{-\comone\brackets{\comthrone(x,y), 1}, 1} = -\comthrone\brackets[\big]{-\comone(y,1) x, 1}.
	\end{align*}
	Here $ -\comone(y,1)x = -\ringcoord{0}{y_1} \ringcoord{x_1}{x_2} = \ringcoord{0}{-y_1 x_2} $. Since $ \comthrone(\mapdot, 1) = \subisom^{-1} $ on $ \prodring_2 $ by~\ref{prop:R1R2iso}, we infer that
	\begin{align*}
		\ringproj_1\brackets[\big]{\comthrone(x,y)} &= -\subisom^{-1}(\ringcoord{0}{-y_1x_2}) = \ringcoord{y_1x_2}{0} = y_1 x_2.
	\end{align*}
	Similarly, by~(35) and~(30),
	\begin{align*}
		\ringproj_2\brackets[\big]{\comthrone(x,y)} &= -\comtwo\brackets[\big]{\comthrone(x,y), 1} = -\comthrone\brackets[\big]{-x, \comtwo(y,1)} = \comthrone\brackets[\big]{-x, \ringcoord{0}{y_2}} \\
		&= \comthrone\brackets[\big]{-x, \comone(1, \ringcoord{y_2}{0})} = - \comone\brackets[\big]{\comthrone(-x,1), x\ringcoord{y_2}{0}}.
	\end{align*}
	Since $ \ringproj_1\brackets{\comthrone(-x,1)} = -x_2 $ by the previous part of the proof, we infer that
	\begin{align*}
		\ringproj_2\brackets[\big]{\comthrone(x,y)} &= -\ringcoord{0}{-x_1 x_2 y_2} = \subisom(x_1 x_2 y_2).
	\end{align*}
	Assertion~\ref{prop:comm-formula:comthrone} follows. Finally, by~\ref{lem:psi-id}~\ref{lem:psi-id:one-three},
	\begin{align*}
		\ringproj_1\brackets[\big]{\comthrtwo(x,y)} &= -\comthrone\brackets[\big]{-\comone\brackets[\big]{\comthrtwo(x,y), 1}, 1} = -\comthrone\brackets[\big]{-y\comtwo(x,1), 1} \\
		&= -\comthrone\brackets{\ringcoord{0}{x_2 y_2}, \ringcoord{1}{1}} = \ringcoord{-x_2 y_2}{0}
	\end{align*}
	and, by~(34),
	\begin{align*}
		\ringproj_2\brackets[\big]{\comthrtwo(x,y)} &= -\comtwo\brackets[\big]{\comthrtwo(x,y),1} = -\comtwo\brackets[\big]{\comone(x,y)yx, 1}  \\
		&= -\comtwo(\ringcoord{0}{x_1 x_2 y_1 y_2} \ringcoord{1}{1}) = \ringcoord{0}{x_1 y_1 x_2 y_2}.
	\end{align*}
	This finishes the proof.
\end{proof}

In the following proof, we apply the strategy from~\ref{rem:conjugate-formulas}.

\begin{prop}\label{prop:comm-formula-H4}
	Assume that $\ell = 4$. Then for all $x_1, x_2, y_1, y_2 \in \prodring_1$,
	\[ \commutator{\risom{\rho_0}(\ringcoord{x_1}{x_2})}{\risom{\rho_1}(\ringcoord{y_1}{y_2})} = \risom{\rho_0 + \rho_1}(\ringcoord{x_1 y_1}{x_2 y_2}). \]
\end{prop}
\begin{proof}
	From the definition of the multiplication on $\prodring$ in~\eqref{eq:def-mult}, we already know that
	\[ \commutator{\risom{\rho_1}(\ringcoord{x_1}{x_2})}{\risom{\rho_2}(\ringcoord{y_1}{y_2})} = \risom{\rho_1 + \rho_2}(\ringcoord{x_1 y_1}{x_2 y_2}). \]
	We want to conjugate this equation by $w_{\rho_0} w_{\rho_1} w_{\rho_2}$ because
	\begin{align*}
		(\rho_1, \rho_2)^{\reflbr{(\rho_0, \rho_1, \rho_2)}} = (\rho_0 +\rho_1, \rho_2)^{\reflbr{(\rho_1, \rho_2)}} = (\rho_0, \rho_1+\rho_2)^{\reflbr{\rho_2}} = (\rho_0, \rho_1).
	\end{align*}
	We can read off from Figures~\ref{fig:parmapex-H4-1} and~\ref{fig:parmapex-H4-2} that
	\begin{align*}
		\inverparbr{\rho_1}{(\rho_0, \rho_1, \rho_2)} &= \inverparbr{\rho_1}{\rho_0} \inverparbr{\rho_0 + \rho_1}{\rho_1} \inverparbr{\rho_0}{\rho_2} \\
		&= (-1,-1) (-1,-1) (1,1) = (1,1), \\
		\inverparbr{\rho_2}{(\rho_0, \rho_1, \rho_2)} &= \inverparbr{\rho_2}{\rho_0} \inverparbr{\rho_2}{\rho_1} \inverparbr{\rho_1 + \rho_2}{\rho_2} \\
		&= (1,1) (-1,-1) (-1,-1) =(1,1), \\
		\inverparbr{\rho_1 + \rho_2}{(\rho_0, \rho_1, \rho_2)} &= \inverparbr{\rho_1 + \rho_2}{\rho_0} \inverparbr{\rho_0 + \rho_1 + \rho_2}{\rho_1} \inverparbr{\rho_0 + \rho_1 + \rho_2}{\rho_2} \\
		&= (-1,1) (1,1) (-1,1) = (1,1).
	\end{align*}
	Hence by conjugating the commutator relation above by $w_{\rho_0} w_{\rho_1} w_{\rho_2}$, we obtain the commutator relation in the assertion.
\end{proof}

\begin{rem}\label{rem:comm-formula:as-in-fold}
	The formulas in~\ref{prop:comm-formula} and~\ref{prop:comm-formula-H4} and the definition of the ring multiplication in~\eqref{eq:def-mult} show that for all $\zeta, \xi \in \paratwopos(\rootbase)$ and all $\rho \in \oprootint{\zeta, \xi}$, the commutation map $\commmap{\zeta, \xi}{\rho}$ of $G$ is the same as in the folding of $\Chev(X, \prodring_1)$ in~\ref{lem:excommrel}.
	Further, it follows from~\ref{ex:twist-is-switch} and~\ref{thm:invocomp} that the involution $ \invosym $ on (the root groups of) $G$ has the same form as in $ \Chev(X, \prodring_1) $.
\end{rem}

\begin{prop}\label{prop:weyl-1}
	For all $\delta \in \rootbase$, we have $w_\delta = \risom{-\delta}(-1_\prodring) \risom{\delta}(1_\prodring) \risom{-\delta}(-1_\prodring)$.
\end{prop}
\begin{proof}
	By the definition of $1_\prodring$ in~\ref{conclusion:para}, we have $w_{\rho_2} = \risom{-\rho_2}(a) \risom{\rho_2}(1_\prodring) \risom{-\rho_2}(b)$ for some $a,b \in \prodring$. Here $a=b$ by ~\ref{prop:a2gradprop}~\ref{prop:a2gradprop:factor-unique} and
	\[ \risom{-\rho_2}(a) = \risom{\rho_2}(1_\prodring)^{w_{\rho_2}} = \risom{\rho_2}(\inverparbr{\rho_2}{\rho_2}.1_\prodring) \]
	by~\ref{prop:a2gradprop}~\ref{prop:a2gradprop:balanced}. Since $\inverparbr{\rho_2}{\rho_2} = (-1,-1)$, it follows that $a=-1_\prodring$.
	
	Put $w_i \defl w_{\rho_i}$ for all $i \in \Set{0,1,2,3}$ (the case $i=0$ only being allowed if $\ell = 4$). For the remaining simple roots, we use the braid relations from~\ref{braidrel}. For example, it follows from $w_1 w_2 w_1 = w_2 w_1 w_2$ that $w_1=w_2^{w_1 w_2}$. Since $\inverparbr{\rho_2}{(\rho_1, \rho_2)} = (1,1) = \inverparbr{-\rho_2}{(\rho_1, \rho_2)}$, we infer that $w_1$ has the desired form. Similarly, $w_0 = w_1^{w_0 w_1}$ and $w_0$ has the desired form (if $\ell = 4$). For $i=3$, we have $w_2 w_3 w_2 w_3 w_2 = w_3 w_2 w_3 w_2 w_3$ and hence $w_3 = w_2^{w_3 w_2 w_3 w_2}$. As $\inverparbr{\rho_2}{(\rho_3, \rho_2, \rho_3, \rho_2)} = (1,1) = \inverparbr{-\rho_2}{(\rho_3, \rho_2, \rho_3, \rho_2)}$, we infer that $w_3$ has the desired form as well.
\end{proof}

We can now state the main result of this paper, which will be reformulated in the language of foldings in Section~\ref{sec:from-fold}.

\begin{theorem}\label{thm:prod-param}
	Let $ \ell \in \Set{3,4} $, let $ G $ be a group with an $ H_\ell $-grading $ (\rootgr{\xi})_{\xi \in H_\ell} $ and fix a family $ (w_\rho)_{\rho \in \rootbase} $ of Weyl elements in $ G $. Then there exist a commutative ring $ \ring $ and a standard coordinatisation of $ G $ by $ \ring \times \ring $ with respect to $ (w_\rho)_{\rho \in \rootbase} $ (in the sense of~\ref{def:stand-coord}).
\end{theorem}
\begin{proof}
	By~\ref{prop:commutative} and~\ref{thm:associative}, $ \ring \defl \prodring_1 $ is a commutative ring. By~\ref{notation:blueh3}, \ref{not:prod} and~\ref{thm:invocomp}, $ (G, (\rootgr{\alpha})_{\alpha \in H_3}) $ is parametrised by the standard parameter system for $ \ring $ with respect to $ \inverparsym $ and $ (w_\rho)_{\rho \in \rootbase} $. The desired commutator relations hold by~\ref{rem:comm-formula:as-in-fold}. Further, the Weyl elements $ (w_\rho)_{\rho \in \rootbase} $ have the desired form by~\ref{prop:weyl-1}.
\end{proof}

\section{All \texorpdfstring{$H_\ell$}{H\_l}-gradings arise from foldings}

\label{sec:from-fold}

\begin{notation}
	In this section, we fix $\ell \in \Set{3,4}$, an $H_\ell$-grading $(\rootgrH{\alpha})_{\alpha \in H_\ell}$ of a group $G$ and a family $(w_\rho)_{\rho \in \rootbaseH}$ of Weyl elements in $G$. We also choose a root base $\rootbaseH$ of $H_\ell$ in standard order, which we denote by $(\rho_0, \rho_1, \rho_2, \rho_3)$ if $\ell = 4$ and by $(\rho_1, \rho_2, \rho_3)$ if $\ell = 3$. As in~\ref{not:fold-goldproj}, we denote by $ \map{\goldfoldproj}{X}{GH_\ell}{}{} $, $\map{\foldproj}{X}{H_\ell}{}{}$ the folding maps from~\ref{pro:fold-abstract-description} and by $ \rootbaseX = (\delta_i)_{i=1}^{2\ell} $ the ordered root base of $ X $ given by Figure~\ref{fig:fold-diags}. By~\ref{thm:prod-param}, there exist a commutative ring $\ring$ and a standard coordinatisation $(\risomH{\zeta})_{\zeta \in H_\ell}$ of $G$ by $\ring \times \ring$ with respect to $(w_\rho)_{\rho \in \rootbaseH}$, which we also fix. We further put $\ring_1 \defl \ring \times \Set{0}, \ring_2 \defl \Set{0} \times \ring \subseteq \ring \times \ring$ and we denote elements $(x,y)$ of $\ring \times \ring$ by $\ringcoord{x}{y}$, in analogy to~\ref{not:prod}.
\end{notation}

It is a corollary of our main result~\ref{thm:prod-param} that the $ H_\ell $-grading $ (\rootgrH{\alpha})_{\alpha \in H_\ell} $ of $ G $ (and hence any $ H_\ell $-grading of any group) is the folding of an $ X $-grading of $ G $ in the sense of~\ref{def:rootgr-fold}. This section deals with the technical details of this deduction. The first intermediate goal is~\ref{prop:opposite-R1R2-comm}.

\begin{lemma}\label{lem:com-form-neg}
	Denote by $(\alpha, \beta, \gamma, \delta, \ep)$ the $H_2$-quintuple associated to $(\rho_2, \rho_3)$. Then for all $ x,y \in \ring \times \ring$, we have
	\begin{align*}
		\commutator{\risomH{-\delta}(x)}{\risomH{\gamma}(y)} &= \risomH{-\ep}\brackets[\big]{\invo{\comthrone(\invo{x}, -y)}} \risomH{\alpha}\brackets[\big]{-\invo{\comthrtwo(\invo{x}, -y)}} \risomH{\beta}\brackets[\big]{-\comthrone(y, \invo{x})}.
	\end{align*}
\end{lemma}
\begin{proof}
	For all $ x,y \in \ring \times \ring $, we know from~\ref{lem:psi-id} that
	\[ \commutator{\risomH{\alpha}(x)}{\risomH{\ep}(y)} = \risomH{\beta}\brackets[\big]{\comthrone(x,y)} \risomH{\gamma}\brackets[\big]{\comthrtwo(x,y)} \risomH{\delta}\brackets[\big]{-\comthrone(-y, x)}. \]
	Using the values of the standard parity map $ \inverparsym $ and~\eqref{eq:H3-param}, we can conjugate this equation by $ w_\alpha w_\ep $. This yields
	\[ \commutator{\risomH{-\delta}(\invo{x})}{\risomH{\gamma}(-y)} = \risomH{-\ep}\brackets[\big]{\invo{\comthrone(x,y)}} \risomH{\alpha}\brackets[\big]{-\invo{\comthrtwo(x,y)}} \risomH{\beta}\brackets[\big]{-\comthrone(-y,x)} \]
	for all $ x,y \in \ring \times \ring $. The assertion follows.
\end{proof}

\begin{lemma}\label{lem:R2-comtwo-zero}
	Let $ (\bar{\alpha}, \bar{\beta}, \bar{\gamma}, \bar{\delta}, \bar{\ep}) $ be any $ H_2 $-quintuple in $ H_\ell $. Then $ \risomH{\bar{\alpha}}(\ring_2) $ and $ \risomH{\bar{\delta}}(\ring_2) $ commute.
\end{lemma}
\begin{proof}
	We have $ \comtwo(\ring_2, \ring_2) = \{0\} $ by~\ref{prop:comm-formula}~\ref{prop:comm-formula:comtwo}. Thus by~\ref{lem:psi-id}~\ref{lem:psi-id:leave-out-2}, $ \commmap{\alpha,\delta}{\beta}(\ring_2, \ring_2) = \Set{0} = \commmap{\alpha,\delta}{\gamma}(\ring_2, \ring_2) $ where $(\alpha, \beta, \gamma, \delta, \ep)$ denotes the $H_2$-quintuple associated to $(\rho_2, \rho_3)$. Hence the assertion holds for $ (\bar{\alpha}, \bar{\beta}, \bar{\gamma}, \bar{\delta}, \bar{\ep}) = (\alpha, \beta, \gamma, \delta, \ep) $. The general assertion follows by conjugation with suitable Weyl elements.
\end{proof}

\begin{rem}\label{rem:hall-witt}
	For all $ x,y,z $ in any group, the \emph{Hall-Witt identity} says that
	\[ \commutator[\big]{\commutator{x}{y}}{z^x} \commutator[\big]{\commutator{z}{x}}{y^z} \commutator[\big]{\commutator{y}{z}}{x^y} = 1. \]
\end{rem}

\begin{prop}\label{prop:opposite-R1R2-comm}
	For all $ \zeta \in H_\ell $, we have $ \commutator{\risomH{\zeta}(\ring_1)}{\risomH{-\zeta}(\ring_2)} = \{1\} $.
\end{prop}
\begin{proof}
	Let $ x,y \in \ring $ and denote by $(\alpha, \beta, \gamma, \delta, \ep)$ the $H_2$-quintuple associated to $(\rho_2, \rho_3)$. By conjugating with suitable Weyl elements, we may assume that $ \zeta = \alpha $. By~\ref{rem:hall-witt}, we have $ ABC = 1_G $ where
	\begin{align*}
		A &= \commutator[\big]{\commutator{\risomH{-\delta}(\ringcoord{0}{1})}{\risomH{\gamma}(\ringcoord{0}{x})}}{\risomH{-\alpha}(\ringcoord{0}{y})^{\risomH{-\delta}(\ringcoord{0}{1})}}, \\
		B &= \commutator[\big]{\commutator{\risomH{-\alpha}(\ringcoord{0}{y})}{\risomH{-\delta}(\ringcoord{0}{1})}}{\risomH{\gamma}(\ringcoord{0}{x})^{\risomH{-\alpha}(\ringcoord{0}{y})}}, \\
		C &= \commutator[\big]{\commutator{\risomH{\gamma}(\ringcoord{0}{x})}{\risomH{-\alpha}(\ringcoord{0}{y})}}{\risomH{-\delta}(\ringcoord{0}{1})^{\risomH{\gamma}(\ringcoord{0}{x})}}.
	\end{align*}
	It follows from~\ref{lem:R2-comtwo-zero} that $ B=1=C $ and that $ \risomH{-\alpha}(\ringcoord{0}{y})^{\risomH{-\delta}(\ringcoord{0}{1})} = \risomH{-\alpha}(\ringcoord{0}{y}) $. Further, since $ \comthrone(\ring_2, \ring_2) = \{0\} $ by~\ref{prop:comm-formula}~\ref{prop:comm-formula:comthrone}, we infer from~\ref{lem:com-form-neg} that
	\begin{align*}
		\commutator{\risomH{-\delta}(\ringcoord{0}{1})}{\risomH{\gamma}(\ringcoord{0}{x})} &= \risomH{\alpha}\brackets[\big]{-\invo{\comthrtwo(\ringcoord{0}{1}, \ringcoord{0}{-x})}} = \risomH{\alpha}(\ringcoord{x}{0}).
	\end{align*}
	Altogether,
	\begin{align*}
		1_G &= ABC = A = \commutator{\risomH{\alpha}(\ringcoord{x}{0})}{\risomH{-\alpha}(\ringcoord{0}{y})}.
	\end{align*}
	This finishes the proof.
\end{proof}

We now define the $ X $-grading of $ G $ whose folding is $ (\rootgrH{\alpha})_{\alpha \in H_\ell} $.

\begin{definition}\label{def:unfold-rootgr}
	For all $ \xi \in X $, we define a root homomorphism $ \map{\risomX{\xi}}{(\ring, +)}{\rootgrH{\foldproj(\xi)}}{}{} $ by $ \risomX{\xi}(s) \defl \risomH{\foldproj(\xi)}(\ringcoord{s}{0}) $ if $ \goldfoldproj(\xi) $ is short in $ GH_\ell $ and by $ \risomX{\xi}(s) \defl \risomH{\foldproj(\xi)}(\ringcoord{0}{s}) $ if $ \goldfoldproj(\xi) $ is long in $ GH_\ell $. Further, we put $ \rootgrunfold{\xi} \defl \risomX{\xi}(\ring) $.
\end{definition}

\begin{lemma}\label{lem:unfold-weyl}
	Let $ (i,j) \in \{(1,6), (2,4), (3,5)\} $ (or possibly $(i,j) = (7,8)$ if $X = E_8$) and put $ \rho \defl \foldproj(\delta_i) = \foldproj(\delta_j) \in \rootbase \subset H_\ell $. Then for all $ \xi \in X $, we have $ (\rootgrunfold\xi)^{w_\rho} = \rootgrunfold{\xi^{\reflbr{\delta_i} \reflbr{\delta_j}}} $.
\end{lemma}
\begin{proof}
	Let $\xi \in X$. Put $k \defl 1$ if $\goldfoldproj(\xi)$ is short in $GH_\ell$ and $k \defl 2$ otherwise, so that $\rootgrunfold{\xi} = \risomH{\foldproj(\xi)}(\ring_k)$. Then
	\begin{align*}
		(\rootgrunfold{\xi})^{w_\rho} = \risomH{\foldproj(\xi)}(\ring_k)^{w_\rho} =\risomH{\foldproj(\xi)^{\reflbr{\rho}}}(\inverpar{\foldproj(\xi)}{\rho}.\ring_k) = \risomH{\foldproj(\xi)^{\reflbr{\rho}}}(\ring_k).
	\end{align*}
	Here $\foldproj(\xi)^{\reflbr{\rho}} = \foldproj(\xi^{\weylemb(\reflbr{\rho})}) = \foldproj(\xi^{\reflbr{\delta_i} \reflbr{\delta_j}})$ by~\ref{pro:fold-abstract-description}~\ref{pro:fold-abstract-description:foldproj}. Further, by~\ref{note:preim-order}, $\goldfoldproj(\xi^{\reflbr{\delta_i} \reflbr{\delta_j}})$ is short in $GH_\ell$ if and only if $\goldfoldproj(\xi)$ is, so that
	\[ \rootgrunfold{\xi^{\reflbr{\delta_i} \reflbr{\delta_j}}} = \risomH{\foldproj(\xi^{\reflbr{\delta_i} \reflbr{\delta_j}})}(\ring_k). \]
	Hence
	\begin{align*}
		(\rootgrunfold{\xi})^{w_\rho} &= \risomH{\foldproj(\xi^{\reflbr{\delta_i} \reflbr{\delta_j}})}(\ring_k) = \rootgrunfold{\xi^{\reflbr{\delta_i} \reflbr{\delta_j}}},
	\end{align*}
	as desired.
\end{proof}

\begin{lemma}\label{lem:D6-commrel-parabolic}
	Let $ \xi_1, \xi_2 \in X $ and assume that $\foldproj(\xi_1)\ne \foldproj(\xi_2)$ and that $\foldproj(\xi_1)$, $\foldproj(\xi_2)$ both lie in $\paratwopos(\rootbase)$, the set defined in~\ref{def:paratwopos}. If $ \xi_1 + \xi_2 $ is a root (in $ X $), then
	\[ \commutator{\risomX{\xi_1}(r)}{\risomX{\xi_2}(s)} = \risomX{\xi_1 +\xi_2}(c_{\xi_1, \xi_2} rs) \]
	for all $r,s \in \ring$ where $c_{\xi_1, \xi_2} \in \Set{\pm 1}$ is the same sign that appears in the Chevalley group $\Chev(X, \ring)$ in~\ref{rem:chev-comm-sign}. If $\xi_1+\xi_2$ is not a root (in $X$), then $ \commutator{\rootgrunfold{\xi_1}}{\rootgrunfold{\xi_2}} = \{1\} $.
\end{lemma}
\begin{proof}
	The $ X $-root groups $ \rootgrunfold{\xi_1}, \rootgrunfold{\xi_2} $ are contained in the $ H_\ell $-root groups $ \rootgrH{\zeta_1} $, $ \rootgrH{\zeta_2} $, respectively. For the latter, we have a commutator formula $ \commutator{\rootgrH{\zeta_1}}{\rootgrH{\zeta_2}} \subseteq \rootgrH{\oprootint{\zeta_1, \zeta_2}} $. By~\ref{rem:comm-formula:as-in-fold}, the commutator maps $ \commmap{\zeta_1, \zeta_2}{\rho} $ in this commutator formula are of the same form as in $ (\Chev(X, \ring), (\rootgrD{\rho})_{\rho \in D_6}) $. The assertion follows.
\end{proof}

\begin{example}
	We illustrate the argument in the proof of~\ref{lem:D6-commrel-parabolic} by the example $ \xi_1 \defl e_2-e_3 $, $ \xi_2 \defl e_3-e_4 \in D_6 $. Then $ \foldproj(\xi_1) = \rho_2 $ and $ \foldproj(\xi_2) = \rho_3 $ by Figure~\ref{fig:foldrootsd6}. Denote by $ (\alpha, \beta, \gamma, \delta, \ep) $ the $ H_2 $-quintuple associated to $ (\rho_2, \rho_3) $. The commutator formulas in Figure~\ref{fig:excommrel}, which hold by~\ref{thm:prod-param}, yield that
	\begin{align*}
		\commutator[\big]{\risomX{\xi_1}(a)}{\risomX{\xi_2}(b)} = \commutator[\big]{\risomH{\alpha}(\ringcoord{a}{0})}{\risomH{\ep}(\ringcoord{0}{b})} = \risomH{\delta}(\ringcoord{ab}{0}) 
	\end{align*}
	for all $ a,b \in \ring_1 $. Since $ \delta = \rootcoord{0,1,\gold} $, we read off from Figure~\ref{fig:foldrootsd6} that the short root in $ \scalmap^{-1}(\delta) $ is $ \goldfoldproj(e_2 - e_4) $. We conclude that $ \commutator{\rootgrunfold{\xi_1}}{\rootgrunfold{\xi_2}} \le \rootgrunfold{\xi_1 + \xi_2} $.
\end{example}

\begin{lemma}\label{lem:unfold-comm}
	The group $ G $ has $ X $-commutator relations with respect to the root groups $ (\rootgrunfold{\xi})_{\xi \in X} $ from~\ref{def:unfold-rootgr} and has the same commutation maps (with respect to $ (\risomX{\xi})_{\xi \in X} $) as the Chevalley group $ \Chev(X, \ring_1) $ (with respect to the root isomorphisms constructed in~\ref{rem:ex-twist-signs}).
\end{lemma}
\begin{proof}
	Let $ \xi_1, \xi_2 \in X $ be non-proportional and put $ \zeta_1 \defl \goldfoldproj(\xi_1) $, $ \zeta_2 \defl \goldfoldproj(\xi_2) \in GH_3 $. If $ \zeta_1 = e\zeta_2 $ for some $ e \in \{\pm 1\} $, then $ \xi_1 = e \xi_2 $ because $ \map{\goldfoldproj}{X}{GH_\ell}{}{} $ is bijective and linear, which contradicts the assumption that $ \xi_1, \xi_2 $ are non-proportional. If $\zeta_1$, $\zeta_2$ differ by a factor of $\gold$, then $\rootgrunfold{\xi_1}$ and $\rootgrunfold{\xi_2}$ are contained in the same (abelian) $H_\ell$-root group $\rootgrH{\foldproj(\xi_1)} = \rootgrH{\foldproj(\xi_2)}$, so that $\commutator{\rootgrunfold{\xi_1}}{\rootgrunfold{\xi_2}} = \Set{1_G}$. If $ \zeta_1$, $ \zeta_2 $ differ by a factor of $-\gold$, then $\foldproj(\xi_1) = -\foldproj(\xi_2)$ and $ \commutator{\rootgrunfold{\xi_1}}{\rootgrunfold{\xi_2}} = \{1_G\} $ by~\ref{prop:opposite-R1R2-comm}. Hence we can assume that $ \zeta_1, \zeta_2 $ are non-proportional.
	
	If $\zeta_1$, $\zeta_2$ both lie in $\paratwopos(\rootbase)$, then the assertion holds by~\ref{lem:D6-commrel-parabolic}. In general, there exists (by~\ref{rem:paratwopos}) some $w \in \Weyl(H_\ell)$ such that $\zeta_1^w$ and $\zeta_2^w$ both lie in $\paratwopos(\rootbase)$. Thus we can conjugate the commutator relation of the pair $(\rootgrunfold{\zeta_1^w}, \rootgrunfold{\zeta_2^w})$ by a suitable sequence of the Weyl elements $(w_\delta')_{\delta' \in \rootbase}$ to obtain the commutator relation of $(\rootgrunfold{\zeta_1}, \rootgrunfold{\zeta_2})$. Since the parametrisations of $G$ and $\Chev(X, \ring_1)$ adhere to the same parity map (namely, the standard parity map for $H_\ell$), we infer that the commutator relation of $(\rootgrunfold{\zeta_1}, \rootgrunfold{\zeta_2})$ has the desired form as well (cf.~\ref{rem:conjugate-formulas}).
\end{proof}

\begin{theorem}\label{thm:from-folding}
	Let $\ell \in \Set{3,4}$, let $X $ denote the root system in $ \Set{D_6, E_8}$ of rank $2\ell$ and let $(\rootgrH{\zeta})_{\zeta \in H_\ell}$ be an $H_\ell$-grading of a group $G$. Then the family $ (\rootgrunfold{\xi})_{\xi \in X} $ defined in~\ref{def:unfold-rootgr} is an $ X $-grading of $ G $, and its folding (in the sense of~\ref{def:rootgr-fold}) is precisely $ (\rootgrH{\zeta})_{\zeta \in H_\ell} $.
\end{theorem}
\begin{proof}
	By~\ref{lem:unfold-comm}, $G$ has $X$-commutator relations with root groups $(\rootgrunfold{\xi})_{\xi \in X}$. Even more, the same commutator formulas as in the $X$-grading of $\Chev(X, \ring)$ hold. This implies that for all $\alpha_1, \alpha_2 \in X$ with $\goldfoldproj(\alpha_2) = \gold \goldfoldproj(\alpha_1) $, the elements
	\begin{align*}
		w_{\alpha_1} &\defl \risomH{-\foldproj(\alpha_1)}(-1,0) \risomH{\foldproj(\alpha_1)}(1,0) \risomH{-\foldproj(\alpha_1)}(-1,0) = \risomX{-\alpha_1}(-1) \risomX{\alpha_1}(1) \risomX{-\alpha_1}(-1), \\
		w_{\alpha_2} &\defl \risomH{-\foldproj(\alpha_1)}(0,-1) \risomH{\foldproj(\alpha_1)}(0,1) \risomH{-\foldproj(\alpha_1)}(0,-1) = \risomX{-\alpha_2}(-1) \risomX{\alpha_2}(1) \risomX{-\alpha_2}(-1)
	\end{align*}
	are Weyl elements with respect to $ (\rootgrunfold{\xi})_{\xi \in X} $. (See \cite[5.6.6]{torben} for the details of this conclusion.) Further, since $\foldproj(\rootbaseX) = \rootbaseH$, we have $\rootgrunfold{\possysX} = \rootgrH{\possysH}$ where $\possysX$, $\possysH$ denote the positive systems corresponding to $\possysX$ and $\possysH$, respectively. Hence $(\rootgrunfold{\xi})_{\xi \in X}$ satisfies Axiom~\ref{def:rgg}~\ref{def:rgg:pos}. We infer that $ (\rootgrunfold{\xi})_{\xi \in X} $ is an $ X $-grading of $ G $. Since $\rootgrunfold{\alpha_1} \rootgrunfold{\alpha_2} = \rootgrH{\foldproj(\alpha_1)}$ for all $\alpha_1, \alpha_2 \in X$ with $\goldfoldproj(\alpha_2) = \gold \goldfoldproj(\alpha_1)$, its folding is precisely $(\rootgrH{\beta})_{\beta \in H_\ell}$.
\end{proof}

\appendix

\section{Commutator relations in foldings}

\renewcommand{\thefigure}{\arabic{figure}} 

\begin{notation}
	In this section, the notation of Section~\ref{sec:fold} (that is,~\ref{not:root-base-order} and~\ref{not:fold-goldproj}) continues to hold. Further, we fix an arbitrary $ X $-graded group $ (G, (\rootgrX{\alpha})_{\alpha \in X}) $.
\end{notation}

In~\ref{folding-is-rgg}, we claimed that the folding $ (\rootgrH{\beta})_{\beta \in H_\ell} $ of $ (\rootgrX{\alpha})_{\alpha \in X} $ (as defined in~\ref{def:rootgr-fold}) satisfies $ H_\ell $-commutator relations, but we only proved this for the special case that $ G $ is a (certain) Chevalley group. The general assertion is essentially a consequence of \cite[2.1.29]{torben}. Unfortunately, \cite[2.1.29]{torben} is not stated in the generality that is needed in our situation, but the same arguments apply. The purpose of this appendix is to formulate a variation of \cite[2.1.29]{torben} in sufficient generality (see~\ref{lem:fold-comm-general}) and then to apply it, in~\ref{pro:fold-comm}~\ref{pro:fold-comm:H}, to complete the proof of~\ref{folding-is-rgg}.

At first, we check that some assumptions in~\ref{lem:fold-comm-general} are satisfied.

\begin{lemma}\label{lem:foldproj-possys:GH}
	Let $ \alpha, \beta \in GH_\ell $ be non-proportional. Then there exists a positive system $ \possys $ in $ X $ which contains $ \goldfoldproj^{-1}(\Set{\alpha,\beta}) $.
\end{lemma}
\begin{proof}
	Since $\goldfoldproj$ is bijective, $ \goldfoldproj^{-1}(\Set{\alpha,\beta}) $ has exactly two elements $\alpha'$, $\beta'$. These two elements must be non-proportional because $\alpha, \beta$ are non-proportional. Thus there exists a positive system in $X$ containing $\Set{\alpha', \beta'}$.
\end{proof}

\begin{lemma}\label{lem:foldproj-possys:H}
	Let $ \alpha, \beta \in H_\ell $ be non-proportional. Then there exists a positive system $ \possys $ in $ X $ which contains $ \foldproj^{-1}(\Set{\alpha,\beta}) $.
\end{lemma}
\begin{proof}
	If $\alpha, \beta$ are contained in the positive system $\possysH$ of $H_\ell$ corresponding to the root base $\rootbaseH$, then the positive system $\possysX$ in $X$ corresponding to the root base $\rootbaseX$ contains $ \foldproj^{-1}(\Set{\alpha,\beta}) $ (because $\possysX = \foldproj^{-1}(\possysH)$). In general, there exists $w \in \Weyl(H_\ell)$ such that $\alpha, \beta$ are contained in $\possysH^w$. Then $\possysX^{\weylemb(w)} = \foldproj^{-1}(\possysH)^{\weylemb(w)} = \foldproj^{-1}(\possysH^w)$ by~\ref{lem:preim-trans}, so $\possysX^{\weylemb(w)}$ contains $\foldproj^{-1}(\Set{\alpha, \beta})$.
\end{proof}

We will need the following generalisation of the notion of root intervals.

\begin{definition}
	Let $\roots$ be a root system and let $ A,B \subseteq \roots $. We denote by $ \oprootint{A,B} $ (respectively, $ \oprootintcry{A,B} $) the set of all roots $ \gamma \in \roots $ which can be written as $ \gamma = \sum_{i=1}^n \lambda_i \alpha_i + \sum_{j=1}^m \mu_j \beta_j $ for $ n,m \in \Z_{>0} $, $ \alpha_1, \ldots, \alpha_n \in A $, $ \beta_1, \ldots, \beta_m \in B $ and $ \lambda_1, \ldots, \lambda_n, \mu_1, \ldots, \mu_m $ positive real (respectively, integral) numbers.
\end{definition}

Now we can state the following slight generalisation of \cite[2.1.29]{torben}.

\begin{lemma}\label{lem:fold-comm-general}
	Let $ \map{\bar{\pi}}{\roots}{\roots'}{}{} $ be a surjective map between root systems $ \roots, \roots' $ such that $ \roots $ is reduced. Assume that for all non-proportional $ \alpha', \beta' \in \roots' $, there exists a positive system $ \possys $ in $ \roots $ which contains $ \bar{\pi}^{-1}(\Set{\alpha', \beta'}) $. Let $ G $ be a group which has $ \roots $-commutator relations with root groups $ (\bar{V}_{\alpha})_{\alpha \in \roots} $ and put $ \bar{U}_{\alpha'} \defl \bar{V}_{\bar{\pi}^{-1}(\Set{\alpha'})} $ for all $ \alpha' \in \roots' $. Assume also that all groups $ (\bar{V}_{\alpha})_{\alpha \in \roots} $ are abelian. If
	\begin{equation}\label{eq:comp-rootint}
		\bar{\pi}(\oprootint{A,B}) \subseteq \oprootint{\bar{\pi}(A), \bar{\pi}(B)} \qquad \text{for all } A,B \subseteq \roots,
	\end{equation}
	then $\commutator{\bar{U}_{\alpha'}}{\bar{U}_{\bar{\beta}'}} \subseteq \bar{U}_{\oprootint{\alpha', \beta'}}$ for all non-proportional $\alpha', \beta' \in \roots'$.
	If the commutator relations of $ G $ are crystallographic and
	\begin{equation}\label{eq:comp-rootintcry}
		\bar{\pi}(\oprootintcry{A,B}) \subseteq \oprootintcry{\bar{\pi}(A), \bar{\pi}(B)} \qquad \text{for all } A,B \subseteq \roots,
	\end{equation}
	then $\commutator{\bar{U}_{\alpha'}}{\bar{U}_{\bar{\beta}'}} \subseteq \bar{U}_{\oprootintcry{\alpha', \beta'}}$ for all non-proportional $\alpha', \beta' \in \roots'$.
\end{lemma}
\begin{proof}
	This follows from the same arguments as in \cite[2.1.29]{torben}.
\end{proof}

\begin{rem}\label{rem:fold-comm-general-crit}
	Since $\map{\goldfoldproj}{X}{H_\ell}{}{}$ is (induced by) a linear map, it satisfies~\eqref{eq:comp-rootint} and~\eqref{eq:comp-rootintcry}. The scaling map $\map{\scalmap}{GH_\ell}{H_\ell}{}{}$ also satisfies~\eqref{eq:comp-rootint}, but as it scales some roots by the non-integral factor $\gold^{-1}$, it does not satisfy~\eqref{eq:comp-rootintcry}. Hence $\map{\foldproj = \scalmap \circ \goldfoldproj}{X}{H_\ell}{}{}$ satisfies~\eqref{eq:comp-rootint} but not~\eqref{eq:comp-rootintcry}. Further, by~\ref{lem:foldproj-possys:GH} and~\ref{lem:foldproj-possys:H}, $\foldproj$ and $\goldfoldproj$ satisfy the requirement on preimages in~\ref{lem:fold-comm-general}.
\end{rem}

\begin{prop}\label{pro:fold-comm}
	Let $G$ be a group which has $X$-commutator relations with root groups $(\rootgrX{\alpha})_{\alpha \in X}$. Then the following hold:
	\begin{enumerate}[(a)]
		\item \label{pro:fold-comm:H}$G$ has $H_\ell$-commutator relations with root groups $(\rootgrH{\beta})_{\beta \in H_\ell}$, as defined in~\ref{def:rootgr-fold}.
		
		\item \label{pro:fold-comm:GH}Put $\rootgrGH{\beta} \defl \rootgrX{\goldfoldproj^{-1}(\beta)}$ for all $\beta \in GH_\ell$. Then $G$ has crystallographic $GH_\ell$-commutator relations with root groups $(\rootgrGH{\beta})_{\beta \in GH_\ell}$. By this we mean that all root groups $(\rootgrGH{\beta})_{\beta \in GH_\ell}$ are abelian, that $\commutator{\rootgrGH{\beta}}{\rootgrGH{\gold \beta}} = \Set{1}$ for all short roots $\beta \in GH_\ell$ and that $\commutator{\rootgrGH{\alpha}}{\rootgrGH{\beta}} \subseteq \rootgrGH{\oprootintcry{\alpha, \beta}}$ for all non-proportional $\alpha, \beta \in GH_\ell$.
	\end{enumerate}
\end{prop}
\begin{proof}
	By~\ref{rem:fold-comm-general-crit} and~\ref{lem:fold-comm-general}, both families $(\rootgrH{\beta})_{\beta \in H_\ell}$ and $(\rootgrGH{\beta})_{\beta \in GH_\ell}$ satisfy the commutator axiom for non-proportional roots and $(\rootgrGH{\beta})_{\beta \in GH_\ell}$ even satisfies the crystallographic version of this axiom. Further, $\gen{\rootgrGH{\beta}, \rootgrGH{\gold \beta}}$ is abelian by~\ref{rem:fold-rootgr-commute} for all short roots $\beta \in GH_\ell$, which proves the remaining part of the assertion.
\end{proof}

\begin{rem}
	One can show that $ (\rootgrGH{\beta})_{\beta \in GH_\ell} $ does not have Weyl elements, so it is not a \enquote{$ GH_\ell $-grading of $ G $}.
\end{rem}

\begin{rem}
\label{rem:exgh3commrel}
	The commutator relations described in~\ref{pro:fold-comm}~\ref{pro:fold-comm:GH} are strictly stronger than the ones in~\ref{pro:fold-comm}~\ref{pro:fold-comm:GH} because they are crystallographic. We illustrate this by the following computation. Let $ (\alpha, \beta, \gamma, \delta, \ep) $ be an $ H_2 $-quintuple in $ GH_\ell $. Then $ \rootgrH{\alpha} = \rootgrGH{\alpha} \rootgrGH{\tau \alpha} $ and $ \rootgrH{\delta} = \rootgrGH{\delta} \rootgrGH{\tau \delta} $, so arbitrary $ y_{\alpha} \in \rootgrH{\alpha} $ and $ y_{\delta} \in \rootgrH{\delta} $ can be written as $ y_{\alpha} = x_\alpha x_{\tau \alpha} $ and $ y_{\delta'} = x_\delta x_{\tau \delta} $ where $ x_\xi \in \rootgrGH{\xi} $ for all $ \xi \in \Set{\alpha, \tau\alpha, \delta, \tau \delta} $. By relations~\ref{rem:commrel}~\ref{rem:commrel:add-left} and~\ref{rem:commrel}~\ref{rem:commrel:add-right},
	\begin{align*}
		\commutator{y_{\alpha}}{y_{\delta}} &= \commutator{x_\alpha x_{\tau \alpha}}{x_\delta x_{\tau \delta}} = \commutator{x_\alpha}{x_\delta x_{\tau \delta}}^{x_{\tau \delta}} \commutator{x_{\tau \alpha}}{x_\delta x_{\tau \delta}} \\
		&= \commutator{x_\alpha}{x_{\tau \delta}}^{x_{\tau \alpha}} \commutator{x_\alpha}{x_\delta}^{x_{\tau \delta} x_{\tau \alpha}} \commutator{x_{\tau \alpha}}{x_{\tau \delta}} \commutator{x_{\tau \alpha}}{x_{\delta}}^{x_{\tau \delta}}.
	\end{align*}
	Using the crystallographic $ GH_\ell $-commutator relations and the crystallographic adjacency relations that we can read off from Figure~\ref{fig:gh2}, we infer that
	\[ \commutator{y_{\alpha}}{y_{\delta}} = \commutator{x_{\alpha}}{x_{\tau \delta}} \commutator{x_{\tau \alpha}}{x_\delta} \in \rootgrGH{\tau \gamma} \rootgrGH{\tau \beta}. \]
	We conclude that $ \commutator{\rootgrH{\alpha}}{\rootgrH{\delta}} \subseteq \rootgrGH{\tau \gamma} \rootgrGH{\tau \beta} $, which is a refinement of the non-crystallographic commutator relation
	\[ \commutator{\rootgrH{\alpha}}{\rootgrH{\delta}} \subseteq \rootgrH{\beta} \rootgrH{\gamma} = \rootgrGH{\beta} \rootgrGH{\tau \beta} \rootgrGH{\gamma} \rootgrGH{\tau \gamma}. \]
	Similar computations can be performed for all other pairs of non-proportional roots. See also the commutator relations in Figure~\ref{fig:excommrel}.
	
	Note that if $ \ell \ge 3 $, then every $ H_\ell $-grading is the folding of an $ X $-grading by~\ref{thm:from-folding}. Hence the observations in the previous paragraph apply to all $ H_\ell $-graded groups in this case.
\end{rem}

\begin{rem}
	For the special case $ G = \Chev(X, \ring) $ (for any commutative ring~$ \ring $), the assertions of~\ref{pro:fold-comm}~\ref{pro:fold-comm:GH} and~\ref{rem:exgh3commrel} can also be read off from the commutator relations in Figure~\ref{fig:excommrel}.
\end{rem}

\premidfigure
\begin{figure}[htb]
	\centering\begin{tabular}{ccc}
		\toprule
		$ \beta $ & $ \alpha_1 $ & $ \alpha_2 $ \\
		\midrule
		$ \rootcoord{0, 0, 0, 1} $ & $ \rootcoord{0, 0, 0, 0, 1, 0, 0, 0} $ & $ \rootcoord{0, 0, 1, 0, 0, 0, 0, 0} $ \\
		$ \rootcoord{0, 0, 1, 0} $ & $ \rootcoord{0, 1, 0, 0, 0, 0, 0, 0} $ & $ \rootcoord{0, 0, 0, 1, 0, 0, 0, 0} $ \\
		$ \rootcoord{0, 1, 0, 0} $ & $ \rootcoord{1, 0, 0, 0, 0, 0, 0, 0} $ & $ \rootcoord{0, 0, 0, 0, 0, 1, 0, 0} $ \\
		$ \rootcoord{1, 0, 0, 0} $ & $ \rootcoord{0, 0, 0, 0, 0, 0, 1, 0} $ & $ \rootcoord{0, 0, 0, 0, 0, 0, 0, 1} $ \\
		$ \rootcoord{0, 0, \gold, 1} $ & $ \rootcoord{0, 0, 0, 1, 1, 0, 0, 0} $ & $ \rootcoord{0, 1, 1, 1, 0, 0, 0, 0} $ \\
		$ \rootcoord{0, 0, \gold, \gold} $ & $ \rootcoord{0, 0, 1, 1, 0, 0, 0, 0} $ & $ \rootcoord{0, 1, 1, 1, 1, 0, 0, 0} $ \\
		$ \rootcoord{0, 0, 1, \gold} $ & $ \rootcoord{0, 1, 1, 0, 0, 0, 0, 0} $ & $ \rootcoord{0, 0, 1, 1, 1, 0, 0, 0} $ \\
		$ \rootcoord{0, 1, 1, 0} $ & $ \rootcoord{1, 1, 0, 0, 0, 0, 0, 0} $ & $ \rootcoord{0, 0, 0, 1, 0, 1, 0, 0} $ \\
		$ \rootcoord{1, 1, 0, 0} $ & $ \rootcoord{1, 0, 0, 0, 0, 0, 1, 0} $ & $ \rootcoord{0, 0, 0, 0, 0, 1, 0, 1} $ \\
		$ \rootcoord{0, \gold, \gold, 1} $ & $ \rootcoord{0, 0, 0, 1, 1, 1, 0, 0} $ & $ \rootcoord{1, 1, 1, 1, 0, 1, 0, 0} $ \\
		$ \rootcoord{0, \gold, \gold, \gold} $ & $ \rootcoord{0, 0, 1, 1, 0, 1, 0, 0} $ & $ \rootcoord{1, 1, 1, 1, 1, 1, 0, 0} $ \\
		$ \rootcoord{0, 1, 1, \gold} $ & $ \rootcoord{1, 1, 1, 0, 0, 0, 0, 0} $ & $ \rootcoord{0, 0, 1, 1, 1, 1, 0, 0} $ \\
		$ \rootcoord{1, 1, 1, 0} $ & $ \rootcoord{1, 1, 0, 0, 0, 0, 1, 0} $ & $ \rootcoord{0, 0, 0, 1, 0, 1, 0, 1} $ \\
		$ \rootcoord{\gold, \gold, \gold, 1} $ & $ \rootcoord{0, 0, 0, 1, 1, 1, 0, 1} $ & $ \rootcoord{1, 1, 1, 1, 0, 1, 1, 1} $ \\
		$ \rootcoord{\gold, \gold, \gold, \gold} $ & $ \rootcoord{0, 0, 1, 1, 0, 1, 0, 1} $ & $ \rootcoord{1, 1, 1, 1, 1, 1, 1, 1} $ \\
		$ \rootcoord{0, \gold, \gold^2, \gold} $ & $ \rootcoord{0, 1, 1, 1, 0, 1, 0, 0} $ & $ \rootcoord{1, 1, 1, 2, 1, 1, 0, 0} $ \\
		$ \rootcoord{0, 1, \gold^2, \gold} $ & $ \rootcoord{1, 1, 1, 1, 0, 0, 0, 0} $ & $ \rootcoord{0, 1, 1, 2, 1, 1, 0, 0} $ \\
		$ \rootcoord{1, 1, 1, \gold} $ & $ \rootcoord{1, 1, 1, 0, 0, 0, 1, 0} $ & $ \rootcoord{0, 0, 1, 1, 1, 1, 0, 1} $ \\
		$ \rootcoord{\gold, \gold, \gold^2, \gold} $ & $ \rootcoord{0, 1, 1, 1, 0, 1, 0, 1} $ & $ \rootcoord{1, 1, 1, 2, 1, 1, 1, 1} $ \\
		$ \rootcoord{0, \gold, 2\gold, \gold^2} $ & $ \rootcoord{0, 0, 1, 2, 1, 1, 0, 0} $ & $ \rootcoord{1, 2, 2, 2, 1, 1, 0, 0} $ \\
		$ \rootcoord{0, \gold, \gold^2, \gold^2} $ & $ \rootcoord{0, 1, 1, 1, 1, 1, 0, 0} $ & $ \rootcoord{1, 1, 2, 2, 1, 1, 0, 0} $ \\
		$ \rootcoord{0, 1, \gold^2, \gold^2} $ & $ \rootcoord{1, 1, 1, 1, 1, 0, 0, 0} $ & $ \rootcoord{0, 1, 2, 2, 1, 1, 0, 0} $ \\
		$ \rootcoord{1, 1, \gold^2, \gold} $ & $ \rootcoord{1, 1, 1, 1, 0, 0, 1, 0} $ & $ \rootcoord{0, 1, 1, 2, 1, 1, 0, 1} $ \\
		$ \rootcoord{\gold, \gold, 2\gold, \gold^2} $ & $ \rootcoord{0, 0, 1, 2, 1, 1, 0, 1} $ & $ \rootcoord{1, 2, 2, 2, 1, 1, 1, 1} $ \\
		$ \rootcoord{\gold, \gold, \gold^2, \gold^2} $ & $ \rootcoord{0, 1, 1, 1, 1, 1, 0, 1} $ & $ \rootcoord{1, 1, 2, 2, 1, 1, 1, 1} $ \\
		$ \rootcoord{\gold, \gold^2, \gold^2, \gold} $ & $ \rootcoord{1, 1, 1, 1, 0, 1, 0, 1} $ & $ \rootcoord{1, 1, 1, 2, 1, 2, 1, 1} $ \\
		$ \rootcoord{1, 1, \gold^2, \gold^2} $ & $ \rootcoord{1, 1, 1, 1, 1, 0, 1, 0} $ & $ \rootcoord{0, 1, 2, 2, 1, 1, 0, 1} $ \\
		$ \rootcoord{1, \gold^2, \gold^2, \gold} $ & $ \rootcoord{1, 1, 1, 1, 0, 1, 1, 0} $ & $ \rootcoord{1, 1, 1, 2, 1, 2, 0, 1} $ \\
		$ \rootcoord{\gold, 2\gold, 2\gold, \gold^2} $ & $ \rootcoord{0, 0, 1, 2, 1, 2, 0, 1} $ & $ \rootcoord{2, 2, 2, 2, 1, 2, 1, 1} $ \\
		$ \rootcoord{\gold, \gold^2, \gold^2, \gold^2} $ & $ \rootcoord{1, 1, 1, 1, 1, 1, 0, 1} $ & $ \rootcoord{1, 1, 2, 2, 1, 2, 1, 1} $ \\
		\bottomrule
	\end{tabular}
	\caption{The preimages of the positive roots in $ H_4 $ under $ \foldproj $, part 1/2. See~\ref{rem:H3-preimages}.}
	\label{fig:foldrootse8-1}
\end{figure}

\begin{figure}[htb]
	\centering\begin{tabular}{ccc}
		\toprule
		$ \beta $ & $ \alpha_1 $ & $ \alpha_2 $ \\
		\midrule
		$ \rootcoord{1, \gold^2, \gold^2, \gold^2} $ & $ \rootcoord{1, 1, 1, 1, 1, 1, 1, 0} $ & $ \rootcoord{1, 1, 2, 2, 1, 2, 0, 1} $ \\
		$ \rootcoord{\gold, 2\gold, 2\gold+1, \gold^2} $ & $ \rootcoord{0, 1, 1, 2, 1, 2, 0, 1} $ & $ \rootcoord{2, 2, 2, 3, 1, 2, 1, 1} $ \\
		$ \rootcoord{\gold, \gold^2, 2\gold+1, \gold^2} $ & $ \rootcoord{1, 1, 1, 2, 1, 1, 0, 1} $ & $ \rootcoord{1, 2, 2, 3, 1, 2, 1, 1} $ \\
		$ \rootcoord{1, \gold^2, 2\gold+1, \gold^2} $ & $ \rootcoord{1, 1, 1, 2, 1, 1, 1, 0} $ & $ \rootcoord{1, 2, 2, 3, 1, 2, 0, 1} $ \\
		$ \rootcoord{\gold, 2\gold, 2\gold+1, 2\gold+1} $ & $ \rootcoord{0, 1, 2, 2, 1, 2, 0, 1} $ & $ \rootcoord{2, 2, 3, 3, 2, 2, 1, 1} $ \\
		$ \rootcoord{\gold, \gold^2, 2\gold+1, 2\gold+1} $ & $ \rootcoord{1, 1, 2, 2, 1, 1, 0, 1} $ & $ \rootcoord{1, 2, 3, 3, 2, 2, 1, 1} $ \\
		$ \rootcoord{1, \gold^2, 2\gold+1, 2\gold+1} $ & $ \rootcoord{1, 1, 2, 2, 1, 1, 1, 0} $ & $ \rootcoord{1, 2, 3, 3, 2, 2, 0, 1} $ \\
		$ \rootcoord{\gold, 2\gold, 3\gold+1, 2\gold+1} $ & $ \rootcoord{0, 1, 2, 3, 1, 2, 0, 1} $ & $ \rootcoord{2, 3, 3, 4, 2, 2, 1, 1} $ \\
		$ \rootcoord{\gold, \gold^2, 2\gold+2, 2\gold+1} $ & $ \rootcoord{1, 2, 2, 2, 1, 1, 0, 1} $ & $ \rootcoord{1, 2, 3, 4, 2, 2, 1, 1} $ \\
		$ \rootcoord{1, \gold^2, 2\gold+2, 2\gold+1} $ & $ \rootcoord{1, 2, 2, 2, 1, 1, 1, 0} $ & $ \rootcoord{1, 2, 3, 4, 2, 2, 0, 1} $ \\
		$ \rootcoord{\gold, 2\gold, 3\gold+1, 2\gold+2} $ & $ \rootcoord{0, 1, 2, 3, 2, 2, 0, 1} $ & $ \rootcoord{2, 3, 4, 4, 2, 2, 1, 1} $ \\
		$ \rootcoord{\gold, 2\gold+1, 3\gold+1, 2\gold+1} $ & $ \rootcoord{1, 1, 2, 3, 1, 2, 0, 1} $ & $ \rootcoord{2, 3, 3, 4, 2, 3, 1, 1} $ \\
		$ \rootcoord{\gold, 2\gold+1, 2\gold+2, 2\gold+1} $ & $ \rootcoord{1, 2, 2, 2, 1, 2, 0, 1} $ & $ \rootcoord{2, 2, 3, 4, 2, 3, 1, 1} $ \\
		$ \rootcoord{1, \gold+2, 2\gold+2, 2\gold+1} $ & $ \rootcoord{2, 2, 2, 2, 1, 1, 1, 0} $ & $ \rootcoord{1, 2, 3, 4, 2, 3, 0, 1} $ \\
		$ \rootcoord{\gold, 2\gold+1, 3\gold+1, 2\gold+2} $ & $ \rootcoord{1, 1, 2, 3, 2, 2, 0, 1} $ & $ \rootcoord{2, 3, 4, 4, 2, 3, 1, 1} $ \\
		$ \rootcoord{\gold^2, 2\gold+1, 3\gold+1, 2\gold+1} $ & $ \rootcoord{1, 1, 2, 3, 1, 2, 1, 1} $ & $ \rootcoord{2, 3, 3, 4, 2, 3, 1, 2} $ \\
		$ \rootcoord{\gold^2, 2\gold+1, 2\gold+2, 2\gold+1} $ & $ \rootcoord{1, 2, 2, 2, 1, 2, 1, 1} $ & $ \rootcoord{2, 2, 3, 4, 2, 3, 1, 2} $ \\
		$ \rootcoord{\gold^2, \gold+2, 2\gold+2, 2\gold+1} $ & $ \rootcoord{2, 2, 2, 2, 1, 1, 1, 1} $ & $ \rootcoord{1, 2, 3, 4, 2, 3, 1, 2} $ \\
		$ \rootcoord{\gold, 2\gold+1, 3\gold+2, 2\gold+2} $ & $ \rootcoord{1, 2, 2, 3, 2, 2, 0, 1} $ & $ \rootcoord{2, 3, 4, 5, 2, 3, 1, 1} $ \\
		$ \rootcoord{\gold^2, 2\gold+1, 3\gold+1, 2\gold+2} $ & $ \rootcoord{1, 1, 2, 3, 2, 2, 1, 1} $ & $ \rootcoord{2, 3, 4, 4, 2, 3, 1, 2} $ \\
		$ \rootcoord{\gold, 2\gold+1, 3\gold+2, 3\gold+1} $ & $ \rootcoord{1, 2, 3, 3, 1, 2, 0, 1} $ & $ \rootcoord{2, 3, 4, 5, 3, 3, 1, 1} $ \\
		$ \rootcoord{\gold^2, 2\gold+1, 3\gold+2, 2\gold+2} $ & $ \rootcoord{1, 2, 2, 3, 2, 2, 1, 1} $ & $ \rootcoord{2, 3, 4, 5, 2, 3, 1, 2} $ \\
		$ \rootcoord{\gold^2, 2\gold+1, 3\gold+2, 3\gold+1} $ & $ \rootcoord{1, 2, 3, 3, 1, 2, 1, 1} $ & $ \rootcoord{2, 3, 4, 5, 3, 3, 1, 2} $ \\
		$ \rootcoord{\gold^2, 2\gold+2, 3\gold+2, 2\gold+2} $ & $ \rootcoord{2, 2, 2, 3, 2, 2, 1, 1} $ & $ \rootcoord{2, 3, 4, 5, 2, 4, 1, 2} $ \\
		$ \rootcoord{\gold^2, 2\gold+2, 3\gold+2, 3\gold+1} $ & $ \rootcoord{2, 2, 3, 3, 1, 2, 1, 1} $ & $ \rootcoord{2, 3, 4, 5, 3, 4, 1, 2} $ \\
		$ \rootcoord{\gold^2, 2\gold+2, 3\gold+3, 3\gold+1} $ & $ \rootcoord{2, 3, 3, 3, 1, 2, 1, 1} $ & $ \rootcoord{2, 3, 4, 6, 3, 4, 1, 2} $ \\
		$ \rootcoord{2\gold, 3\gold+1, 4\gold+2, 3\gold+2} $ & $ \rootcoord{1, 2, 3, 4, 2, 3, 0, 2} $ & $ \rootcoord{3, 4, 5, 6, 3, 4, 2, 2} $ \\
		$ \rootcoord{\gold^2, 3\gold+1, 4\gold+2, 3\gold+2} $ & $ \rootcoord{1, 2, 3, 4, 2, 3, 1, 1} $ & $ \rootcoord{3, 4, 5, 6, 3, 4, 1, 2} $ \\
		$ \rootcoord{\gold^2, 2\gold+2, 4\gold+2, 3\gold+2} $ & $ \rootcoord{2, 2, 3, 4, 2, 2, 1, 1} $ & $ \rootcoord{2, 4, 5, 6, 3, 4, 1, 2} $ \\
		$ \rootcoord{\gold^2, 2\gold+2, 3\gold+3, 3\gold+2} $ & $ \rootcoord{2, 3, 3, 3, 2, 2, 1, 1} $ & $ \rootcoord{2, 3, 5, 6, 3, 4, 1, 2} $ \\
		\bottomrule
	\end{tabular}
	\caption{The preimages of the positive roots in $ H_4 $ under $ \foldproj $, part 2/2. See~\ref{rem:H3-preimages}.}
	\label{fig:foldrootse8-2}
\end{figure}

\begin{figure}[htb]
	\centering
	\begin{tabular}{ccccc}
		\toprule
		$ \alpha $ & $ \inverparbr{\alpha}{\rho_0} $ & $\inverparbr{\alpha}{\rho_1}$ & $\inverparbr{\alpha}{\rho_2}$ & $\inverparbr{\alpha}{\rho_3}$ \\
		\midrule
		$ \rootcoord{0, 0, 0, 1} $ & $ (1, 1) $ & $ (1, 1) $ & $ (-1, -1) $ & $ (-1, -1) $ \\
		$ \rootcoord{0, 0, 1, 0} $ & $ (1, 1) $ & $ (-1, -1) $ & $ (-1, -1) $ & $ (1, -1) $ \\
		$ \rootcoord{0, 1, 0, 0} $ & $ (-1, -1) $ & $ (-1, -1) $ & $ (1, 1) $ & $ (1, 1) $ \\
		$ \rootcoord{1, 0, 0, 0} $ & $ (-1, -1) $ & $ (1, 1) $ & $ (1, 1) $ & $ (1, 1) $ \\
		$ \rootcoord{0, 0, \gold, 1} $ & $ (1, 1) $ & $ (1, 1) $ & $ (1, -1) $ & $ (1, 1) $ \\
		$ \rootcoord{0, 0, \gold, \gold} $ & $ (1, 1) $ & $ (-1, -1) $ & $ (1, 1) $ & $ (1, -1) $ \\
		$ \rootcoord{0, 0, 1, \gold} $ & $ (1, 1) $ & $ (-1, 1) $ & $ (1, -1) $ & $ (-1, -1) $ \\
		$ \rootcoord{0, 1, 1, 0} $ & $ (-1, 1) $ & $ (1, 1) $ & $ (-1, -1) $ & $ (1, 1) $ \\
		$ \rootcoord{1, 1, 0, 0} $ & $ (1, 1) $ & $ (-1, -1) $ & $ (1, -1) $ & $ (1, 1) $ \\
		$ \rootcoord{0, \gold, \gold, 1} $ & $ (-1, -1) $ & $ (-1, 1) $ & $ (1, 1) $ & $ (-1, -1) $ \\
		$ \rootcoord{0, \gold, \gold, \gold} $ & $ (1, 1) $ & $ (1, -1) $ & $ (-1, -1) $ & $ (-1, 1) $ \\
		$ \rootcoord{0, 1, 1, \gold} $ & $ (-1, -1) $ & $ (1, -1) $ & $ (1, 1) $ & $ (-1, 1) $ \\
		$ \rootcoord{1, 1, 1, 0} $ & $ (1, -1) $ & $ (1, 1) $ & $ (-1, 1) $ & $ (1, -1) $ \\
		$ \rootcoord{\gold, \gold, \gold, 1} $ & $ (1, -1) $ & $ (1, 1) $ & $ (1, 1) $ & $ (1, 1) $ \\
		$ \rootcoord{\gold, \gold, \gold, \gold} $ & $ (-1, 1) $ & $ (1, 1) $ & $ (-1, -1) $ & $ (1, -1) $ \\
		$ \rootcoord{0, \gold, \gold^2, \gold} $ & $ (1, 1) $ & $ (-1, 1) $ & $ (1, 1) $ & $ (-1, -1) $ \\
		$ \rootcoord{0, 1, \gold^2, \gold} $ & $ (-1, -1) $ & $ (-1, -1) $ & $ (-1, 1) $ & $ (1, -1) $ \\
		$ \rootcoord{1, 1, 1, \gold} $ & $ (1, 1) $ & $ (1, 1) $ & $ (1, 1) $ & $ (-1, -1) $ \\
		$ \rootcoord{\gold, \gold, \gold^2, \gold} $ & $ (-1, 1) $ & $ (-1, -1) $ & $ (1, 1) $ & $ (1, -1) $ \\
		$ \rootcoord{0, \gold, 2\gold, \gold^2} $ & $ (-1, 1) $ & $ (1, 1) $ & $ (-1, 1) $ & $ (1, 1) $ \\
		$ \rootcoord{0, \gold, \gold^2, \gold^2} $ & $ (-1, 1) $ & $ (1, 1) $ & $ (-1, -1) $ & $ (1, 1) $ \\
		$ \rootcoord{0, 1, \gold^2, \gold^2} $ & $ (-1, -1) $ & $ (1, -1) $ & $ (1, 1) $ & $ (-1, 1) $ \\
		$ \rootcoord{1, 1, \gold^2, \gold} $ & $ (1, 1) $ & $ (-1, 1) $ & $ (-1, 1) $ & $ (1, -1) $ \\
		$ \rootcoord{\gold, \gold, 2\gold, \gold^2} $ & $ (1, 1) $ & $ (-1, 1) $ & $ (-1, 1) $ & $ (1, 1) $ \\
		$ \rootcoord{\gold, \gold, \gold^2, \gold^2} $ & $ (1, 1) $ & $ (-1, -1) $ & $ (-1, -1) $ & $ (-1, 1) $ \\
		$ \rootcoord{\gold, \gold^2, \gold^2, \gold} $ & $ (1, 1) $ & $ (1, 1) $ & $ (1, 1) $ & $ (1, -1) $ \\
		$ \rootcoord{1, 1, \gold^2, \gold^2} $ & $ (1, 1) $ & $ (-1, 1) $ & $ (1, 1) $ & $ (-1, 1) $ \\
		$ \rootcoord{1, \gold^2, \gold^2, \gold} $ & $ (1, -1) $ & $ (1, 1) $ & $ (1, 1) $ & $ (1, -1) $ \\
		$ \rootcoord{\gold, 2\gold, 2\gold, \gold^2} $ & $ (1, 1) $ & $ (1, 1) $ & $ (-1, -1) $ & $ (1, 1) $ \\
		$ \rootcoord{\gold, \gold^2, \gold^2, \gold^2} $ & $ (1, 1) $ & $ (1, 1) $ & $ (-1, 1) $ & $ (-1, 1) $ \\
		\bottomrule
	\end{tabular}
	\caption{Values of the parity map $ \map{\inverparsym}{H_4 \times \Set{\rho_0, \rho_1, \rho_2, \rho_3}}{\Set{\pm 1} \times \Set{\pm 1}}{}{} $ in~\ref{rem:fold-ex-parmap}, part 1/2. The values $ \inverpar{\alpha}{\delta} $ for negative roots $ \alpha $ are given by the formula $ \inverpar{-\alpha}{\delta} = \inverpar{\alpha}{\delta} $.}
	\label{fig:parmapex-H4-1}
\end{figure}

\begin{figure}[htb]
	\centering
	\begin{tabular}{ccccc}
		\toprule
		$ \alpha $ & $ \inverparbr{\alpha}{\rho_0} $ & $\inverparbr{\alpha}{\rho_1}$ & $\inverparbr{\alpha}{\rho_2}$ & $\inverparbr{\alpha}{\rho_3}$ \\
		\midrule
		$ \rootcoord{1, \gold^2, \gold^2, \gold^2} $ & $ (1, -1) $ & $ (1, 1) $ & $ (-1, 1) $ & $ (-1, 1) $ \\
		$ \rootcoord{\gold, 2\gold, 2\gold+1, \gold^2} $ & $ (1, 1) $ & $ (-1, 1) $ & $ (1, 1) $ & $ (-1, 1) $ \\
		$ \rootcoord{\gold, \gold^2, 2\gold+1, \gold^2} $ & $ (1, 1) $ & $ (-1, -1) $ & $ (1, 1) $ & $ (-1, 1) $ \\
		$ \rootcoord{1, \gold^2, 2\gold+1, \gold^2} $ & $ (1, -1) $ & $ (1, 1) $ & $ (1, 1) $ & $ (-1, 1) $ \\
		$ \rootcoord{\gold, 2\gold, 2\gold+1, 2\gold+1} $ & $ (1, 1) $ & $ (-1, 1) $ & $ (-1, 1) $ & $ (1, 1) $ \\
		$ \rootcoord{\gold, \gold^2, 2\gold+1, 2\gold+1} $ & $ (1, 1) $ & $ (-1, -1) $ & $ (-1, -1) $ & $ (1, 1) $ \\
		$ \rootcoord{1, \gold^2, 2\gold+1, 2\gold+1} $ & $ (1, -1) $ & $ (1, 1) $ & $ (-1, -1) $ & $ (1, 1) $ \\
		$ \rootcoord{\gold, 2\gold, 3\gold+1, 2\gold+1} $ & $ (1, 1) $ & $ (-1, 1) $ & $ (1, 1) $ & $ (-1, -1) $ \\
		$ \rootcoord{\gold, \gold^2, 2\gold+2, 2\gold+1} $ & $ (1, 1) $ & $ (-1, -1) $ & $ (1, 1) $ & $ (1, 1) $ \\
		$ \rootcoord{1, \gold^2, 2\gold+2, 2\gold+1} $ & $ (1, -1) $ & $ (-1, 1) $ & $ (1, 1) $ & $ (1, 1) $ \\
		$ \rootcoord{\gold, 2\gold, 3\gold+1, 2\gold+2} $ & $ (1, 1) $ & $ (-1, 1) $ & $ (1, 1) $ & $ (1, 1) $ \\
		$ \rootcoord{\gold, 2\gold+1, 3\gold+1, 2\gold+1} $ & $ (-1, -1) $ & $ (1, -1) $ & $ (-1, 1) $ & $ (-1, -1) $ \\
		$ \rootcoord{\gold, 2\gold+1, 2\gold+2, 2\gold+1} $ & $ (-1, -1) $ & $ (1, -1) $ & $ (-1, -1) $ & $ (1, 1) $ \\
		$ \rootcoord{1, \gold+2, 2\gold+2, 2\gold+1} $ & $ (1, 1) $ & $ (1, -1) $ & $ (1, 1) $ & $ (1, 1) $ \\
		$ \rootcoord{\gold, 2\gold+1, 3\gold+1, 2\gold+2} $ & $ (-1, -1) $ & $ (1, -1) $ & $ (-1, -1) $ & $ (1, 1) $ \\
		$ \rootcoord{\gold^2, 2\gold+1, 3\gold+1, 2\gold+1} $ & $ (1, 1) $ & $ (1, 1) $ & $ (-1, 1) $ & $ (-1, -1) $ \\
		$ \rootcoord{\gold^2, 2\gold+1, 2\gold+2, 2\gold+1} $ & $ (1, 1) $ & $ (-1, 1) $ & $ (-1, -1) $ & $ (1, 1) $ \\
		$ \rootcoord{\gold^2, \gold+2, 2\gold+2, 2\gold+1} $ & $ (-1, 1) $ & $ (-1, -1) $ & $ (1, 1) $ & $ (1, 1) $ \\
		$ \rootcoord{\gold, 2\gold+1, 3\gold+2, 2\gold+2} $ & $ (-1, -1) $ & $ (1, 1) $ & $ (1, 1) $ & $ (-1, 1) $ \\
		$ \rootcoord{\gold^2, 2\gold+1, 3\gold+1, 2\gold+2} $ & $ (1, 1) $ & $ (1, 1) $ & $ (-1, -1) $ & $ (1, 1) $ \\
		$ \rootcoord{\gold, 2\gold+1, 3\gold+2, 3\gold+1} $ & $ (-1, -1) $ & $ (1, 1) $ & $ (1, 1) $ & $ (-1, -1) $ \\
		$ \rootcoord{\gold^2, 2\gold+1, 3\gold+2, 2\gold+2} $ & $ (1, 1) $ & $ (-1, 1) $ & $ (1, 1) $ & $ (-1, 1) $ \\
		$ \rootcoord{\gold^2, 2\gold+1, 3\gold+2, 3\gold+1} $ & $ (1, 1) $ & $ (-1, 1) $ & $ (1, 1) $ & $ (-1, -1) $ \\
		$ \rootcoord{\gold^2, 2\gold+2, 3\gold+2, 2\gold+2} $ & $ (1, 1) $ & $ (1, -1) $ & $ (1, 1) $ & $ (-1, 1) $ \\
		$ \rootcoord{\gold^2, 2\gold+2, 3\gold+2, 3\gold+1} $ & $ (1, 1) $ & $ (1, -1) $ & $ (-1, 1) $ & $ (-1, -1) $ \\
		$ \rootcoord{\gold^2, 2\gold+2, 3\gold+3, 3\gold+1} $ & $ (1, 1) $ & $ (1, 1) $ & $ (1, -1) $ & $ (-1, 1) $ \\
		$ \rootcoord{2\gold, 3\gold+1, 4\gold+2, 3\gold+2} $ & $ (-1, -1) $ & $ (1, 1) $ & $ (1, 1) $ & $ (1, 1) $ \\
		$ \rootcoord{\gold^2, 3\gold+1, 4\gold+2, 3\gold+2} $ & $ (-1, 1) $ & $ (1, -1) $ & $ (1, 1) $ & $ (1, 1) $ \\
		$ \rootcoord{\gold^2, 2\gold+2, 4\gold+2, 3\gold+2} $ & $ (1, 1) $ & $ (1, 1) $ & $ (-1, -1) $ & $ (1, 1) $ \\
		$ \rootcoord{\gold^2, 2\gold+2, 3\gold+3, 3\gold+2} $ & $ (1, 1) $ & $ (1, 1) $ & $ (-1, 1) $ & $ (1, -1) $ \\
		\bottomrule
	\end{tabular}
	\caption{Values of the parity map $ \map{\inverparsym}{H_4 \times \Set{\rho_0, \rho_1, \rho_2, \rho_3}}{\Set{\pm 1} \times \Set{\pm 1}}{}{} $ in~\ref{rem:fold-ex-parmap}, part 2/2. The values $ \inverpar{\alpha}{\delta} $ for negative roots $ \alpha $ are given by the formula $ \inverpar{-\alpha}{\delta} = \inverpar{\alpha}{\delta} $.}
	\label{fig:parmapex-H4-2}
\end{figure}
\postmidfigure

\clearpage

\bibliographystyle{amsalpha}
\bibliography{references}
\addcontentsline{toc}{section}{\refname}

\end{document}